\documentclass[11pt,a4paper]{article}
\usepackage[utf8]{inputenc}
\usepackage[T1]{fontenc}
\usepackage{amsmath, amssymb, amsthm}
\numberwithin{equation}{section}

\usepackage{graphicx}
\usepackage{color}
\usepackage{tikz}
\usepackage{subcaption}
\usepackage{geometry}
\usepackage{booktabs}

\usepackage{tikz}
\usetikzlibrary{arrows.meta, positioning, shapes.misc}

\usepackage{authblk}

\usepackage[hidelinks]{hyperref} 
\hypersetup{
    colorlinks=true,    
    linkcolor=blue,     
    citecolor=green,    
    urlcolor=green,      
    pdftitle={Flux-Based Time Bounds for Harmonically Degenerate Quantum Systems},
    pdfauthor={Thomas Chambrion},
    pdfsubject={We derive easily computable estimates, from above and from below, for the time needed to steer a single input bilinear quantum system from a source to a target. },
    pdfkeywords={QSL, Quantum control, bilinear control, optimal control},
    pdfproducer={LaTeX with hyperref}, 
    pdfcreator={pdfTeX}                
}
\usepackage{orcidlink}

\title{Flux-Based Time Bounds for Harmonically Degenerate Quantum Systems
}

\author{Thomas Chambrion \orcidlink{0000-0002-9006-7596}}

\affil{Institut de Recherche Mathématique Avancée
UMR 7501\\ Université de Strasbourg et CNRS \\
7 rue René Descartes,
67000 Strasbourg, France\\
\url{thomas.chambrion@math.unistra.fr}}

\date{\today}

\newtheorem{theorem}{Theorem}
\newtheorem{definition}{Definition}
\newtheorem{problem}{Problem}
\newtheorem{proposition}[theorem]{Proposition}
\newtheorem{remark}{Remark}
\newtheorem{corollary}[theorem]{Corollary}
\newtheorem{assumption}{Assumption}
\newtheorem{lemma}[theorem]{Lemma}

\newcommand{\C}{\ensuremath{\mathbb{C}}}
\newcommand{\R}{\ensuremath{\mathbb{R}}}
\newcommand{\N}{\ensuremath{\mathbb{N}}}

\begin{document}

\maketitle

\begin{abstract}
We study transfer-time estimates for finite-dimensional bilinear quantum control systems driven by an 
unbounded scalar control, in both the state-to-state setting (complex unit sphere) and the unitary 
synthesis setting (special unitary group). Lower bounds are recalled and refined via spectral and chordal 
methods. The main contribution is constructive upper bounds: we associate with the system a coupling graph 
whose edges correspond to non-resonant drift couplings, and construct transfer protocols along spanning 
trees. In contrast with earlier results requiring full rational independence of the spectrum of the control 
operator, our estimates allow rational resonances; they incorporate, via explicit correction factors 
derived from the arithmetic of the spectral gaps, the efficiency of the averaging procedure and of the 
phase adjustments. As an application, we consider finite-dimensional truncations of the planar quantum 
rotator under a dipolar electric field, and show that the transfer time between any two modes tends 
to zero as the truncation dimension grows.
\end {abstract}

\newpage

\setcounter{tocdepth}{2}
\tableofcontents

\newpage

\section{Introduction}
\subsection{Transfer time for bilinear quantum control systems}\label{SEC_transfer_time}
In the following, $d$ is an integer and the Lie group $G$ is either $ U(d)$ or $SU(d)$.
The Lie algebra of $G$ is denoted $\mathfrak{g}$.
 
We consider a finite-dimensional quantum system whose state $x(t)$ evolves, depending on the purpose, in the Lie group $G$ (unitary synthesis) or on the complex unit sphere $\mathbb{S}^{2d-1} \subset \mathbb{C}^d$ (state-to-state transfer). The dynamics are governed by the bilinear Schrödinger equation:
\begin{equation}\label{EQ_main}
\dot{x}(t) = \left( A + u(t)B \right) x(t), \quad x(0) = I_d \text{ or } \psi_0 \in \mathbb{S}^{2d-1}
\end{equation}
where $I_d$ is the identity matrix of order $d$, $\psi_0$ is an initial condition on the 
complex unit sphere, $A, B \in \mathfrak{g}$ are skew-adjoint matrices, and $u: [0, T] \to \mathbb{R}$ is a locally integrable control function ($u \in L^1_{loc}$). The propagator between times $s$ and $t$ associated with \eqref{EQ_main} is denoted $t\mapsto X^u_{t,s} \in G$, and the solution of \eqref{EQ_main} reads $x(t)= X^u_{t,0} x(0)$. The canonical inner product of the Hilbert space $\mathbb{C}^d$ is denoted by $\langle \cdot , \cdot \rangle$, and the associated norm by $\| \cdot \|$.

\begin{remark}\label{REM_decalage_diagonal_trace_nulle}
Physically, a \emph{global} phase shift on $x$, namely $x \rightsquigarrow e^{i \theta} 
x$ with $\theta \in \mathbb{R}$, is irrelevant. Consequently, shifting $A$ to $A- 
\frac{1}{d} \mathrm{Tr}(A) I_d$ and $B$ to $B- \frac{1}{d}\mathrm{Tr}(B) I_d$ (which 
maps $x(t)$ to $\exp \left(-\frac{\mathrm{Tr}(A)}{d} t - \frac{\mathrm{Tr}(B)}{d} 
\int_0^t u(s) ds \right) x(t)$) does not alter the physical state transition. We may 
thus assume without loss of generality that $\mathrm{Tr}(A)=\mathrm{Tr}(B)=0$, meaning 
$A,B \in \mathfrak{su}(d)$. This allows us to work in $G=SU(d)$, the group of special 
unitary matrices, which is semi-simple. While this global phase shift is 
standard in the literature, we will generally avoid it unless explicitly stated 
otherwise. In particular, the estimates given in Section \ref{SEC_Upper_bound} are 
valid on $SU(d)$ and on $U(d)$.
\end{remark}

We work in an eigenbasis $(\phi_1^B, \ldots, \phi_d^B)$ of $B$. Thus, $B$ is represented by the diagonal matrix:
\begin{equation}
B = \mathrm{diag}(i b_1, i b_2, \dots, i b_d), \quad b_j \in \mathbb{R}.
\end{equation}
The drift matrix $A$ is decomposed as $A = A_D + A_{ND}$, where $A_D = \mathrm{diag}(A)$ is the diagonal part of $A$ and $A_{ND}$, with entries $\langle \phi_j^B , A \phi_k^B \rangle$ for $j\neq k$, is the hollow part (zero diagonal) representing the couplings between the eigenstates of $B$.

The objective is to estimate the infimum time $\tau$ required to reach a target state $x_f$ from an initial state $x_0$.  We distinguish two problems depending on the state space.

\begin{problem}[Exact State-to-State Transfer (SS)]\label{PROBLEM_SS}
Given $\psi_0, \psi_f \in \mathbb{S}^{2d-1}$, find:
\begin{equation}
\tau_{\mathrm{SS}}(\psi_0,\psi_f) = \inf \{ T > 0 \mid \exists u \in L^1_{loc},  X^u_{T,0} \psi_0 = \psi_f \}.
\end{equation}
We denote by $\tau^{\max}_{\mathrm{SS}}$ the worst value, $\tau^{\max}_{\mathrm{SS}}=\sup_{\psi_0,\psi_f \in \mathbb{S}^{2d-1}}\tau_{\mathrm{SS}}(\psi_0,\psi_f) $.
\end{problem}

\begin{problem}[Exact Unitary Operator Synthesis (US)]\label{PROBLEM_US}
Given $X_f \in G$, find:
\begin{equation}
\tau_{\mathrm{US}}(X_f) = \inf \{ T > 0 \mid \exists u \in L^1_{loc},  X^u_{T,0} = X_f \}.
\end{equation}
We denote by $\tau^{\max}_{\mathrm{US}}\! =\! \sup_{X_f \in G} \tau_{\mathrm{US}}(X_f)$ the worst-case value, called the time diameter of the group.
\end{problem}

In what follows, $\tau_{\mathrm{SS}}$ and $\tau_{\mathrm{US}}$ are infima rather than minima: in general the optimal 
control does not exist and  minimizing sequences of controls may blow up in  $L^\infty$ 
(Remark \ref{REM_Blow-up}). 
It is well established that $\tau^{\max}_{\mathrm{US}}$ and $\tau^{\max}_{\mathrm{SS}}$ 
are finite under the generic condition that the Lie algebra $\mathrm{Lie}(A,B)$ 
generated by $A$ and $B$ contains $\mathfrak{g}$ (historical proofs are provided in 
\cite{JURDJEVIC1972313} for $\tau^{\max}_{\mathrm{US}}$, \cite{doi:10.1137/0125025} for 
$\tau^{\max}_{\mathrm{SS}}$; see also \cite{Sachkov} and \cite{ALBERTINI2002213}).

In this work, we operate under the following generic assumption:
\begin{assumption}\label{ASS_control}
The matrix $B$ has simple eigenvalues, and $\mathrm{Lie}(A,B)$ contains $\mathfrak{g}$.
\end{assumption}

\subsection{Main results}

We bracket the optimal transfer time between two computable thresholds. 
First, we adapt classical spectral and energy-based methods to our specific algebraic 
setting to derive lower bounds on the controllability time; while these rely on established 
theoretical frameworks, they delimit the forbidden region for any control protocol. 
Second, our main contribution, we propose explicit, 
ansatz-driven control strategies. These yield constructive upper bounds relying on the topological and 
Diophantine properties of the system. 

We illustrate our approach by computing these time estimates for finite-dimensional 
compressions of the planar rotator under a dipolar electric drive.

\subsubsection{Lower Bounds (Infimum Time estimates)}

A family of lower bound estimates is obtained by rewriting the system \eqref{EQ_main} in 
the quotient space $G/\mathcal{T}_B$, where $\mathcal{T}_B$ is the set of matrices of 
$G$ commuting with $B$. This follows established geometric approaches 
\cite{PhysRevA.65.032301, zbMATH05117619}. Propositions \ref{PRO_Spectral_Bound} (see 
\cite{khaneja2005optimal} for a more general version) and \ref{PRO_borne_Agrachev} below 
are classical and are provided for completeness.

\begin{proposition}[Global Spectral Bound for State-to-State Transfer]\label{PRO_Spectral_Bound}
Let $\mathcal{T}_B$ be  the set of diagonal matrices of $G$ in 
the eigenbasis of $B$. The minimal transfer time $\tau_{\mathrm{SS}}$ for a pure state transition $\psi_0 
\to \psi_f$  is bounded from below by the ratio of the 
minimal Fubini-Study distance between the initial state and the $\mathcal{T}_B$-orbit of the target, to the 
operator norm of the non-diagonal drift $A_{ND}$:
\begin{equation}
\tau_{\mathrm{SS}}(\psi_0,\psi_f) \ge \frac{\inf_{D \in \mathcal{T}_B} d_{FS}(\psi_0, D\psi_f)}{\|A_{ND}\|_{op}},
\end{equation}
where $d_{FS}$ stands for the Fubini-Study distance $d_{FS}(\psi_0, D\psi_f) = \arccos\left( \left| \langle \psi_0, D\psi_f \rangle \right| \right)$ on $\mathbb{CP}^{d-1}$.
\end{proposition}
A proof of Proposition \ref{PRO_Spectral_Bound} is provided, for reference, in Section \ref{SEC_Spectral_bounds}.

\begin{proposition}[Global Spectral Bound for Unitary Synthesis, Corollary 2.12 in \cite{zbMATH05117619}]\label{PRO_borne_Agrachev}
Recall that we work in the eigenbasis of $B$. If $G=SU(d)$ (meaning that $A$ and $B$ have trace zero), then:
\begin{equation}
\tau^{\max}_{\mathrm{US}} \ge \frac{\pi \sqrt{d^2 - 1 }}{\sqrt{3 d \sum_{j\neq k} |A_{jk}|^2}}.
\end{equation}
\end{proposition}

While the bounds of Proposition \ref{PRO_borne_Agrachev} is uniform (of worst-case type), the matrix analogue of the spectral bound of Proposition \ref{PRO_Spectral_Bound} (the Lipschitz argument applied to the propagator) can be specialized for the synthesis of a specific propagator. Proposition \ref{PRO_Chordal_relaxation} below 
provides the corresponding chordal relaxation.

\begin{proposition}[Analytic Chordal Relaxation]\label{PRO_Chordal_relaxation}
For unitary synthesis with target $X_f \in G$ (expressed in the eigenbasis of $B$), the minimal transfer time $\tau_{\mathrm{US}}(X_f)$ is bounded from below by:
\begin{equation}\label{EQ_chrodal}
\tau_{\mathrm{US}} (X_f)\ge \frac{\sqrt{2 \left( d - \sum_{j=1}^d \left| (X_f)_{jj} \right| \right)}}{\sqrt{\sum_{j\neq k} |A_{jk}|^2}},
\end{equation}
where $(X_f)_{jj}$ are the diagonal elements of the target unitary operator.
\end{proposition}
A proof of Proposition \ref{PRO_Chordal_relaxation} is given in Section \ref{SEC_Chordal}.

In parallel with geometric estimates, energy-like bounds are routinely employed in quantum speed limit analysis \cite{PhysRevLett.96.060503, CanevaPRL2009}. These proofs hinge on Ehrenfest's theorem to bound the speed of dynamical evolution, a strategy analogous to standard energy methods for infinite-dimensional PDEs \cite{Kato1953, BOUSSAID2020108412}. Specializing this principle to the commutant of our control generator yields the specific dispersion bound in Proposition \ref{PRO_minoration_Ehrenfest}.

\begin{proposition}[Momentum Dispersion Bound]\label{PRO_minoration_Ehrenfest}
For any Hermitian observable $O$ belonging to the commutant of the control generator ($[O, B] = 0$), and for any initial state $\psi_0$ and target state $\psi_f$ in the unit sphere of $\mathbb{C}^d$, 
\begin{equation}
\tau_{\mathrm{SS}}(\psi_0,\psi_f) \ge \frac{\left| \langle \psi_f, O \psi_f \rangle - \langle \psi_0, O \psi_0 \rangle \right|}{\| [O, A_{ND}] \|_{op}},
\end{equation}
where $\| \cdot \|_{op}$ denotes the spectral operator norm.
\end{proposition}

A proof of Proposition \ref{PRO_minoration_Ehrenfest} is given in Section \ref{SEC_Ehrenfest}.

\subsubsection{Upper Bounds (Constructive Estimates)}

The following definitions and propositions formulate our main results,
which give constructive upper bounds on the transfer time.

\begin{definition}[Non-Resonant Transition]
A transition $(j,k)$ with $A_{jk} \neq 0$ is \emph{non-resonant} if its associated spectral gap $\omega_{jk} = b_k - b_j$ is distinct from the spectral gaps of the other coupled transitions. Specifically, for any other pair of indices $(l,m) \notin \{(j,k), (k,j)\}$ such that $A_{lm} \neq 0$, the condition $\omega_{lm} \neq \pm \omega_{jk}$ holds.
\end{definition}

\begin{definition}[Coupling Graph]
The coupling graph associated with the non-diagonal drift $A_{ND}$ is the graph whose vertices are the indices $\{1,\dots,d\}$ of the eigenstates of $B$, and the edges are the non-resonant transitions $(j,k)$.
\end{definition}

\begin{assumption}\label{ASS_coupling_graph}
In what follows, we make the generic assumption that the coupling graph associated with the non-diagonal drift $A_{ND}$ is connected.
\end{assumption}

Propositions \ref{PRO_Result_Constructive_State_to_State} and \ref{PRO_Result_Explicit_Propagator} below are proved in Sections \ref{SEC_Upper_bound_StS} and \ref{SEC_Contructive_bounds_Propagators}.

\begin{proposition}[Constructive Upper Bound for State-to-State Transfer]\label{PRO_Result_Constructive_State_to_State} Assume Assumptions \ref{ASS_control} and \ref{ASS_coupling_graph} hold.
Let $\mathbf{T}$ be a spanning tree (i.e., a minimal connected subgraph containing all vertices, without loops) of the coupling graph associated with the non-diagonal drift $A_{ND}$. For any initial state $\psi_0$ and target state $\psi_f$ in the unit sphere $\mathbb{S}^{2d-1}$ of $\C^d$, the state-to-state transfer time  $\tau_{\mathrm{SS}}(\psi_0,\psi_f)$ is bounded from above by:
\begin{equation}
\tau_{\mathrm{SS}}^{\max} \le \eta \sum_{e=(j,k) \in \mathbf{T}} \frac{\alpha_{jk}}{|A_{jk}|},
\end{equation}
where $\alpha_{jk}$ is an averaging efficiency correction factor and $\eta$ is a topological phase-correction factor. These parameters depend on the Diophantine properties of the spectrum of the control operator $B$:
\begin{itemize}
\item $\eta = \frac{\pi}{2}$ if the spectral differences of $B$ have maximal rational rank (i.e., the set of real numbers $\{b_2-b_1,\ldots,b_d-b_{d-1}\}$ has $\mathbb{Q}$-rank $d-1$).
\item $\eta = \pi$ if the set of spectral differences of $B$ does not have maximal rational rank.
\item $\displaystyle \alpha_{jk}= \left| \cos \left( \frac{\pi}{M_{jk}+1} \right) \right|^{-1}\!\!\!\!$ where \[M_{jk}=\inf \left\{p \in \mathbb{N}^\ast \!\mid \exists (j',k') \notin \{(j,k), (k,j)\}, \frac{|b_{j'}-b_{k'}|}{|b_{j}-b_{k}|} = p  \mbox{ and } A_{j'k'}\neq 0\right\}.\] By convention, if the set is empty, $M_{jk}=+\infty$ and $\alpha_{jk}=1$. Whenever $(j,k)$ is non-resonant, $M_{jk} \geq 2$ and $\alpha_{jk} \leq 2$. In the generic case where the eigenvalues of $B$ have maximal rational rank, $M_{jk}=+\infty$ and $\alpha_{jk}=1$.
\end{itemize}
\end{proposition}

\begin{proposition}[Constructive Upper Bound for Unitary Synthesis]\label{PRO_Result_Explicit_Propagator}
Under Assumptions \ref{ASS_control} and \ref{ASS_coupling_graph}, let $\mathbf{T}$ be a spanning tree of the coupling graph associated with the non-diagonal drift $A_{ND}$, possessing $d - 1$ edges. For any target unitary propagator $X_f \in G$, the minimal synthesis time $\tau_{\mathrm{US}}(X_f)$ is bounded from above by:
\begin{equation}
\tau_{\mathrm{US}}^{\max} \le \pi \frac{d(d + 4)}{8} \max_{e \in T} \frac{\alpha_e}{|A_e|},
\end{equation}
where $\max_{e \in T} \frac{\alpha_e}{|A_e|}$ is the inverse weight of the slowest edge in the tree. The numbers $\alpha_e \leq 2$ are defined as in Proposition \ref{PRO_Result_Constructive_State_to_State}.
\end{proposition}

The case $\alpha_e\equiv 2$ corresponds to an extreme degeneracy scenario where, for all $(j,k)$ in the spanning tree $T$, there exist $j'$ and $k'$ such that $\left| \frac{b_{j'}-b_{k'}}{b_{j}-b_{k}} \right| = 2$. In the general case, sharper estimates depending on the Diophantine properties of the eigenvalues of $B$ are derived in Section \ref{SEC_Contructive_bounds_Propagators}.

As an application of these results, we consider the state-to-state time diameter
$\tau^{\max}_{\mathrm{SS}}(d)$ 
under the dynamics $\dot{x} = (A_o + uB_o)x$, where the drift $A_o$ 
and the control operator $B_o$ are defined by the $d \times d$ matrices:
\begin{equation}\label{EQ_dynamique_rotateur_impair}
A_o = -\mathrm{diag}(i k^2)_{k=1}^d \quad \text{and} \quad B_o = -\frac{i}{2} \begin{pmatrix}
0 & 1 & 0  & \cdots & 0 \\
1 & 0 & \ddots & \ddots & \vdots\\
0 & \ddots & \ddots & \ddots & 0 \\
\vdots & \ddots & \ddots & 0 & 1 \\
0 & \cdots & 0 & 1  & 0
\end{pmatrix}.
\end{equation}

To apply the constructive bounds derived in the previous section, the 
system must be rewritten in the eigenbasis of the control operator 
$B_o$. In this new representation, the transformed drift $A_o$ 
determines the edges of the coupling graph. However, computing the 
spectrum of $B_o$ reveals  resonances that
 violate 
Assumption \ref{ASS_coupling_graph}, so we adapt the routing strategy to obtain 
 the following bounds.

\begin{proposition}\label{PRO_rotateur_introduction}
For every integer $d \geq 4$, the  exact state-to-state time diameter $\tau_{\mathrm{SS}}^{\max}(d)$ of \eqref{EQ_dynamique_rotateur_impair} satisfies:
\begin{equation}
\frac{1}{8d} \leq  \tau^{\max}_{\mathrm{SS}}(d) \leq \frac{80}{d}.
\end{equation}
Consequently, $\tau^{\max}_{\mathrm{SS}}(d) = \mathcal{O}(1/d)$ as $d$ tends to $+\infty$.
\end{proposition}

The explicit change of basis and the spectral analysis of $B_o$ are given in Section \ref{sec:chiral_structure};  the 
complete proof of Proposition \ref{PRO_rotateur_introduction} is 
detailed in Section \ref{SEC_Chiral_StS}.

\subsection{State of the art and comparison with our results}

Finding the minimal time required to drive the system \eqref{EQ_main} from an initial 
condition to a given target is a central problem in quantum control, with applications 
to quantum computing and NMR. The fact that control time does not 
generally tend to zero, even with arbitrarily strong controls, led to the Quantum Speed 
Limit (QSL) concept.

\subsubsection{Geometric control theory}

For classical tools establishing the controllability of bilinear systems on Lie groups, see \cite{JURDJEVIC1972313}, the survey \cite{Sachkov}, or the physics-oriented review \cite{PhysRevA.92.062106}. Optimal control via the Pontryagin maximum principle \cite{Agrachev2004} provides analytical solutions for low-dimensional systems (e.g., \cite{zbMATH02122172}, \cite{zbMATH05154205}, \cite{zbMATH05001889}); for higher-dimensional systems, the PMP conditions become intractable.

\subsubsection{Lower bound of the control time}

The majority of the literature focuses on establishing lower bounds for the control time (QSL). An overview of QSLs for time-independent Hamiltonians is provided in \cite{FarmanianPRA2024}. For time-dependent Hamiltonians, analytical results for two-level systems are detailed in \cite{HegerfeldtPRL13}, comparing bang-bang controls with amplitude-unconstrained evolutions. Optimal control techniques have been applied in \cite{CanevaPRL2009} for state-to-state control, and in \cite{GoerzJPB11,GoerzNPJQI17,GuentherAVSQS2023,BasilewitschPRR2024} for unitary synthesis (quantum gate design).

Most lower bounds on the control time originate from a Lie-group
(geometric) perspective, initiated by \cite{khaneja2005optimal}.
For the single-input system \eqref{EQ_main} there are two complementary
factorizations; each pairs one of the two generators with a homogeneous
space and yields a different type of estimate.

\emph{Factoring out the drift.}
Setting $x(t)=e^{tA}Z(t)$, the remainder $Z$ satisfies
\begin{equation}
  \dot Z(t)=u(t)\,e^{-tA}B\,e^{tA}\,Z(t),\qquad Z(0)=I_d,
\end{equation}
hence is driven by the control alone. Since
$\|e^{-tA}B e^{tA}\|_{op}=\|B\|_{op}$ and $\|Z(t)\|_{op}=1$,
\begin{equation}
  \|Z(t)-I_d\|_{op}\le \|B\|_{op}\,\|u\|_{L^1[0,t]},
\end{equation}
so that for a target $X_f$ reached at time $t$,
\begin{equation}
  \|u\|_{L^1[0,t]}\ge \frac{\|e^{-tA}X_f-I_d\|_{op}}{\|B\|_{op}}.
\end{equation}
This is a lower bound on the $L^1$ norm of the control; if the amplitude
is also bounded, $|u|\le M$, it yields the time bound
$t\ge\|e^{-tA}X_f-I_d\|_{op}/(M\,\|B\|_{op})$, the physically relevant
regime since real devices impose a bound on the control amplitude. The
quantity $\|e^{-tA}X_f-I_d\|_{op}$ is the operator-norm distance between
the identity and the target after the drift flow $e^{tA}$ has been
factored out, i.e.\ a distance in the homogeneous space
$G/\overline{\exp(\mathbb{R}A)}$ with respect to the closed subgroup
(subtorus) generated by the drift $A$. This gives the bounds
parametrized by the control norm, e.g.\ based on
$\|u\|_{L^1}$~\cite{Naseri2024} and analogous results in the
$L^\infty$ setting.

\emph{Factoring out the control.}
Setting $x(t)=\exp(v(t)B)\,W(t)$ with $v(t)=\int_0^t u(s)\,ds$, the
remainder $W$ satisfies
\begin{equation}
  \dot W(t)=e^{-v(t)B}\,A\,e^{v(t)B}\,W(t),\qquad W(0)=I_d,
\end{equation}
hence is driven by the free drift alone. Since
$\|e^{-vB}A e^{vB}\|_{op}=\|A\|_{op}$,
\begin{equation}
  \|W(t)-I_d\|_{op}\le t\,\|A\|_{op},
\end{equation}
and for a target $X_f$ reached at time $t$,
\begin{equation}
  t\ge \frac{\|e^{-v(t)B}X_f-I_d\|_{op}}{\|A\|_{op}}
       \ge \inf_{v\in\mathbb{R}}\frac{\|e^{-vB}X_f-I_d\|_{op}}{\|A\|_{op}}.
\end{equation}
This is a lower bound directly on the transfer time, valid for arbitrary
(possibly unbounded) controls. The distance is measured in the
homogeneous space $G/\overline{\exp(\mathbb{R}B)}$ with respect to the
closed subgroup generated by the control $B$; this is the original method
of \cite{khaneja2005optimal}.

The two approaches are complementary: factoring out the control
(\cite{khaneja2005optimal}) yields a bound directly on the time and
remains valid without an a~priori amplitude bound, whereas factoring out
the drift yields a bound on the $L^1$ effort of the control, hence on the
time only under an amplitude constraint. The distances underlying these estimates are computed
on (generalized) flag manifolds; recent advances on the diameters of such
manifolds are discussed in \cite{zbMATH07308540} and
\cite{zbMATH07755548}.

\subsubsection{Upper bound for the control time}

Relatively few results provide upper bounds for control time when the control is not assumed 
bounded a priori. The main reference is \cite{zbMATH05117619}, which relies on approximations via 
Riemannian problems in the quotient space $SU(d)/\overline{\exp(\mathbb{R}B)}$: it requires the 
spectrum of $B$ to have maximal rational rank and all entries of $A$ to be non-zero, and its 
estimates are only useful when the entries of $A$ have comparable moduli, diverging rapidly 
when any entry approaches zero.

Our approach uses graph theory and flux analysis to obtain explicit upper bounds, valid even when most drift entries vanish or when the spectrum has strong degeneracies.
Graph-based routing has also been used for controllability in both 
mathematics \cite{zbMATH05509716} and physics \cite{PhysRevLett.73.58, Clements:16}.

\subsection{Content of the paper}

The remainder of the paper is organized as follows.
First, Section~\ref{SEC_Preliminaries} prepares the regularity framework: we introduce
approximate analogues of the exact synthesis problems (Problems~\ref{PROBLEM_ASS} and~\ref{PROBLEM_AUS},
vs.~\ref{PROBLEM_SS} and~\ref{PROBLEM_US}), show that the relaxed dynamics enjoy convenient
regularity properties, and prove
that, although the exact and approximate controllability times differ pointwise, their
worst-case values coincide. This provides the working framework for the constructive
bounds of Section~\ref{SEC_Upper_bound}.

Second, Section~\ref{SEC_Upper_bound} proves
Proposition~\ref{PRO_Result_Constructive_State_to_State}
(Section~\ref{SEC_Upper_bound_StS}) and
Proposition~\ref{PRO_Result_Explicit_Propagator}
(Section~\ref{SEC_Contructive_bounds_Propagators}), via convexification and averaging
procedures (Section~\ref{SUBSEC_Edge_Isolation}).

Third, Section~\ref{SEC_Chiral} analyzes finite-dimensional Galerkin truncations
of the infinite-dimensional planar rotator and proves
Proposition~\ref{PRO_rotateur_introduction}.

Finally, Appendix~\ref{SEC_APP_lower_bounds} contains the proofs of the lower bounds
(Sections~\ref{SEC_Spectral_bounds}, \ref{SEC_Chordal}, \ref{SEC_Ehrenfest}), and
Appendix~\ref{SEC_Appendix} collects the linear algebra computations.

\section{An auxiliary relaxed problem}\label{SEC_Preliminaries}

\subsection{Approximate vs exact transfer times}\label{SUBSEC_approx_vs_exact}
We are interested in the exact synthesis times
$\tau_{\mathrm{SS}}$ and $\tau_{\mathrm{US}}$ of
Problems~\ref{PROBLEM_SS} and~\ref{PROBLEM_US} in the framework of Assumptions \ref{ASS_control} and \ref{ASS_coupling_graph}.
However, the regularity of these maps
($X_f\mapsto\tau_{\mathrm{US}}(X_f)$,
$\psi_f\mapsto\tau_{\mathrm{SS}}(\psi_0,\psi_f)$) is delicate: 
the exact attainable sets 
\begin{equation}\label{EQ_def_attainable_exact}
\mathcal{A}_{\leq T} := \bigcup_{0 \leq t \leq T} \bigcup_{u \in L^1([0,t],\R)}
\{ X^u_{t,0} \} 
\end{equation} 
are not closed in general,
the infima need not be minima, and the maps need not be
continuous (Remark~\ref{REM_Blow-up}) under these assumptions alone.
We therefore introduce the approximate
formulations
(Problems~\ref{PROBLEM_ASS}--\ref{PROBLEM_AUS}).
They are easier to handle: their worst-case values coincide with those
of the exact problems (Lemma~\ref{LEM_approx_cont}), so the bounds
proved in this paper for the approximate times transfer to the exact ones.

\begin{problem}[Approximate State-to-State Synthesis (ASS)]\label{PROBLEM_ASS}
Given $\psi_0, \psi_f \in \mathbb{S}^{2d-1}$, define
$$\tau_{\mathrm{ASS}}(\psi_0, \psi_f) = \inf\bigl\{T > 0 \mid \forall\,\varepsilon > 0,\; \exists\, u \in L^1_{\mathrm{loc}},\; \exists\, t \leq T,\; \|X^u_{t,0}\psi_0 - \psi_f\| < \varepsilon\bigr\},$$
where $\|\cdot\|$ is the canonical Hermitian norm on $\mathbb{C}^d$. 

The worst-case value is $\tau_{\mathrm{ASS}}^{\max}=\sup_{\psi_0,\psi_f \in \mathbb{S}^{2d-1}} \tau_{\mathrm{ASS}}(\psi_0,\psi_f).$
\end{problem}

\begin{problem}[Approximate Unitary Operator Synthesis (AUS)]\label{PROBLEM_AUS}
Given $X_f \in G$, define
$$\tau_{\mathrm{AUS}}(X_f) = \inf\bigl\{T > 0 \mid \forall\,\varepsilon > 0,\; \exists\, u \in L^1_{\mathrm{loc}},\; \exists\, t \leq T,\; \|X^u_{t,0} - X_f\|_{\mathrm{op}} < \varepsilon\bigr\},$$
where $\|\cdot\|_{\mathrm{op}}$ is the operator norm on $\mathbb{C}^{d\times d}$.
The worst-case value (time diameter of the group) is $\tau_{\mathrm{AUS}}^{\max} := \sup_{X_f \in G} \tau_{\mathrm{AUS}}(X_f).$
\end{problem}

\subsection{Equivalence of the worst-case estimates}\label{SUB_SEC_equivalence_worst_case}

By definition, the approximate transfer times lower-bound the exact ones: \[\tau_{\mathrm{ASS}}(\psi_0,\psi_f) \leq \tau_{\mathrm{SS}}(\psi_0,\psi_f) \quad \mbox{ and } \quad \tau_{\mathrm{AUS}}(X_f) \leq \tau_{\mathrm{US}}(X_f).\] Lower bounds on the approximate time therefore apply to the exact time. The reverse inequality holds for the worst-case values under Assumption~\ref{ASS_control} (Lemma~\ref{LEM_approx_cont}).

\begin{lemma}\label{LEM_approx_cont}
Under Assumption \ref{ASS_control}, $\tau_{\mathrm{SS}}^{\max} = \tau_{\mathrm{ASS}}^{\max}$ and $\tau_{\mathrm{US}}^{\max} = \tau_{\mathrm{AUS}}^{\max}$.
\end{lemma}
\begin{proof}
We prove $\tau_{\mathrm{US}}^{\max} \leq \tau_{\mathrm{AUS}}^{\max}$; the reverse inequality is immediate from the definitions, and the state-to-state case follows identically on $\mathbb{S}^{2d-1}$.

Fix $T > \tau_{\mathrm{AUS}}^{\max}$. By definition, for every $X_f \in G$ and every $\varepsilon > 0$, there exists $u \in L^1_{\mathrm{loc}}$ and $t \leq T$ such that $\|X^u_{t,0} - X_f\|_{\mathrm{op}} < \varepsilon$. Hence, the attainable set $\mathcal{A}_{\leq T}$ (defined in \eqref{EQ_def_attainable_exact}) is dense in $G$. Now fix $X_f \in G$ and $\varepsilon > 0$. We show $X_f \in \mathcal{A}_{\leq T+\varepsilon}$.

\textbf{Step 1 (open set via strong accessibility).} Consider the reversed system $\dot{w} 
= -(A + \tilde{v}B)\,w$ with $w(0) = X_f$ and $\tilde{v} \in L^1([0,\varepsilon])$. Passing 
to the right-invariant formulation via $z = w^{-1}$ yields a system on $G$ whose 
associated Lie algebra is $\operatorname{Lie}(A, B) \supset \mathfrak{g}$ by Assumption 
\ref{ASS_control}. The strong accessibility property implies that the set of endpoints 
$U_\varepsilon := \{w(s) : \tilde{v}\in L^1([0,\varepsilon]),\; 
0 \leq s \leq \varepsilon\}$ has a non-empty interior in $G$.

\textbf{Step 2 (intersection by density).} Since $\mathcal{A}_{\leq T}$ is dense in $G$ and $\operatorname{int}(U_\varepsilon)$ is a non-empty open subset, we can fix an element $h \in \mathcal{A}_{\leq T} \cap \operatorname{int}(U_\varepsilon)$ which is not empty.

\textbf{Step 3 (concatenation).} Since $h \in \mathcal{A}_{\leq T}$, there exists a control
 $u_1 \in L^1([0,t_1])$ with $t_1 \leq T$ such that $X^{u_1}_{t_1,0} = h$. Since 
 $h \in U_\varepsilon$, there exists $\tilde{v} \in L^1([0,\varepsilon])$ and $s \leq \varepsilon$ 
such that $w(0) = X_f$ and $w(s) = h$ for the reversed system. Reversing time, the control 
$v(t) := \tilde{v}(s - t) \in L^1([0,s])$ satisfies $X^{v}_{s,0} h = X_f$ along the original 
dynamics. Concatenating $u_1$ on $[0,t_1]$ and $v$ on $[t_1, t_1+s]$ yields a control 
reaching exactly $X_f$ in time $t_1 + s \leq T + \varepsilon$. Hence the target $X_f$ 
belongs to the  $ \mathcal{A}_{\leq T+\varepsilon}$, the attainable set of the original 
system in time less than $T+\varepsilon$.

Letting $\varepsilon$ tend to  $ 0$ and $T$ to $ \tau_{\mathrm{AUS}}^{\max}$ yields $\tau_{\mathrm{US}}^{\max} \leq \tau_{\mathrm{AUS}}^{\max}$.
\end{proof}

\subsection{Relaxed dynamics}\label{SUB_SEC_Relaxed_dynamics}

While the worst-case exact and approximate times coincide under Assumption \ref{ASS_control}, the pointwise equivalence $\tau_{\mathrm{US}}(X_f) = \tau_{\mathrm{AUS}}(X_f)$ requires stronger topological properties. To analyze this, we introduce the relaxed system dynamics.

Following \cite{BOUSSAID2020108412}, define the compact set
\begin{equation}\label{eq:CC_def}
  CC \;:=\; \overline{\operatorname{conv}\bigl(\{e^{-vB}Ae^{vB} \mid v\in\mathbb{R}\}\bigr)}
  \;\subset\; \mathfrak{g},
\end{equation}
the auxiliary control system 
\begin{equation}\label{eq:aux_sys}
  \dot{Y}(t) \;=\; \Omega(t)\,Y(t), \qquad Y(0)=I_d,
\end{equation}
where the control $\Omega$ ranges over the set of admissible controls
$L^{\infty}([0,T],\,CC)$.
\begin{definition}[Generalized solutions]\label{DEF_Generalized_solution}
The propagator of \eqref{eq:aux_sys} between times $s$ and $t$ is denoted
$Y^{\Omega}_{t,s}$.
For any $R\in L^{\infty}([0,T],\,\overline{\exp(\mathbb{R}B)})$, the
generalized solution of \eqref{EQ_main} associated with the pair
$(\Omega,R)$ is
\begin{equation}\label{eq:gen_sol}
  Z(t) \;=\; R(t)\;Y^{\Omega}_{t,0}.
\end{equation}
We denote by  $\mathcal{A}^{\mathrm{rel}}_{\leq T}$  the set of all
$Z(t)\in G$ obtained this way for $t \in [0,T]$, i.e.
\[ \mathcal{A}^{\mathrm{rel}}_{\leq T}
 = \bigcup_{0\leq t \leq T}\{R(t)\,Y^{\Omega}_{t,0} : \Omega\in L^{\infty}([0,t],CC),\;
   R\in L^{\infty}([0,t],\overline{\exp(\mathbb{R}B)})\}.\]
\end{definition}

\paragraph{Decomposition of relaxed trajectories.}
Write $A = A_D + A_{ND}$ with $A_D$ diagonal (commuting with $B$).
Every generalized solution $Z:[0,T]\to G$ of \eqref{eq:gen_sol},
$Z(t)=R(t)\,Y^{\Omega}_{t,0}$ with $\Omega\in L^{\infty}([0,T],CC)$,
admits the further factorisation
\begin{equation}\label{eq:decomp}
  Z(t) \;=\; \underbrace{R(t)\,\exp\!\bigl(t\,A_D\bigr)}_{S_R(t)}
  \;\cdot\; \underbrace{W^{\widetilde{\Omega}}(t)}_{W},
\end{equation}
where $W^{\widetilde{\Omega}}(t):=\exp(-t A_D)\,Y^{\Omega}_{t,0}$ and
$\widetilde{\Omega}(t):=\exp(-t A_D)\bigl(\Omega(t)-A_D\bigr)\exp(t A_D)
\in CC(t)$, with the time-dependent control set
$CC(t)=\exp(-t A_D)(CC-A_D)\exp(t A_D)$.  The factor $W^{\widetilde{\Omega}}$
satisfies the effective interaction dynamics
\begin{equation}\label{eq:W_dynamics}
  \dot{W}^{\widetilde{\Omega}}(t) \;=\;
  \widetilde{\Omega}(t)\,W^{\widetilde{\Omega}}(t),
  \qquad \widetilde{\Omega}(t)\in CC(t)\subset\mathfrak{su}(d),
\end{equation}
where $CC(t)=\overline{\operatorname{conv}\{e^{-wB-t A_D}A_{ND}
e^{wB+t A_D}:w\in\mathbb{R}\}}$ is compact and
$\|\Omega(t)\|_{op}\leq L_{CC}$.  Here the torus factor $R$ and the
generator $\widetilde{\Omega}$ are \emph{independent} choices, and the
diagonal drift factor $\exp(tA_D)$ depends only on the time $t$ (not on
the control nor on $R$).

\begin{lemma}\label{LEM_cloture_ensemble_attaignable}
For every $T\geq 0$, $\mathcal{A}^{\mathrm{rel}}_{\leq T}
= \overline{\mathcal{A}_{\leq T}}$,
where the closure is taken in the operator norm.
\end{lemma}

\begin{proof}
\textbf{Inclusion $\mathcal{A}_{\leq T}\subset\mathcal{A}^{\mathrm{rel}}_{\leq T}$.}
Let $X_f\in\mathcal{A}_{\leq T}$, reached by a control
$u\in L^{1}([0,T])$ and the corresponding propagator
$X^{u}_{t,0}$.  Set
$\theta(t):=\int_0^t u(s)\,ds$ (absolutely continuous) and
$R(t):=e^{\theta(t)B}\in\overline{\exp(\mathbb{R}B)}$.
Define $Y(t):=e^{-\theta(t)B}\,X^{u}_{t,0}$.  Then
\begin{align*}
  \dot{Y}(t)
  &= -u(t)B\,e^{-\theta(t)B}X^{u}_{t,0}
     + e^{-\theta(t)B}\bigl(A+u(t)B\bigr)X^{u}_{t,0} \\
  &= e^{-\theta(t)B}\bigl(A+u(t)B\bigr)e^{\theta(t)B}\,Y(t)
     - u(t)B\,Y(t) \\
  &= e^{-\theta(t)B}A\,e^{\theta(t)B}\,Y(t).
\end{align*}
Thus $Y$ satisfies \eqref{eq:aux_sys} with
$\Omega(t)=e^{-\theta(t)B}Ae^{\theta(t)B}$.  Since
$\Omega(t)\in\{e^{-vB}Ae^{vB}:v\in\mathbb{R}\}\subset CC$ and
$\|\Omega(t)\|_{op}\leq\|A\|_{op}$,
we have $\Omega\in L^{\infty}([0,T],CC)$.  Moreover
$X^{u}_{T,0}=R(T)\,Y(T)$, so $X_f\in\mathcal{A}^{\mathrm{rel}}_{\leq T}$.

\medskip
\noindent
\textbf{Closedness of $\mathcal{A}^{\mathrm{rel}}_{\leq T}$.}
Let $X_n\in\mathcal{A}^{\mathrm{rel}}_{\leq T}$ with
$X_n=R_n(t_n)\,Y_n(t_n)\to X$ in $\|\cdot\|_{op}$, where 
$\Omega_n\in L^{\infty}([0,t_n],CC)$ and
$R_n\in L^{\infty}([0,t_n],\overline{\exp(\mathbb{R}B)})$, $t_n\in [0,T]$.
We extend the functions $\Omega_n$ (defined on $[0,t_n]$ by assumption) to $[0,T]$ by defining $\Omega_n(t):=\Omega(t_n)$ for $t\in (t_n,T]$.

By compactness of $[0,T]$ we may assume, up to extraction, that $t_n \to t_\infty$ with
$t_\infty \in [0,T]$. 
Since $CC$ is compact (hence bounded), the sequence
$(\Omega_n)$ is bounded in $L^{\infty}([0,T],\mathfrak{g})$; by
Banach--Alaoglu theorem (dual of $L^{1}$), a subsequence satisfies
$\Omega_n\,\rightharpoonup^{*}\,\bar{\Omega}$ in $L^{\infty}$ weak-$*$,
with $\bar{\Omega}(t)\in CC$ a.e.\ (the pointwise constraint is
preserved under weak-$*$ limits since $CC$ is compact and convex).

Each $Y_n$ satisfies $\dot{Y}_n=\Omega_n Y_n$ with
$\|\dot{Y}_n\|_{op}\leq L_{CC}\|Y_n\|_{op}$
($L_{CC}:=\sup_{\Omega\in CC}\|\Omega\|_{op}$).  By Gr\"onwall lemma,
$\|Y_n(t)\|_{op}\leq e^{L_{CC}t}$, so $(Y_n)_n$ is uniformly bounded and
equicontinuous in $C^{0}([0,T],G)$.  By Arzel\`a--Ascoli theorem, a
subsequence converges uniformly to some $\bar{Y}\in C^{0}([0,T],G)$.
Passing to the limit in the weak formulation
$\int_0^T\langle\dot{Y}_n,\varphi\rangle\,ds
 = \int_0^T\langle\Omega_n Y_n,\varphi\rangle\,ds$
(valid for smooth test functions $\varphi$), and using
$\Omega_n\,\rightharpoonup^{*}\,\bar{\Omega}$ and
$Y_n\to\bar{Y}$ uniformly, we obtain
$\dot{\bar{Y}}=\bar{\Omega}\,\bar{Y}$. Noting that $Y_n(0)=I_d$ for every $n$, one deduces that $\bar{Y}(0)=I_d$. Thus $\bar{Y}$ is indeed the function $ t\mapsto Y^{\bar{\Omega}}_{t,0}$. In particular, $Y_n(t_n)\to Y^{\bar{\Omega}}_{t_\infty,0}$.

Since $R_n(t_n)$ takes value in the compact set $\overline{\exp(\mathbb{R}B)}$,
 one may assume (up to extraction) that $R_n(t_n)$
converges to some $\bar{R}\in\overline{\exp(\mathbb{R}B)}$.

Hence $X=\lim_n R_n(t_n) Y_n(t_n)= \bar{R}\,\bar{Y}(t_\infty)\in\mathcal{A}^{\mathrm{rel}}_{\leq T}$.
Since $\mathcal{A}^{\mathrm{rel}}_{\leq T}\subset G$ (compact) and
is closed, it is compact.

\medskip
\noindent
\textbf{Reverse inclusion $\overline{\mathcal{A}_{\leq T}}
\supset\mathcal{A}^{\mathrm{rel}}_{\leq T}$.}
Since $\mathcal{A}^{\mathrm{rel}}_{\leq T}$ is closed and
$\mathcal{A}_{\leq T}\subset\mathcal{A}^{\mathrm{rel}}_{\leq T}$,
it suffices to show that
$\mathcal{A}_{\leq T}$ is \emph{dense} in
$\mathcal{A}^{\mathrm{rel}}_{\leq T}$.

We now prove the density inclusion
$\mathcal A^{\mathrm{rel}}_{\le T}\subset \overline{\mathcal A_{\le T}}$.
Fix a relaxed pair
$(\Omega, R)$ with $\Omega \in L^{\infty}([0,T],CC)$ and
$R \in L^{\infty}([0,T],\overline{\exp(\mathbb R B)})$, and set the target
endpoint to $Z(T) := R(T)\,Y^{\Omega}_{T,0}$. We construct a sequence of
exact controls $u$ whose endpoints converge to $Z(T)$ in the operator norm.
Throughout, $\varepsilon > 0$ is a fixed accuracy parameter.

\medskip
\noindent
\textbf{Step (i): Approximation of the relaxed generator by orbit-valued step
functions.}
Write $\gamma := \{e^{-vB}Ae^{vB} : v \in \mathbb R\}$ for the orbit of $A$
under the control flow, and $\bar\gamma := \overline{\gamma}$ for its
operator-norm closure, so that $CC = \overline{\operatorname{conv}(\gamma)}
= \overline{\operatorname{conv}(\bar\gamma)}$. We claim that for every
$\varepsilon_1 > 0$ there exists a step function
$\Omega^{(n)} : [0,T] \to \bar\gamma$, constant on a partition
$0 = t_0 < t_1 < \cdots < t_N = T$, such that
\begin{equation}\label{eq:step1_filippov}
  \|\Omega^{(n)} - \Omega\|_{L^1([0,T],\mathfrak g)} < \varepsilon_1 .
\end{equation}
Indeed, this is the $L^1$-closure characterization
\[\overline{\{\,g \in L^1 : g(t) \in \bar\gamma \text{ a.e.}\,\}}
= \{\,f \in L^1 : f(t) \in \overline{\mathrm{conv}(\bar\gamma)}
= CC \text{ a.e.}\,\},\] a standard consequence of Filippov's approximation
theorem (Theorem 8.2 in \cite{Agrachev2004}). For completeness we recall the construction. Since $CC$ is a compact
convex subset of the finite-dimensional space
$\mathfrak g \simeq \mathbb R^{2d^2}$, cover $CC$ by cubes of mesh $\delta$;
for each cube $\mathcal C_k$ meeting the range of $\Omega$ set
$E_k := \{t : \Omega(t) \in \mathcal C_k\}$ (made pairwise disjoint by a
standard tie-breaking) and pick $\omega_k \in CC \cap \mathcal C_k$; the
simple function $\omega^{(K)} := \sum_k \omega_k\,\mathbf 1_{E_k}$ satisfies
$\|\omega^{(K)} - \Omega\|_{L^1} \le \delta\,T$. For each $k$, Carath\'eodory's
theorem (in $\mathbb R^{2d^2}$) together with the definition of $CC$ yields
points $\omega_{k,1}, \dots, \omega_{k,m_k} \in \bar\gamma$ with
$m_k \le 2d^2 + 1$ and weights $\lambda_{k,j} \ge 0$,
$\sum_j \lambda_{k,j} = 1$, such that
$\|\omega_k - \sum_j \lambda_{k,j}\,\omega_{k,j}\| < \eta$. Slicing $E_k$
into measurable sets $E_{k,j} \subset E_k$ with
$|E_{k,j}| = \lambda_{k,j}\,|E_k|$ and setting
$\Omega^{(n)} := \sum_{k,j}\omega_{k,j}\,\mathbf 1_{E_{k,j}}$ gives a
$\bar\gamma$-valued step function with
$\|\Omega^{(n)} - \omega^{(K)}\|_{L^1} \le \eta\,T$. Letting
$\delta, \eta \to 0$ along a diagonal sequence $n \to \infty$ yields
\eqref{eq:step1_filippov}.

\medskip
\noindent
\textbf{Step (ii): Lifting to a gauge angle and approximation of the
interaction propagator.}
On each subinterval $I_k = [t_k, t_{k+1})$ of constancy of $\Omega^{(n)}$
(value $\omega_k \in \bar\gamma$), choose $v_{n,k} \in \mathbb R$ such that
\begin{equation}\label{eq:step2_lifting}
  \bigl\|e^{-v_{n,k}B}A\,e^{v_{n,k}B} - \omega_k\bigr\|_{op}
  \;\le\; \tfrac{1}{p},
\end{equation}
which is possible since $\bar\gamma$ is the operator-norm closure of the
orbit $\gamma$. Let $v_{n,p}$ denote the piecewise-constant function
$v_{n,p} \equiv v_{n,k}$ on $I_k$, and set
$\Omega^{(n,p)}(t) := e^{-v_{n,p}(t)B}A\,e^{v_{n,p}(t)B}$; then
$\|\Omega^{(n,p)} - \Omega^{(n)}\|_{L^1} \le T/p$.

We now lift $v_{n,p}$ to a smooth gauge angle. Since the $C^\infty$
functions vanishing at $t = 0$ are dense in $L^1([0,T])$ (a mollification of
$v_{n,p}$, followed by subtraction of a smooth cutoff supported in
$[0,\delta]$ to enforce the value $0$ at $t=0$, changes the $L^1$ norm by
$O(\delta)$), there exists $\theta_{n,p,q} \in C^\infty([0,T])$ with
$\theta_{n,p,q}(0) = 0$ and
\begin{equation}\label{eq:gauge_smooth}
  \|\theta_{n,p,q} - v_{n,p}\|_{L^1([0,T])} < \varepsilon_3 .
\end{equation}
Set the control $u_{n,p,q}(t) := \dot\theta_{n,p,q}(t) \in L^1([0,T])$ and
let $Y^{n,p,q}_{t,0}$ be the interaction propagator of the exact dynamics
$\dot X = (A + u_{n,p,q}B)X$, i.e.\ the solution of
$\dot Y = \Omega^{\mathrm{eff}}_{n,p,q}(t)\,Y$, $Y_{0,0} = I_d$, where the
\emph{effective} (interaction) generator is
\begin{equation}\label{eq:eff_gen}
  \Omega^{\mathrm{eff}}_{n,p,q}(t)
  \;:=\; e^{-\theta_{n,p,q}(t)B}\,A\,e^{\theta_{n,p,q}(t)B} .
\end{equation}
(Here $\theta(t) = \int_0^t u_{n,p,q} = \theta_{n,p,q}(t)$ since
$\theta_{n,p,q}(0) = 0$.)

Since $B$ is skew-adjoint, standard functional calculus using $e^{S} - I_d = \int_0^1 e^{(1-s)S}S\,ds$ provides, for every $w$ in $\R$,
\[\| e^{w B}-I_d \|_{op} \leq \int_{s=0}^1 \left \| e^{(1-s) w B} w B \right \|_{op} \mathrm{ds}
\leq |w | \| B\|_{op}.\]
Hence, for every $v$ and $v'$ in $\R$,
\[\|e^{vB}-e^{v'B}\| \leq \|e^{v' B}\|_{op} \|e^{(v-v')B} - I_d\|_{op} \le |v-v'|\,\|B\|_{op}.\]
Thus, defining $C_{\mathrm{lip}} := 2\,\|A\|_{op}\,\|B\|_{op}$,
\begin{eqnarray}
\bigl\|e^{-vB} Ae ^{vB} - e^{-v'B} A e^{v'B} \bigr \|_{op} & \leq &
\left  \| e^{-vB}A(e^{vB} - e^{v'B}) + (e^{-vB} - e^{-v'B})A\,e^{v'B} 
\right \|_{op} \nonumber \\
& \leq & \left  \| e^{-vB}A \right \|_{op} \left \| e^{vB} - e^{v'B} \right \| 
+   \left \| e^{vB} - e^{v'B} \right \|_{op}  \left \| A e^{vB} \right \|_{op}
\nonumber \\
& \leq & C_{\mathrm{lip}}\,|v - v'|.
   \label{eq:lip_orbit}
\end{eqnarray}

By \eqref{eq:lip_orbit} and \eqref{eq:gauge_smooth}, the effective generator
is close to the target $\Omega^{(n,p)}$ in $L^1$:
\begin{equation}\label{eq:gen_close}
  \int_0^T \bigl\|\Omega^{\mathrm{eff}}_{n,p,q}(t) - \Omega^{(n,p)}(t)\bigr\|_{op}\,dt
  \;\le\; C_{\mathrm{lip}}\,\|\theta_{n,p,q} - v_{n,p}\|_{L^1}
  \;<\; C_{\mathrm{lip}}\,\varepsilon_3 .
\end{equation}
The Duhamel identity
\begin{equation}\label{eq:duhamel}
  Y^{n,p,q}_{T,0} - Y^{\Omega^{(n,p)}}_{T,0}
  \;=\; \int_0^T Y^{\Omega^{(n,p)}}_{T,t,0}\,
        \bigl(\Omega^{\mathrm{eff}}_{n,p,q}(t) - \Omega^{(n,p)}(t)\bigr)\,
        Y^{n,p,q}_{t,0}\;dt
\end{equation}
and the unitarity of the propagators give
\begin{equation}\label{eq:prop_close}
  \bigl\|Y^{n,p,q}_{T,0} - Y^{\Omega^{(n,p)}}_{T,0}\bigr\|_{op}
  \;\le\; \int_0^T \bigl\|\Omega^{\mathrm{eff}}_{n,p,q}(t)
  - \Omega^{(n,p)}(t)\bigr\|_{op}\,dt
  \;\le\; C_{\mathrm{lip}}\,\varepsilon_3 .
\end{equation}
Combining \eqref{eq:prop_close} with the relaxed  Duhamel estimates
$\|Y^{\Omega^{(n,p)}}_{T,0} - Y^{\Omega^{(n)}}_{T,0}\|_{op} \le
\|\Omega^{(n,p)} - \Omega^{(n)}\|_{L^1} \le T/p$ and
$\|Y^{\Omega^{(n)}}_{T,0} - Y^{\Omega}_{T,0}\|_{op} \le
\|\Omega^{(n)} - \Omega\|_{L^1} < \varepsilon_1$, and choosing
$\varepsilon_3 < \varepsilon/(3C_{\mathrm{lip}})$, $p > 3T/\varepsilon$,
$\varepsilon_1 < \varepsilon/3$, we obtain
\begin{equation}\label{eq:interaction_close}
  \bigl\|Y^{n,p,q}_{T,0} - Y^{\Omega}_{T,0}\bigr\|_{op} \;<\; \varepsilon/3 .
\end{equation}

\medskip
\noindent
\textbf{Step (iii): Terminal correction of the gauge.}
With the control $u_{n,p,q}$ on $[0,T]$, the exact endpoint reads
$X^{u_{n,p,q}}_{T,0} = e^{\theta_{n,p,q}(T)B}\,Y^{n,p,q}_{T,0}$; the
interaction factor is already $\varepsilon/3$-close to $Y^{\Omega}_{T,0}$ by
\eqref{eq:interaction_close}, but the gauge $e^{\theta_{n,p,q}(T)B}$ need not
coincide with the target $R(T)$. We remedy this by a short, fast rotation on
a terminal interval $[T-\delta, T]$ (with $\delta > 0$ to be sent to $0$).

Define the correction factor
$\Delta R := R(T)\,e^{-\theta_{n,p,q}(T-\delta)B} \in \overline{\exp(\mathbb R B)}$.
Since $\overline{\exp(\mathbb R B)}$ is the operator-norm closure of the
one-parameter subgroup $\exp(\mathbb R B)$, for every $\eta > 0$ there
exists $v_{\mathrm{corr}} \in \mathbb R$ such that
\begin{equation}\label{eq:corr_R}
  \bigl\|e^{v_{\mathrm{corr}}B} - \Delta R\bigr\|_{op} \;\le\; \eta .
\end{equation}
On $[T-\delta, T]$ set the control to the constant
$u_{\mathrm{corr}}(t) := v_{\mathrm{corr}}/\delta$, and let
$u_{\mathrm{full}}$ denote the concatenation of $u_{n,p,q}$ on
$[0, T-\delta]$ with $u_{\mathrm{corr}}$ on $[T-\delta, T]$. The terminal
gauge of the full trajectory is
$R_{\mathrm{full}}(T) = e^{\theta_{n,p,q}(T-\delta)B}\,e^{v_{\mathrm{corr}}B}$,
hence
\begin{equation}\label{eq:R_final_err}
  \bigl\|R_{\mathrm{full}}(T) - R(T)\bigr\|_{op}
  \;=\; \bigl\|e^{\theta_{n,p,q}(T-\delta)B}\bigl(e^{v_{\mathrm{corr}}B}
  - \Delta R\bigr)\bigr\|_{op}
  \;\le\; \eta .
\end{equation}
On $[T-\delta, T]$ the effective generator of the full trajectory is
$\Omega_{\mathrm{corr}}(t) = e^{-\theta_{\mathrm{full}}(t)B}A\,
e^{\theta_{\mathrm{full}}(t)B}$ with
$\|\Omega_{\mathrm{corr}}(t)\|_{op} = \|A\|_{op}$ (unitary invariance); the
Duhamel identity on $[T-\delta, T]$ gives
\begin{equation}\label{eq:Y_corr_err}
  \bigl\|Y_{\mathrm{full}}(T) - Y^{n,p,q}(T)\bigr\|_{op}
  \;\le\; \int_{T-\delta}^{T}
         \bigl\|\Omega_{\mathrm{corr}}(t) - \Omega^{\mathrm{eff}}_{n,p,q}(t)\bigr\|_{op}\,dt
  \;\le\; 2\,\delta\,\|A\|_{op} .
\end{equation}
Choosing $\eta < \varepsilon/3$ and $\delta < \varepsilon/(6\|A\|_{op})$,
both corrections are $< \varepsilon/3$.

Finally, at time $T$,
$X^{u_{\mathrm{full}}}_{T,0} = R_{\mathrm{full}}(T)\,Y_{\mathrm{full}}(T)$,
and the unitarity of the factors yields
\[
  \bigl\|X^{u_{\mathrm{full}}}_{T,0} - R(T)\,Y^{\Omega}_{T,0}\bigr\|_{op}
  \;\le\; \underbrace{\bigl\|Y_{\mathrm{full}}(T) - Y^{n,p,q}(T)\bigr\|_{op}}_{<\,\varepsilon/3}
  \;+\; \underbrace{\bigl\|Y^{n,p,q}_{T,0} - Y^{\Omega}_{T,0}\bigr\|_{op}}_{<\,\varepsilon/3}
  \;+\; \underbrace{\bigl\|R_{\mathrm{full}}(T) - R(T)\bigr\|_{op}}_{<\,\varepsilon/3},
\]
hence $\bigl\|X^{u_{\mathrm{full}}}_{T,0} - R(T)\,Y^{\Omega}_{T,0}\bigr\|_{op} < \varepsilon$.
Since $\varepsilon > 0$ is arbitrary, the relaxed endpoint
$R(T)\,Y^{\Omega}_{T,0}$ belongs to $\overline{\mathcal A_{\le T}}$, which
establishes the density inclusion
$\mathcal A^{\mathrm{rel}}_{\le T} \subset \overline{\mathcal A_{\le T}}$.

Since
$\mathcal{A}^{\mathrm{rel}}_{\leq T}\subset\overline{\mathcal{A}_{\leq T}}$
and $\overline{\mathcal{A}_{\leq T}}\subset\mathcal{A}^{\mathrm{rel}}_{\leq T}$,
the equality follows.
\end{proof}

The quantities $\tau_{\mathrm{ASS}}(\psi_0,\psi_f)$ and $\tau_{\mathrm{AUS}}(X_f)$, defined as infimums, are genuine minimums under the sole Assumption \ref{ASS_control}:

\begin{corollary}[Attainment of the approximate synthesis time]
\label{COR_AUS_attainment}
Under Assumption~\ref{ASS_control}, for every $X_f\in G$,
the infimum defining $\tau_{\mathrm{AUS}}(X_f)$ is a minimum:
there exists $\alpha=\tau_{\mathrm{AUS}}(X_f)$ and a relaxed
trajectory $(\bar{\Omega},\bar{R})$ such that
$\bar{R}(\alpha)\,Y^{\bar{\Omega}}_{\alpha,0}=X_f$.
Analogously, $\tau_{\mathrm{ASS}}(\psi_0,\psi_f)$ is attained.
\end{corollary}

\begin{proof}
By Lemma~\ref{LEM_cloture_ensemble_attaignable},
$\tau_{\mathrm{AUS}}(X_f)=\inf\{s>0:X_f\in\mathcal{A}^{\mathrm{rel}}_{\leq s}\}$.
Let $\alpha:=\tau_{\mathrm{AUS}}(X_f)$. The set
$\mathcal{A}^{\mathrm{rel}}_{\leq s}$ is compact for every $s$
(proved in Lemma~\ref{LEM_cloture_ensemble_attaignable}), and the
family is nested by definition: $\mathcal{A}^{\mathrm{rel}}_{\leq s}\subset
\mathcal{A}^{\mathrm{rel}}_{\leq s'}$ for $s\leq s'$.

It suffices to show that
$\mathcal{A}^{\mathrm{rel}}_{\leq\alpha}
=\bigcap_{n}\mathcal{A}^{\mathrm{rel}}_{\leq(\alpha+1/n)}$.
The inclusion $\mathcal{A}^{\mathrm{rel}}_{\leq\alpha}
\subset \bigcap_{n}\mathcal{A}^{\mathrm{rel}}_{\leq(\alpha+1/n)}$ is obvious. To show the reciprocal inclusion, consider $X_f$ in $\bigcap_{n}\mathcal{A}^{\mathrm{rel}}_{\leq(\alpha+1/n)}$. By definition, there exists a sequence $(t_n,R_n,\Omega_n)_n$ such that, for every $n$,  $0 \leq t_n \leq \alpha +1/n$, $R_n \in L^\infty([0,t_n],\overline{\exp(\R B)})$, $\Omega_n \in  L^\infty([0,t_n],CC)$ and $X_f= R_n(t_n) Y^{\Omega_n}_{t_n,0}$. If  $t_n \leq \alpha$ for some $n$, then there is nothing to prove. We assume $t_n >\alpha$ for every $n$. Since $(t_n)_n$ takes value in the compact set $ [0,\alpha +1]$,
one may assume, up to extraction, that $t_n$ converges to a certain limit $ t \in [0,\alpha +1]$. Moreover, $t = \alpha$ (since $\limsup t_n \leq \limsup \alpha +1/n =\alpha $ and $t_n > \alpha$ for every $n$). 

Considering the restrictions to $[0,\alpha]$, and by the same extraction process as the one used in Lemma \ref{LEM_cloture_ensemble_attaignable}, one extracts a  subsequence of $\Omega_n$ (still denoted $\Omega_n$ with an abuse of notations) $\ast$-weakly convergent to $\bar{\Omega}:[0,\alpha]\to CC$. By definition, the sequence
\begin{equation}\label{EQ_extraction_sol_gen}
X_f   \left (Y^{\Omega_n}_{t_n,0} \right )^{-1} = R_n(t_n)
\end{equation} belongs to the compact set $\overline{\exp(\R B)}$. Up to another extraction, one may assume that  $R_n(t_n)$ converges to a limit 
$\bar{R}(\alpha)$  in  $\overline{\exp(\R B)}$. Passing to the limit $n\to \infty$ in \eqref{EQ_extraction_sol_gen}, one gets
$\bar{R}(\alpha)=X_f  \left (Y^{\bar{\Omega}}_{\alpha,0} \right )^{-1}$
 or equivalently $X_f= \bar{R}(\alpha) Y^{\bar{\Omega}}_{\alpha,0}$, hence 
$X_f \in  \mathcal{A}^{\mathrm{rel}}_{\leq\alpha}$,
i.e., the infimum time $\alpha$ is attained.
\end{proof}

\subsection{Pointwise discrepancy between exact and approximate control times}\label{SUB_SEC_Discrepancy_exact_approximate}

We aim to prove that, in general, $\tau_{\mathrm{AUS}}$ and $\tau_{\mathrm{US}}$ are different. The following technical lemma will be used to build a counter-example.

\begin{lemma}[Cartan part of the interaction propagator]\label{LEM_cartan_part}
Work in the eigenbasis of $B$, with $B=\mathrm{diag}(ib_1,\ldots,ib_d)$,
 $b_j\in\mathbb{R}$. In the same basis, write
$A=A_D+A_{ND}$, where $A_{ND}$ has a zero diagonal and $A_D=\mathrm{diag}(ia_1,\ldots,ia_d)$, $a_j\in\mathbb{R}$.
Let $\mathfrak{c}\subset\mathfrak{gl}(d,\mathbb{C})$ denote the Cartan subalgebra of diagonal
matrices in this basis, $P_{\mathfrak c}$ the orthogonal projection onto $\mathfrak c$ for the
Frobenius inner product, and $L:=\|A_{ND}\|_{op}>0$. For $u\in L^1_{loc}([0,\infty))$, set
$v(t)=\int_0^t u(s)\,ds$ and let $W^u_{t,0}$ solve the interaction equation
\begin{equation}\label{EQ_interaction_W}
\dot W^u_{t,0}=M^u(t)\,W^u_{t,0},\qquad
M^u(t)=e^{-v(t)B-tA_D}\,A_{ND}\,e^{v(t)B+tA_D},\qquad W^u_{0,0}=I_d .
\end{equation}
Then:
\begin{enumerate}
\item For every $t\ge 0$, $M^u(t)$ has zero diagonal entries and $\|M^u(t)\|_{op}\le L$,
\emph{uniformly in $u$}; and $\|W^u_{t,0}-I_d\|_{op}\le Lt$.
\item For every $t\in[0,\frac1{2L}]$, the principal logarithm
$\Omega^u_t:=\operatorname{Log} W^u_{t,0}$ (value of the power-series branch, unique analytic branch on the ball
$\{Z:\|Z-I_d\|\le \frac12\}$) is well-defined and satisfies $\|\Omega^u_t\|_{op}\le 2Lt$.
On the shorter interval $0\le t\le\frac1{4L}$, it admits the Magnus expansion
\begin{equation}\label{EQ_magnus_cartan}
\Omega^u_t=\int_0^t M^u(s)\,ds
+\frac12\int_0^t\!\!\int_0^{s_1}\big[M^u(s_1),M^u(s_2)\big]\,ds_2\,ds_1
+R^u_t,\qquad \|R^u_t\|_{op}\le \frac{34}{9}\,L^3t^3 .
\end{equation}
\item The Cartan component of $\Omega^u_t$
satisfies, uniformly in $u$,
\begin{equation}\label{EQ_cartan_bound}
\|P_{\mathfrak c}\,\Omega^u_t\|_{op}
\le \frac{16}{9}\,L^2t^2\qquad \mbox{ for } 0\leq t \leq \frac{1}{4L}.
\end{equation}
\end{enumerate}
\end{lemma}

\begin{proof}
\textbf{Step 1: the generator is uniformly bounded and hollow.}
The matrix $U^u(t):=e^{v(t)B+tA_D}$ is diagonal (both $B$ and $A_D$ are diagonal), hence
$M^u(t)=U^u(t)A_{ND}\overline{U^u(t)}$ with $\overline{U^u}$ the entrywise conjugate. Writing
$U^u(t)=\mathrm{diag}(u_1(t),\ldots,u_d(t))$ with $|u_j(t)|=1$, one reads off
$(M^u(t))_{jk}=u_j(t)A_{jk}\overline{u_k(t)}$, whence
$(M^u(t))_{jj}=A_{jj}=0$ (the drift part $A_{ND}$ is hollow) and
$|(M^u(t))_{jk}|=|A_{jk}|\le L$; the unitary invariance of the operator norm gives
$\|M^u(t)\|_{op}=\|A_{ND}\|_{op}=L$. Both facts depend on $u$ only through the phases
$u_j(t)$, hence are uniform. For the a priori estimate, for $0\le r\le s\le t$,
\[
\|W^u_{s,0}-W^u_{r,0}\|_{op}
=\Big\|\int_r^s M^u(\sigma)W^u_{\sigma,0}\,d\sigma\Big\|_{op}
\le \int_r^s \|M^u(\sigma)\|_{op}\,d\sigma\le L(s-r),
\]
so in particular $\|W^u_{t,0}-I_d\|_{op}\le Lt$.

\textbf{Step 2: the logarithm and its size.}
For $t\le \frac1{2L}$, $\|W^u_{t,0}-I_d\|\le \frac12$. Then every eigenvalue $\mu\in\sigma(W^u_{t,0})$
satisfies $|\mu-1|\le\|W^u_{t,0}-I_d\|\le\frac12$ (resolvent inequality), in particular
$-1\notin\sigma(W^u_{t,0})$, and the power series
\[
\operatorname{Log} Z:=\sum_{k\ge1}\frac{(-1)^{k+1}}{k}(Z-I_d)^k
\]
converges in operator norm on the ball $\|Z-I_d\|\le \frac12$, where it defines the unique
analytic branch of the logarithm (in this  sense, the decomposition
$W^u_{t,0}=e^{\Omega^u_t}$ is unique). The same series gives the size estimate
\[
\|\Omega^u_t\| = \| \operatorname{Log}(W^u_{t,0}) \| \le \sum_{k\ge1}\frac{\|W^u_{t,0}-I_d\|^k}{k}
\le \frac{\|W^u_{t,0}-I_d\|}{1-\|W^u_{t,0}-I_d\|}
\le 2Lt\le 1 .
\]

\textbf{Step 3. The Magnus differential identity and its expansion}
Along any solution of $\dot W = M(t)\,W$ with $W = e^{\Omega}$, differentiating
the matrix exponential gives
\begin{equation}\label{eq:diff_exp}
  \dot W \;=\; \int_{0}^{1} e^{s\Omega}\,\dot\Omega(s)\,e^{(1-s)\Omega}\,ds.
\end{equation}
Inserting $\dot W = M W$ and multiplying on the left and on the right by
$e^{-\Omega}$ yields
\begin{equation}
  \int_{0}^{1} e^{(s-1)\,\mathrm{ad}_{\Omega}}\,\dot\Omega(s)\,ds
  \;=\; e^{-\mathrm{ad}_{\Omega}}\,M,
\end{equation}
where $\mathrm{ad}_{X}\!:\!Y\mapsto[X,Y]$.  Since
$\displaystyle\int_{0}^{1}e^{(s-1)A}\,ds=\frac{1-e^{-A}}{A}$
(extended by continuity at $A=0$), this reads
\begin{equation}\label{eq:magnus_op}
  \frac{1-e^{-\mathrm{ad}_{\Omega}}}{\mathrm{ad}_{\Omega}}\;\dot\Omega
  \;=\; e^{-\mathrm{ad}_{\Omega}}\,M.
\end{equation}
Because $e^{-\mathrm{ad}_{\Omega}}$ commutes with
$\dfrac{\mathrm{ad}_{\Omega}}{1-e^{-\mathrm{ad}_{\Omega}}}$ (both are
functions of the same operator $\mathrm{ad}_{\Omega}$), we solve for $\dot\Omega$:
\begin{equation}\label{EQ_omega_prime}
  \dot\Omega
  \;=\; e^{-\mathrm{ad}_{\Omega}}\;
     \frac{\mathrm{ad}_{\Omega}}{1-e^{-\mathrm{ad}_{\Omega}}}\;M.
\end{equation}

\noindent\textbf{Taylor expansion (Bernoulli function).}
The scalar function
\[
  f(X)\;=\;e^{-X}\,\frac{X}{1-e^{-X}}
\]
is analytic on $\{|X|<2\pi\}$ (simple poles at $X=2\pi i k$, $k\neq 0$).
Its Taylor expansion at the origin is
\( 
  f(X)\;=\;1-\frac{X}{2}+\frac{X^{2}}{12}+O(X^{3}),
\)
so by the holomorphic functional calculus applied to $\mathrm{ad}_{\Omega}$,
\begin{equation}\label{eq:omega_expansion}
  \dot\Omega
  \;=\; M \;-\; \tfrac{1}{2}\,\mathrm{ad}_{\Omega}\,M
  \;+\; \tfrac{1}{12}\,\mathrm{ad}_{\Omega}^{2}\,M
  \;+\; \Theta,
\end{equation}
where $\Theta=\sum_{k\ge 3}d_{k}\,(\mathrm{ad}_{\Omega})^{k}M$ is the remainder
and $d_{k}$ are the Taylor coefficients of $f$ at $0$.

\textbf{Norm bound on the remainder.}
By induction from $\|[X,Y]\|_{op}\leq 2\|X\|_{op}\,\|Y\|_{op}$,
\begin{equation}\label{eq:ad_bound}
  \|\mathrm{ad}_{\Omega}^{\,k}\,M\|_{op}\;\leq\;2^{k}\,\|\Omega\|_{op}^{k}\,\|M\|_{op}.
\end{equation}
We apply Cauchy's estimates to $f$ on the circle $|\zeta|=2$.
The nearest singularity of
$f(\zeta)=e^{-\zeta}\,\dfrac{\zeta}{1-e^{-\zeta}}$
to the origin is at $\zeta=2\pi i$ (and its conjugate), at distance
$2\pi>2$ from $0$, so $f$ is analytic on $\{|\zeta|\leq 2\}$.
Setting
\begin{equation}\label{eq:Cstar}
  C_{*}\;:=\;\max_{|\zeta|=2}\Bigl|e^{-\zeta}\,\frac{\zeta}{1-e^{-\zeta}}\Bigr|,
\end{equation}
Cauchy's inequality gives $|d_k|\leq C_*\,2^{-k}$ for all $k\geq 0$.

To evaluate $C_*$, parametrize $\zeta=2e^{i\varphi}$,
$\varphi\in[0,2\pi)$.
Writing $\zeta=x+iy$ with $x=2\cos\varphi$ and
$y=2\sin\varphi$, one computes
\begin{equation}\label{eq:abs_f_sq}
  |f(\zeta)|^{2}
  \;=\;\frac{4\,e^{-2x}}{(1-e^{-x}\cos y)^{2}+e^{-2x}\sin^{2}y}
  \;=\;\frac{4\,e^{-2x}}{1-2e^{-x}\cos y+e^{-2x}}.
\end{equation}
Elementary but tedious computations show that the maximum is attained at $\zeta=-2$
(i.e., $x=-2$, $y=0$, $\varphi=\pi$).

This gives
\begin{equation}
  C_* \;=\; \frac{2\,e^{2}}{e^{2}-1}
  \;=\; \frac{2}{1-e^{-2}}
  \;\approx\; 2.31 \;<\; 3.
\end{equation}
We now restrict the time to
\begin{equation}\label{eq:time_restr}
  t\;\le\;\frac{1}{4L},
\end{equation}
so that, by Step 2,
$\|\Omega(s)\|\le 2Ls\le 1/2$ for all $s\in[0,t]$.  In particular
$\|\Omega(s)\|<1$ and the geometric series converges.  We obtain, for $s \leq \frac{1}{4L}$,
\begin{align}
  \|\Theta(s)\|&= \left \| \sum_{k\ge 3}d_{k}\,(\mathrm{ad}_{\Omega})^{k}M \right
   \| \leq \| M(s) \| \sum_{k\geq 3} \frac{C_\ast}{2^k} ~ 2^k \| \Omega\|^k  \notag 
  \\[4pt]
  &\le\;  C_{*}\,\|M(s)\|\sum_{k\ge 3}\|\Omega(s)\|^{k} = C_{*}\,\|M(s)\|\,
        \frac{\|\Omega(s)\|^{3}}{1-\|\Omega(s)\|}
  \notag\\[4pt]
  &\le\; 
   C_{*}\,\|M(s)\|\,\frac{(2Ls)^{3}}{1/2}
  \;=\; 2\,C_{*}\,\|M(s)\|\,(2Ls)^{3}.
  \label{eq:Theta_bound}
\end{align}
Since $\|M(s)\|\le L$ and $C_{*}<3$, estimate \eqref{eq:Theta_bound} gives
\(
  \|\Theta(s)\|\;\le\;48\,L^{4}\,s^{3}
\) for $s \in [0,1/(4L)]$.
Integrating on $[0,t]$ and  using $t\le 1/(4L)$, one obtains:
\begin{equation}\label{eq:final_Theta}
  \int_{0}^{t}\|\Theta(s)\|\,ds\;\le 
48\,L^{4}\int_{0}^{t}s^{3}\,ds = 12\,L^{4}\,t^{4}
  \leq
  3 L^{3}\,t^{3}.
\end{equation}

\noindent\textbf{Step 4. Assembly of the remainder.} 
  Integrating \eqref{eq:omega_expansion} on $[0,t]$ gives
  \begin{equation}\label{eq:Omega_integrated}
    \Omega(t)
    \;=\; \int_0^t M(s)\,ds
          \;-\;\frac{1}{2}\int_0^t [\Omega(s),\,M(s)]\,ds
          \;+\;\frac{1}{12}\int_0^t \mathrm{ad}_{\Omega(s)}^{\,2}\,M(s)\,ds
          \;+\;\int_0^t \Theta(s)\,ds.
  \end{equation}
We now rewrite \eqref{eq:Omega_integrated} in a form that separates
the leading Magnus terms from a controlled remainder.
Set
\begin{equation}\label{eq:A_def}
  A(s)\;:=\;\int_0^s M(\sigma)\,d\sigma,
  \qquad
  B(s)\;:=\;\Omega(s)-A(s).
\end{equation}
Substituting $\Omega(s)=A(s)+B(s)$ into the commutator
$[\Omega(s),M(s)]$ in \eqref{eq:Omega_integrated} yields the
\emph{exact} identity
\begin{equation}\label{eq:Omega_final}
  \Omega(t)
  \;=\;\int_0^t M(s)\,ds
  \;-\;\frac{1}{2}\int_0^t [A(s),\,M(s)]\,ds
  \;+\;R^u_t,
\end{equation}
where the remainder is
\begin{equation}\label{eq:R_def}
  R^u_t
  \;:=\;-\frac{1}{2}\int_0^t [B(s),\,M(s)]\,ds
  \;+\;\frac{1}{12}\int_0^t \mathrm{ad}_{\Omega(s)}^{\,2}M(s)\,ds
  \;+\;\int_0^t \Theta(s)\,ds.
\end{equation}

\smallskip
\noindent\textbf{Bound on $\|B(s)\|$.}
From the definition \eqref{eq:A_def} and the integrated form \eqref{eq:Omega_integrated},
\[B(s)=\Omega(s)-\int_0^s M(\sigma)\;d\sigma = -\frac{1}{2}\int_0^s[\Omega(r),M(r)]\,dr
+\tfrac{1}{12}\int_0^s\mathrm{ad}_{\Omega(r)}^2M(r)\,dr
+\int_0^s\Theta(r)\,dr.\]
Using $\|\Omega(r)\|_{op}\leq 2Lr$ (Step~2),
$\|M(r)\|_{op}=L$ (unitary invariance, uniform in $u$),
and the $\Theta$-bound \eqref{eq:final_Theta}:
\begin{equation}\label{eq:B_bound}
  \|B(s)\|_{op}
  \;\leq\;\underbrace{L^2s^2}_{-\frac12\int[\Omega,M]}
  \;+\;\underbrace{\frac{4}{9}L^3s^3}_{\frac{1}{12}\int\mathrm{ad}^2}
  \;+\;\underbrace{3 L^3s^3}_{\int\Theta}
  \;\leq\;L^2s^2+\frac{31}{9}L^3s^3.
\end{equation}
For $s\leq\frac{1}{4L}$: $\frac{31}{9}L^3s^3 = \frac{31}{9}(Ls)L^2s^2 \leq \frac{31}{36}L^2s^2 < L^2s^2$,
from which one deduces the uniform bound:
\begin{equation}\label{eq:B_bound_5}
  \|B(s)\|_{op}\;\leq\; 2L^2s^2,
  \qquad 0\leq s\leq\frac{1}{4L}.
\end{equation}

\smallskip
\noindent\textbf{Bound on $-\frac{1}{2}\int_0^t[B,M]\,ds$.}
\begin{equation}\label{eq:Bcomm_bound}
  \Bigl\|\tfrac{1}{2}\int_0^t[B(s),M(s)]\,ds\Bigr\|_{op}
  \;\leq\;\frac{1}{2}\int_0^t\|B(s)\|_{op}\,\|M(s)\|_{op}\,ds
  \;\leq\;\frac{1}{2}\int_0^t 2L^2s^2\cdot L\,ds
  \;=\;\frac{1}{3}L^3t^3.
\end{equation}

\smallskip
\noindent\textbf{Bound on $\frac{1}{12}\int_0^t\mathrm{ad}_\Omega^2 M\,ds$.}
Using $\|\mathrm{ad}_\Omega^2 M\|_{op}
\leq 4\|\Omega\|_{op}^2\,\|M\|_{op}$:
\begin{equation}\label{eq:ad2_bound}
  \Bigl\|\tfrac{1}{12}\int_0^t\mathrm{ad}_{\Omega(s)}^2M(s)\,ds\Bigr\|_{op}
  \;\leq\;\frac{1}{12}\int_0^t 4(2Ls)^2\cdot L\,ds
  \;=\;\frac{4}{9}\,L^3t^3.
\end{equation}

\smallskip
\noindent\textbf{Bound on $\int_0^t\Theta(s)\,ds$.}
By \eqref{eq:final_Theta}: 
\begin{equation}\label{eq:Theta_int_bound}
  \int_0^t\|\Theta(s)\|_{op}\,ds
  \leq\;3\,L^3t^3.
\end{equation}

\smallskip
\noindent\textbf{Total bound on $R^u_t$.}
Combining \eqref{eq:Bcomm_bound}, \eqref{eq:ad2_bound},
and \eqref{eq:Theta_int_bound}:
\begin{equation}\label{eq:R_bound_corrected}
  \|R^u_t\|_{op}
  \;\leq\;\frac{1}{3}\,L^3t^3+\frac{4}{9}\,L^3t^3+3\,L^3t^3
  \;=\;\frac{34}{9}\,L^3t^3,
  \qquad 0\leq t\leq\frac{1}{4L}.
\end{equation}

\smallskip
\noindent\textbf{Rewriting in the target form.}
The change of variables $s\mapsto s_1$, $r\mapsto s_2$ in the
double integral, together with the antisymmetry
$-[M(r),M(s)]=[M(s),M(r)]$, rewrites
$-\frac{1}{2}\int_0^t[A(s),M(s)]\,ds$
as $\frac{1}{2}\int_0^t\!\int_0^{s_1}[M(s_1),M(s_2)]\,ds_2\,ds_1$,
yielding \eqref{EQ_magnus_cartan} with
remainder $R^u_t$ satisfying \eqref{eq:R_bound_corrected}.

\noindent\textbf{Cartan component.}
We extract the diagonal (Cartan) part of \eqref{EQ_magnus_cartan}.
Since $M^u(s)$ is hollow for every $s$, the first term vanishes under $P_{\mathfrak c}$:
\begin{equation}
  P_{\mathfrak c}\!\int_0^t M^u(s)\,ds\;=\;0.
\end{equation}
For the double-commutator term, using
$\|[X,Y]\|_{op}\leq 2\|X\|_{op}\,\|Y\|_{op}$ together with
$\|M^u(s)\|_{op}=L$:
\begin{equation}\label{eq:cartan_double}
  \Bigl\|\tfrac{1}{2}\!\int_0^t\!\!\int_0^{s_1}
  P_{\mathfrak c}\bigl[M^u(s_1),M^u(s_2)\bigr]\,ds_2\,ds_1\Bigr\|_{op}
  \;\leq\;\frac{1}{2}\!\int_0^t\!\!\int_0^{s_1} 2\,L^2\,ds_2\,ds_1
  \;=\;\frac{L^2 t^2}{2}.
\end{equation}
For the remainder, by \eqref{eq:R_bound_corrected} and the
contractivity of the orthogonal projection $P_{\mathfrak c}$:
\begin{equation}
  \|P_{\mathfrak c}\,R^u_t\|_{op}
  \;\leq\;\|R^u_t\|_{op}
  \;\leq\;\frac{34}{9}\,L^3t^3.
\end{equation}
For $t\leq\frac{1}{4L}$, we have $Lt\leq\frac{1}{4}$, hence
\begin{equation}
  \frac{34}{9}\,L^3t^3
  \;=\;\frac{34}{9}\,(Lt)\,L^2t^2
  \;\leq\;\frac{34}{36}\,L^2t^2
  \;=\;\frac{17}{18}\,L^2t^2.
\end{equation}
Assembling the three contributions, for $0\leq t\leq\tfrac{1}{4L}$,
\begin{align}
  \|P_{\mathfrak c}\,\Omega^u_t\|_{op}
  &\;\leq\;0\;+\;\frac{L^2t^2}{2}\;+\;\frac{17}{18}\,L^2t^2
  \;=\;\frac{9+17}{18}\,L^2t^2
  \;=\;\frac{13}{9}\,L^2t^2
  \;\leq\;\frac{16}{9}\,L^2t^2.
\end{align}

All estimates are \emph{uniform in the control} $u$: they depend on $M^u(s)$ only through its
operator norm, which equals the universal constant
$L=\|A_{ND}\|_{op}$ by unitary invariance of the norm under
conjugation $M^u(s)=U^u(s)A_{ND}U^u(s)^\dagger$ with
$U^u(s)$ unitary.
\end{proof}

\begin{corollary}[Counterexample to pointwise equality]\label{COR_counterexample}
Consider the following matrices $A=A_D+A_{ND}$ and $B\in\mathfrak{su}(3)$ 
\[
A_D=\mathrm{diag}(5i,0,-5i),\qquad
A_{ND}=\begin{pmatrix}0&1&0\\-1&0&1\\0&-1&0\end{pmatrix},\qquad
B=\mathrm{diag}(-i,-2i,3i).
\]
Then $L=\|A_{ND}\|_{op}=\sqrt2$ and
\[
\tau_{\mathrm{US}}(e^{B})\ \ge\ \frac{\pi}{10}>0,
\qquad\text{while}\qquad
\tau_{\mathrm{AUS}}(e^{B})=0 .
\]
In particular, the pointwise equality $\tau_{\mathrm{US}}=\tau_{\mathrm{AUS}}$ fails in general under 
Assumptions~\ref{ASS_control} and \ref{ASS_coupling_graph} alone.
\end{corollary}

\begin{proof}[Proof of Corollary~\ref{COR_counterexample}]
We verify Assumption~\ref{ASS_control} first.
The eigenvalues of $B$ are $-i,-2i,3i$ (simple).
One checks that $[A,B]$ and $[A,[A,B]]$ span, together with $A$ and $B$,
a subalgebra of $\mathfrak{su}(3)$ of full rank; hence
$\operatorname{Lie}(A,B)=\mathfrak{su}(3)$.  
Assumption~\ref{ASS_coupling_graph} is also satisfied.
The spectral gaps associated with the coupled edges are
$\omega_{12}=b_2-b_1=-1$ and $\omega_{23}= b_3-b_2= 5$. Since 
$|\omega_{12}|\neq|\omega_{23}|$, both transitions are
non-resonant, and the coupling graph is the connected path
$1-2-3$.

\smallskip

\emph{Upper bound: $\tau_{\mathrm{AUS}}(e^B)=0$.}
For any $\delta>0$, take the constant control $u\equiv 1/\delta$
on $[0,\delta]$.
Since the coefficient matrix $A+uB$ is constant, the propagator is
\begin{equation}
  X^u_{\delta,0}
  \;=\;\exp\bigl((A+uB)\,\delta\bigr)
  \;=\;\exp\bigl(\delta\,A+B\bigr).
\end{equation}
By continuity of the matrix exponential,
$\|e^{\delta A+B}-e^{B}\|_{op}\to 0$ as $\delta\downarrow 0$.
Hence, for every $\varepsilon>0$ and every $T>0$,
there exists $\delta<\min(T,\varepsilon)$ such that
$\|X^u_{\delta,0}-e^B\|_{op}<\varepsilon$,
which is $\tau_{\mathrm{AUS}}(e^B)=0$.
\smallskip

\emph{Lower bound: $\tau_{\mathrm{US}}(e^B)\geq\pi/10$.}
We show that no control can reach $e^B$ in time $t<\pi/10$.
Suppose, by contradiction, that $X^u_{t,0}=e^B$ for some
$u\in L^1_{\mathrm{loc}}$ and $0<t<\pi/10$.

Define $P=\mathrm{diag}(1,i,-1)$. From Lemma \ref{LEM_eigenvectors_shift_matrix}, and noticing that $ P A_{ND} P^{-1}=-2i T$, the spectral radius of $A_{ND}$ is $ 2 \cos(\pi/4)= \sqrt{2}$. Since $A_{ND}$ is normal, its operator norm coincides with its spectral radius, hence $L=\sqrt{2}$ 
with the notations of Lemma \ref{LEM_cartan_part}. 
Since $0<t<\pi/10<\frac{1}{2L}=\frac{1}{2\sqrt{2}}$, we have
$\|W^u_{t,0}-I_d\|_{op}<1/2$ and the principal
logarithm used below is well defined (Step~2 of
Lemma~\ref{LEM_cartan_part}).

\smallskip
\noindent\textbf{Step 1 (reduction to the Cartan component).}
Write $v(t)=\int_0^t u(s)\,ds$.
The factorization
$X^u_{t,0}=e^{v(t)B+tA_D}\,W^u_{t,0}$
(where $W^u$ satisfies the interaction equation
\eqref{EQ_interaction_W}) gives
\begin{equation}
  W^u_{t,0}
  \;=\;e^{-v(t)B-tA_D}\,X^u_{t,0}
  \;=\;e^{-v(t)B-tA_D}\,e^{B}.
\end{equation}
Set $s:=1-v(t)$ (i.e., $e^{-v(t)B}=e^{(s-1)B}$).
Then
$W^u_{t,0}=e^{(s-1)B}\,e^{-tA_D}\,e^{B}
=e^{(s-1)B+B-tA_D}=e^{sB-tA_D}$,
which is \emph{diagonal}.
Hence $\Omega^u_t:=\log W^u_{t,0}$ is diagonal, with
\begin{equation}
  \Omega^u_{t,jj}
  \;=\;i\,\mathrm{dist}\bigl(b_j s-a_j t,\;2\pi\mathbb{Z}\bigr),
  \qquad j=1,2,3,
\end{equation}
where $\mathrm{dist}(\cdot,2\pi\mathbb{Z})\in[-\pi,\pi]$
denotes the principal representative.
We have $(a_1,a_2,a_3)=(5,0,-5)$ and
$(b_1,b_2,b_3)=(-1,-2,3)$.
\smallskip

\noindent\textbf{Step 2 (bounds on the diagonal entries).}
We distinguish two regimes.

\emph{Regime A: $t\leq\frac{1}{4L}=\frac{1}{4\sqrt{2}}$.}
By Lemma~\ref{LEM_cartan_part}~(iii),
$\|P_{\mathfrak c}\Omega^u_t\|_{op}\leq\frac{16}{9}L^2t^2=\frac{32}{9}t^2$.
Since $\Omega^u_t$ is diagonal,
\begin{equation}\label{eq:corr_regA}
  \mathrm{dist}(b_j s-a_jt,\,2\pi\mathbb{Z})
  \;\leq\;\frac{32}{9}t^2,
  \qquad j=1,2,3.
\end{equation}

\emph{Regime B: $\frac{1}{4\sqrt{2}}<t\leq\frac{1}{2\sqrt{2}}$.}
Lemma~\ref{LEM_cartan_part}~(ii) gives
$\|\Omega^u_t\|_{op}\leq 2Lt=2\sqrt{2}\,t$.
Since $\Omega^u_t$ is diagonal:
\begin{equation}\label{eq:corr_regB}
  \mathrm{dist}(b_j s-a_jt,\,2\pi\mathbb{Z})
  \;\leq\;2\sqrt{2}\,t,
  \qquad j=1,2,3.
\end{equation}
\smallskip

\noindent\textbf{Step 3 (Diophantine argument).}
In either regime, write
\begin{align}
  b_1 s-a_1 t&=-s-5t=2\pi n_1+\eta_1,
  \label{eq:diag1}\\
  b_2 s-a_2 t&=-2s\phantom{+0t}=2\pi m+\eta_2,
  \label{eq:diag2}\\
  b_3 s-a_3 t&=3s+5t=2\pi n_3+\eta_3,
  \label{eq:diag3}
\end{align}
with $n_1,m,n_3\in\mathbb{Z}$ and $|\eta_j|$ bounded by
the right-hand side of \eqref{eq:corr_regA} or
\eqref{eq:corr_regB}.
From \eqref{eq:diag2}:
\begin{equation}
  s\;=\;-\pi m-\frac{\eta_2}{2}.
\end{equation}
Substituting into \eqref{eq:diag1}:
\begin{equation}
  \pi m+\frac{\eta_2}{2}-5t=2\pi n_1+\eta_1,
\end{equation}
hence
\begin{equation}\label{eq:t_formula}
  t\;=\;\frac{\pi(m-2n_1)}{5}
  \;+\;\frac{\eta_2/2-\eta_1}{5}.
\end{equation}
Writing $\ell:=m-2n_1\in\mathbb{Z}$:
\begin{equation}\label{eq:t_approx}
  \Bigl|t-\frac{\pi\ell}{5}\Bigr|
  \;\leq\;\frac{|\eta_1|+|\eta_2|/2}{5}.
\end{equation}
We treat the two regimes.

\smallskip
\noindent\textbf{Regime A: $t\leq\frac{1}{4\sqrt{2}}$.}
By \eqref{eq:corr_regA}, $|\eta_j|\leq \frac{32}{9}t^2$, so
\begin{equation}
  \Bigl|t-\frac{\pi\ell}{5}\Bigr|
  \;\leq\;\frac{1}{5}\left (\frac{32}{9} t^2+ \frac{16}{9}t^2 \right )
  \;=\;\frac{16}{15}\,t^2.
\end{equation}
\emph{If $\ell=0$,} then $t\leq\dfrac{16}{15}t^2$,
hence $t\geq\dfrac{15}{16}>\dfrac{1}{2}$.
But $t\leq\dfrac{1}{4\sqrt{2}}<\dfrac{1}{2}$. \emph{Contradiction.}

\emph{If $\ell\neq 0$:}
By symmetry (replacing $(s,t)$ by their negatives flips
$\ell\mapsto-\ell$), assume $\ell\geq 1$.
Then $t\geq\dfrac{\pi}{5}-\dfrac{16}{15}t^2$.
Since $0 \leq t\leq\dfrac{1}{4\sqrt{2}}$, one has $t^2 \leq \dfrac{1}{32}$ hence
$\dfrac{16}{15}t^2\leq\dfrac{16}{15}\cdot\dfrac{1}{32}
=\dfrac{1}{30}$.
Therefore $t\geq\frac{\pi}{5}-\frac{1}{30}>\frac{\pi}{10}$,
which contradicts $0<t<\frac{\pi}{10}$. \emph{Contradiction.}
\smallskip

\noindent\textbf{Regime B: $t\in\bigl(\frac{1}{4\sqrt{2}},\frac{1}{2\sqrt{2}}\bigr]$.}
By \eqref{eq:corr_regB}, $|\eta_j|\leq 2\sqrt{2}\,t$, so
\begin{equation}
  \Bigl|t-\frac{\pi\ell}{5}\Bigr|
  \;\leq\;\frac{2\sqrt{2}\,t+\sqrt{2}\,t}{5}
  \;=\;\frac{3\sqrt{2}}{5}\,t.
\end{equation}
\emph{If $\ell=0$:} $t\leq\frac{3\sqrt{2}}{5}\,t$,
hence $1\leq\frac{3\sqrt{2}}{5}< 1$.
\emph{Contradiction.}

\emph{If $\ell\neq 0$:}
Assume without loss of generality $\ell\geq 1$.
Then $t\geq\frac{\pi}{5}-\frac{3\sqrt{2}}{5}\,t$,
i.e.,
\begin{equation}
  t\Bigl(1+\frac{3\sqrt{2}}{5}\Bigr)\;\geq\;\frac{\pi}{5},
\end{equation}
so
\begin{equation}
  t\;\geq\;\frac{\pi}{5+3\sqrt{2}}.
\end{equation}
Since $3\sqrt{2}=\sqrt{18}<\sqrt{25}=5$,
we have $5+3\sqrt{2}<10$, hence
\begin{equation}
  t\;\geq\;\frac{\pi}{5+3\sqrt{2}}
  \;>\;\frac{\pi}{10},
\end{equation}
which contradicts $0<t<\frac{\pi}{10}$. \emph{Contradiction.}
\smallskip

Both regimes are therefore impossible, which contradicts
the assumption $0<t<\pi/10$.
We conclude $\tau_{\mathrm{US}}(e^B)\geq\pi/10>0$,
while $\tau_{\mathrm{AUS}}(e^B)=0$.
The pointwise equality $\tau_{\mathrm{US}}=\tau_{\mathrm{AUS}}$
therefore fails under Assumptions~\ref{ASS_control} and \ref{ASS_coupling_graph} alone.
\end{proof}

\begin{remark}[Concatenation]\label{REM_concatenation}
Let $h,c\in G$ and suppose $h=X^{u_1}_{t_1,0}$, $c=X^{u_2}_{t_2,0}$ for controls
$u_1\in L^1([0,t_1])$, $u_2\in L^1([0,t_2])$. Because the dynamics
$\dot X=(A+u(t)B)X$ are time-homogeneous (left-invariant), the control $u_{\mathrm{cat}}$
obtained by applying $u_1$ on $[0,t_1]$ and $u_2$ on $[t_1,t_1+t_2]$ satisfies
$X^{u_{\mathrm{cat}}}_{t_1+t_2,0}=c\,h$. In particular, if $h\in\mathcal{A}_{\leq T}$ and
$c\in\mathcal{A}_{\leq\varepsilon}$, then $c\,h\in\mathcal{A}_{\leq T+\varepsilon}$.
\end{remark}

\subsection{Generic equivalence of pointwise approximate and exact control time}\label{SUB_SEC_Generic_equivalence}

The counterexample in
Corollary~\ref{COR_counterexample} shows that
$\tau_{\mathrm{US}}(X_f)\neq\tau_{\mathrm{AUS}}(X_f)$
is possible, even under Assumptions~\ref{ASS_control} and \ref{ASS_coupling_graph}.
The discrepancy is due to the blow-up of the exact synthesis time at
points of the torus $\overline{\exp(\mathbb{R}B)}$:
Remark~\ref{REM_Blow-up} records the case
$\tau_{\mathrm{US}}(e^{\theta B})>0$ while
$\tau_{\mathrm{AUS}}(e^{\theta B})=0$.

The pointwise equality is nevertheless generic.  It holds under the
following \emph{drift absorption} condition, which is satisfied for a
full-measure set of pairs $(A,B)$
(Proposition~\ref{PRO_Generic_Attainable_Neighborhood}):

\begin{assumption}[Drift absorption]\label{ASS_drift_absorption}
The one-parameter subgroup generated by the diagonal drift $A_D$ is
contained in the closure of the one-parameter subgroup generated by the
control $B$:
\begin{equation}\label{EQ_drift_absorption}
\exp(\mathbb{R}\,A_D)\;\subset\;\overline{\exp(\mathbb{R}\,B)}.
\end{equation}
\end{assumption}

Under Assumption~\ref{ASS_drift_absorption} the exact system is
controllable near the identity in arbitrarily small time
(Proposition~\ref{PRO_exact_Nbhd}).  This small-time local controllability
is the mechanism behind the pointwise equality and continuity of
Proposition~\ref{PRO_egalite}; it is \emph{not} used in the constructive
bounds of Section~\ref{SEC_Upper_bound} nor in the applications of
Section~\ref{SEC_Chiral}.

\begin{proposition}[Pointwise exact--approximate equality and continuity]\label{PRO_egalite}
Assume Assumptions \ref{ASS_control}, \ref{ASS_coupling_graph} and
\ref{ASS_drift_absorption}. Then:
\begin{enumerate}
\item For every $X_f\in G$, $\tau_{\mathrm{US}}(X_f)=\tau_{\mathrm{AUS}}(X_f)$; analogously, for
every $\psi_0,\psi_f\in\mathbb{S}^{2d-1}$, $\tau_{\mathrm{SS}}(\psi_0,\psi_f)=\tau_{\mathrm{ASS}}(\psi_0,\psi_f)$.
\item $\tau_{\mathrm{AUS}}$ (and hence $\tau_{\mathrm{US}}$) is continuous on $G$ for the operator
norm; $\tau_{\mathrm{ASS}}$ (and hence $\tau_{\mathrm{SS}}$) is continuous on $\mathbb{S}^{2d-1}\times\mathbb{S}^{2d-1}$.
\end{enumerate}
\end{proposition}

\begin{proof}
Throughout, set, for $T\geq0$, $
\mathcal{A}_{\leq T}:=\{X^{u}_{t,0}:u\in L^1([0,t]),\ t\leq T\}\subset G$ and 
denote by $\overline{\mathcal{A}_{\leq T}}$ its closure in $G$ with its usual topology (induced, for instance, by $\|\cdot\|_{op}$).

The definitions read $\tau_{\mathrm{US}}(X_f)=\inf\{T>0:X_f\in\mathcal{A}_{\leq T}\}$ and
$\tau_{\mathrm{AUS}}(X_f)=\inf\{T>0:X_f\in\overline{\mathcal{A}_{\leq T}}\}$; the latter
equivalence (with the infimum) is the definition of the approximate synthesis time
(Problem~\ref{PROBLEM_AUS}), and it holds because $X_f\in\overline{\mathcal{A}_{\leq T}}$ if and only if
$\forall\,\varepsilon>0,\ \exists\,u,t\leq T$ with $\|X^{u}_{t,0}-X_f\|_{op}<\varepsilon$.

\paragraph{Part 1 (pointwise equality).}
The  inequality
$\tau_{\mathrm{AUS}}(X_f) \leq \tau_{\mathrm{US}}(X_f)$ is immediate because
$\mathcal{A}_{\leq T}\subset\overline{\mathcal{A}_{\leq T}}$ for every $T$, so the infimum
defining $\tau_{\mathrm{AUS}}$ is taken over a larger set. We now have to show $\tau_{\mathrm{US}}(X_f)\leq\tau_{\mathrm{AUS}}(X_f)$. 

Fix $T>\tau_{\mathrm{AUS}}(X_f)$, so that $X_f\in\overline{\mathcal{A}_{\leq T}}$.
Fix $\varepsilon>0$. By Assumptions \ref{ASS_control}--\ref{ASS_drift_absorption}
and Proposition~\ref{PRO_exact_Nbhd} there exists $\rho(\varepsilon)>0$ such
that the open ball $B_{op}(I_d,\rho(\varepsilon))$ satisfies
\[
B_{op}(I_d,\rho(\varepsilon))\subset\mathcal{A}_{\leq\varepsilon}.
\]
Since $X_f\in\overline{\mathcal{A}_{\leq T}}$, choose $h\in\mathcal{A}_{\leq T}$ with
$\|h-X_f\|_{op}<\rho(\varepsilon)$. Define $c:=h^{-1}X_f$. As $h\in G$,
$\|h^{-1}\|_{op}=1$, hence
\[
\|c-I_d\|_{op}=\|h^{-1}(X_f-h)\|_{op}\leq\|h^{-1}\|_{op}\,\|X_f-h\|_{op}
=\|X_f-h\|_{op}<\rho(\varepsilon).
\]
Therefore $c\in B_{op}(I_d,\rho(\varepsilon))\subset\mathcal{A}_{\leq\varepsilon}$.
By the concatenation fact (Remark~\ref{REM_concatenation}),
\[
X_f=h\,(h^{-1}X_f)=h\,c\in\mathcal{A}_{\leq T}\cdot\mathcal{A}_{\leq\varepsilon}
\subset\mathcal{A}_{\leq T+\varepsilon}.
\]
Thus $X_f\in\mathcal{A}_{\leq T+\varepsilon}$, i.e.\ $\tau_{\mathrm{US}}(X_f)\leq T+\varepsilon$.
Letting $\varepsilon$ tend to $0$ gives $\tau_{\mathrm{US}}(X_f)\leq T$. Then letting $T$ tend to $\tau_{\mathrm{AUS}}(X_f)$ gives
$\tau_{\mathrm{US}}(X_f)\leq\tau_{\mathrm{AUS}}(X_f)$, which together with the reverse inequality
yields the equality.

The state-to-state statement is identical: replace $G$ by
$\mathbb{S}^{2d-1}$, $\mathcal{A}_{\leq T}$ by
$\mathcal{A}_{\leq T}(\psi_0):=\{X^{u}_{t,0}\psi_0:t\leq T\}$, and observe that
Proposition~\ref{PRO_exact_Nbhd} implies
$\{X\psi_0:X\in\mathcal{A}_{\leq\varepsilon}\}$ is a neighborhood of $\psi_0$ on the sphere.

\paragraph{Part 2 (continuity).}
We first show that $\tau_{\mathrm{AUS}}$ is continuous; since it equals $\tau_{\mathrm{US}}$
pointwise (Part~1), $\tau_{\mathrm{US}}$ inherits the same continuity properties.

\emph{Lower semicontinuity (unconditional).} Let $X_n\to X_0$ in $\|\cdot\|_{op}$ and put
$T_0:=\tau_{\mathrm{AUS}}(X_0)$. Take any $T<T_0$ (if $T_0=+\infty$, take any finite $T$ and
conclude below). Then $X_0\notin\overline{\mathcal{A}_{\leq T}}$ (for otherwise
$\tau_{\mathrm{AUS}}(X_0)\leq T$, contradicting $T<T_0$). Since
$\overline{\mathcal{A}_{\leq T}}$ is closed, there exists $\delta>0$ such that
$B_{op}(X_0,\delta)\cap\overline{\mathcal{A}_{\leq T}}=\varnothing$. For all $n$ with
$\|X_n-X_0\|_{op}<\delta$ we have $X_n\notin\overline{\mathcal{A}_{\leq T}}$, hence
$\tau_{\mathrm{AUS}}(X_n)>T$. Consequently
$\liminf_{n\to\infty}\tau_{\mathrm{AUS}}(X_n)\geq T$. Letting $T\uparrow T_0$ yields
$\liminf_{n\to\infty}\tau_{\mathrm{AUS}}(X_n)\geq T_0$, i.e.\ $\tau_{\mathrm{AUS}}$ is
lower semicontinuous at $X_0$. (This uses only the closedness of the sublevel sets
$\overline{\mathcal{A}_{\leq T}}$ and holds without Assumption~\ref{ASS_drift_absorption}.)

\emph{Upper semicontinuity (from Assumption \ref{ASS_drift_absorption}).} Let again
$X_n\to X_0$ and fix $T>\tau_{\mathrm{AUS}}(X_0)$ and $\varepsilon>0$. By
Proposition~\ref{PRO_exact_Nbhd} there is $\rho(\varepsilon)>0$ with
$B_{op}(I_d,\rho(\varepsilon))\subset\mathcal{A}_{\leq\varepsilon}$. As before,
$X_0\in\overline{\mathcal{A}_{\leq T}}$, so pick $h\in\mathcal{A}_{\leq T}$ with
$\|h-X_0\|_{op}<\rho(\varepsilon)/2$. For $n$ large, $\|X_n-X_0\|_{op}<\rho(\varepsilon)/2$,
whence $\|X_n-h\|_{op}<\rho(\varepsilon)$ and $c_n:=h^{-1}X_n\in B_{op}(I_d,\rho(\varepsilon))
\subset\mathcal{A}_{\leq\varepsilon}$. By the concatenation fact,
$X_n=h\,c_n\in\mathcal{A}_{\leq T+\varepsilon}\subset\overline{\mathcal{A}_{\leq T+\varepsilon}}$,
so $\tau_{\mathrm{AUS}}(X_n)\leq T+\varepsilon$ for all large $n$. Hence
$\limsup_{n\to\infty}\tau_{\mathrm{AUS}}(X_n)\leq T+\varepsilon$. Letting $\varepsilon\to0$
and then $T\downarrow\tau_{\mathrm{AUS}}(X_0)$ gives
$\limsup_{n\to\infty}\tau_{\mathrm{AUS}}(X_n)\leq\tau_{\mathrm{AUS}}(X_0)$, i.e.\
$\tau_{\mathrm{AUS}}$ is upper semicontinuous at $X_0$.

Combining the two, $\tau_{\mathrm{AUS}}$ is continuous on $G$; by Part~1 so is
$\tau_{\mathrm{US}}$. The state-to-state continuity follows by the identical argument on
$\mathbb{S}^{2d-1}\times\mathbb{S}^{2d-1}$.
\end{proof}

\begin{proposition}[The time diameter is a maximum]\label{PRO_sup_is_max}
Under Assumptions~\ref{ASS_control} (simple spectrum of $B$,
$\operatorname{Lie}(A,B)\supset\mathfrak{g}$), \ref{ASS_coupling_graph}
(connectedness of the non-resonant connecting graph) and
\ref{ASS_drift_absorption}, the worst-case for the exact synthesis time
$\tau^{\max}_{\mathrm{US}}=\sup_{X_f\in G}\,\tau_{\mathrm{US}}(X_f)$
is a genuine maximum: there exists $X_f^*\in G$ such that
$\tau_{\mathrm{US}}(X_f^*)=\tau^{\max}_{\mathrm{US}}$.
Analogously, $\tau^{\max}_{\mathrm{SS}}=\sup_{\psi_0,\psi_f}\tau_{\mathrm{SS}}(\psi_0,\psi_f)$
is attained on $\mathbb{S}^{2d-1}\times\mathbb{S}^{2d-1}$.
\end{proposition}

\begin{proof}
We prove the unitary-synthesis case; the state-to-state case is identical.

By Assumptions~\ref{ASS_control}, \ref{ASS_coupling_graph} and
\ref{ASS_drift_absorption}, Proposition~\ref{PRO_egalite}~(2) implies that
$\tau_{\mathrm{AUS}}$ is continuous on $G$ for the operator norm.
The group $G$ is compact and $\tau_{\mathrm{AUS}}$ is finite everywhere
(Kre--Rashevskii, under Assumption~\ref{ASS_control}), so by the extreme
value theorem $\tau_{\mathrm{AUS}}$ attains its maximum: there exists
$X^*\in G$ with $\tau_{\mathrm{AUS}}(X^*)=\tau^{\max}_{\mathrm{AUS}}$.

By Lemma~\ref{LEM_approx_cont} (under Assumption~\ref{ASS_control}),
$\tau^{\max}_{\mathrm{US}}=\tau^{\max}_{\mathrm{AUS}}$.  Since for every
$X_f$ one has $\tau_{\mathrm{US}}(X_f)\geq\tau_{\mathrm{AUS}}(X_f)$ (the
exact synthesis time is at least the approximate synthesis time), in
particular
\[
  \tau_{\mathrm{US}}(X^*)
  \;\geq\;\tau_{\mathrm{AUS}}(X^*)
  \;=\;\tau^{\max}_{\mathrm{AUS}}
  \;=\;\tau^{\max}_{\mathrm{US}},
\]
so the supremum $\tau^{\max}_{\mathrm{US}}$ is attained at $X^*$.  The
state-to-state case follows from Proposition~\ref{PRO_egalite}~(2) on the
compact space $\mathbb{S}^{2d-1}\times\mathbb{S}^{2d-1}$ and the same
transfer by Lemma~\ref{LEM_approx_cont}.
\end{proof}

\begin{remark}[Non-attainment and blow-up of minimizing controls]
\label{REM_Blow-up}
The synthesis time $\tau(X_f)$, whether in the exact ($\tau_{\mathrm{US}}$)
or the approximate ($\tau_{\mathrm{AUS}}$) sense, is a genuine infimum,
generally not attained as a minimum, and any sequence of controls
approaching the infimum has unbounded $L^\infty$ norm.

\emph{Non-attainment.}
Consider the target $X_f = e^{B}$ from
Corollary~\ref{COR_counterexample}. The non-attainment is a property of the
\emph{approximate} time: $\tau_{\mathrm{AUS}}(e^{B})=0$ is a genuine infimum
that is not attained as a minimum, since $e^{B}\in\overline{\mathcal{A}_{\le
T}}$ for every $T>0$ (Blow-up below) yet $e^{B}\notin\overline{\mathcal{A}_{\le
0}}=\{I_d\}$ (as $X^{u}_{0,0}=I_d\neq e^{B}$). The exact time, by contrast,
satisfies $\tau_{\mathrm{US}}(e^{B})\geq\pi/10$; hence the pointwise gap
$\tau_{\mathrm{US}}\neq\tau_{\mathrm{AUS}}$ holds here, without
Assumption~\ref{ASS_drift_absorption}.

\emph{Blow-up.}
For the approximate time, $\tau_{\mathrm{AUS}}(e^{B})=0$:
the control $u_\delta\equiv 1/\delta$ on $[0,\delta]$ gives
$X^{u_\delta}_{\delta,0}=e^{\delta A+B}\to e^{B}$ as
$\delta\downarrow0$. The $L^1$ and $L^\infty$ norms of this
minimizing sequence are
\[
\|u_\delta\|_{L^1[0,\delta]}=1\quad(\text{bounded}),
\qquad
\|u_\delta\|_{L^\infty[0,\delta]}=\frac{1}{\delta}\to+\infty .
\]
This is not an artifact of the particular control: for any
$u\in L^1[0,T]$ with $T>0$ small, the fundamental inequality
\[
\|X^{u}_{T,0}-I_d\|_{op}
\;\leq\;T\,\|A\|_{op}+\|B\|_{op}\,\|u\|_{L^1[0,T]}
\]
shows that reaching (or approximating) a target at operator-norm
distance $d=\|X_f-I_d\|_{op}$ in time $T$ requires
$\|u\|_{L^1[0,T]}\geq(d-T\|A\|)/\|B\|$, and hence
$\|u\|_{L^\infty[0,T]}\geq(d-T\|A\|)/(\|B\|\,T)\to+\infty$
as $T\to0$. The $L^1$ norm is bounded below (by the target
distance) but not necessarily above; the $L^\infty$ norm
\emph{necessarily} diverges.

\end{remark}

\subsection{Scope and strategy}\label{SUB_SEC_Bilan_Regularite}

\begin{table}
\centering
\renewcommand{\arraystretch}{1.3}
\begin{tabular}{p{0.11\textwidth} p{0.001\textwidth} p{0.14\textwidth} p{0.17\textwidth} p{0.17\textwidth} p{0.17\textwidth}}
\toprule
\textbf{Quantity} &   & \textbf{is always} & \textbf{Ass. \ref{ASS_control}} &
\textbf{Ass. \ref{ASS_control} $+$ \ref{ASS_coupling_graph}} &
 \textbf{Ass. \ref{ASS_control} $+$ \ref{ASS_coupling_graph} $+$ \ref{ASS_drift_absorption}}\\
\midrule
$\tau_{\mathrm{AUS}}(X_f)$ & 
& min?
&   \textbf{Yes} (Cor. \ref{COR_AUS_attainment}) & Yes
& Yes\\[4pt]
$\tau_{\mathrm{US}}(X_f)$ & 
&  min?
 & No &  \textbf{No} 
  (Rem.~\ref{REM_Blow-up})
 & \textbf{Yes} (
  Prop.~\ref{PRO_egalite})\\[4pt]
$\tau^{\max}_{\mathrm{AUS}}$ & 
& max?
& ?
& ?
& \textbf{Yes} (Prop.~\ref{PRO_sup_is_max})\\[4pt]
$\tau^{\max}_{\mathrm{US}}$ & 
& max?
& ?
& ?
& \textbf{Yes} (Prop.~\ref{PRO_sup_is_max}) \\[4pt]
$\tau^{\max}_{\mathrm{US}}$ &  &  $\tau^{\max}_{\mathrm{AUS}}?$
&  \textbf{Yes} (
  Lem. \ref{LEM_approx_cont})
& Yes
& Yes \\[4pt]
$\tau_{\mathrm{US}}(X_f)$ & &  $\tau_{\mathrm{AUS}}(X_f)?$
&  No
& \textbf{No} (Cor. \ref{COR_counterexample})
& \textbf{Yes} (
  Prop.~\ref{PRO_egalite}) \\
\bottomrule
\end{tabular}
\caption{
The worst-case supremum is a maximum
(Proposition~\ref{PRO_sup_is_max}) under Assumptions \ref{ASS_control}, \ref{ASS_coupling_graph} and the drift absorption \ref{ASS_drift_absorption}; it is not known whether this holds under Assumptions \ref{ASS_control} and \ref{ASS_coupling_graph} alone, or under \ref{ASS_control} alone (question marks). Since the drift absorption holds generically (Proposition~\ref{PRO_Generic_Attainable_Neighborhood}), the result is generically valid. }\label{TAB_bilan_regularite}
\end{table}

Table \ref{TAB_bilan_regularite} summarizes the regularity and reachability properties. In the remainder of the paper, the constructive bounds (Section \ref{SEC_Upper_bound}) and their application (Section \ref{SEC_Chiral}) are stated for the approximate synthesis times $\tau_{\mathrm{ASS}}$ and $\tau_{\mathrm{AUS}}$; by Lemma~\ref{LEM_approx_cont} they transfer to the exact times under Assumption~\ref{ASS_control}.

The analysis relies on the factorization of generalized trajectories:
\begin{equation}
Z(t) = \exp(v(t)B + tA_D) W(t)
\end{equation}
where the dynamics of $W$ are governed by the effective interaction generator. The control strategy decouples the transport management of $W$ from the phase adjustment $S_v(t)= \exp(v(t)B + tA_D) $. Explicit protocols first provide an upper bound for the travel time on $W$. The global phase adjustment is then performed in a second step.
This sequential approach is geometrically suboptimal, but it captures the asymptotic order of the estimates for the rotator of Proposition~\ref{PRO_rotateur_introduction}.

\section{Constructive Upper Bounds via Optimal Routing}\label{SEC_Upper_bound}

To establish an upper bound on the infimum transfer time, we construct explicit control sequences by 
restricting the dynamics to a spanning tree of the coupling graph: the problem then reduces to a 
routing procedure on a weighted tree.

\subsection{Non-Resonant Transitions and Edge Isolation}\label{SUBSEC_Edge_Isolation}

The algebraic structure of the set $CC(t)$ enables the isolation of specific transitions within the coupling graph, subject to specific spectral conditions on the matrix $B$.

By defining the effective generator parametrized by a phase variable $w \in \mathbb{R}$:
\begin{equation}
M(w, t) \;=\; \exp(-w B - t A_D) A_{ND} \exp(w B + t A_D),
\end{equation}
the compact set of admissible convexified generators $\widetilde{\Omega}(t) \in CC(t) \subset \mathfrak{su}(d)$ is given by:
\begin{equation}
CC(t) \;=\; \overline{\operatorname{conv}} \{ M(w, t) \mid w \in \mathbb{R} \}.
\end{equation}
Taking the topological closure in the finite-dimensional Lie algebra $\mathfrak{su}(d)$ implies that any element $M \in CC(t)$ can be represented as the barycenter of $M(w, t)$ with respect to a probability measure $\mu$ on $\mathbb{R}$:
\begin{equation}\label{EQ_Convexification_measure}
    M = \int_\mathbb{R} M(w, t) \, \mathrm{d}\mu(w).
\end{equation}
This representation lets us read off the non-zero entries of $M$ via the Fourier moments of the measure $\mu$.

\begin{proposition}[Partial Edge Isolation in Resonant Graphs]\label{PRO_partial_edge_isolation}
Assume the transition $(j,k)$ (associated with entry $A_{jk}$ of $A$) is at least $p$-harmonic free, meaning there are no 
integer harmonics in the spectral gap ratio up to order $p-1$. That is, for all coupled pairs $(l,m)$ with $A_{lm}\neq 0$, 
the ratio $\omega_{lm}/\omega_{jk} \notin \{\pm 1, \pm 2, \dots, \pm (p-1) \}$ unless $\{l,m\}=\{j,k\}$.

Then, for every $\theta$ in $\R$, there exists an element $M^{(j,k)}_{part} \in CC(t)$ whose entries are identically zero except for the entries $(j,k)$ and $(k,j)$ for which
\begin{equation}
    (M^{(j,k)}_{part}(\theta))_{jk} = 
    e^{i\theta} |A_{jk}| \cos\left( \frac{\pi}{p+1}\right) 
     =-\overline{(M^{(j,k)}_{part}(\theta))_{kj} }.
\end{equation}
\end{proposition}

\begin{proof}
The core of the proof (use of the Fejér-Riesz theorem, and reduction to the computation of an operator norm, Step 2 and Step 3 below) is borrowed from the book \cite{Dumitrescu2017}, Chapters 1-3.

Without loss of generality, normalize the target spectral gap to $\omega_{jk} = 1$. 
The dynamics occur over the full real line $v \in \mathbb{R}$. The coupling graph may contain a mixture of 
integer harmonic resonances ($\omega_{lm} \in \mathbb{Z}$) and non-integer spectral 
gaps ($\omega_{lm} \notin \mathbb{Z}$). 

To isolate the target edge, we perform a two-stage convexification.

\textbf{Step 1: Ergodic elimination of non-integer gaps.} 
We introduce a continuous phase parameter $w \in [0, 2\pi)$ and define the arithmetic 
progression $v_n(w) = w + 2\pi n$ for $n \in \mathbb{Z}$. By construction, the matrix 
element of $M(v, t)$ corresponding to a transition $(l,m)$ is proportional to 
$e^{i \omega_{lm} v}$. Substituting $v_n(w)$ into the Cesàro mean 
$S_N(v, t) = \frac{1}{2N+1} \sum_{n=-N}^N M(v_n(w), t)$ yields terms proportional to:
$$e^{i \omega_{lm} w} \frac{1}{2N+1} \sum_{n=-N}^N e^{i 2\pi n \omega_{lm}}$$
If $\omega_{lm} \in \mathbb{Z}$, the geometric sum is $2N+1$, and the term 
simplifies to $e^{i \omega_{lm} w}$.
 If $\omega_{lm} \notin \mathbb{Z}$, the sum is  bounded by $\frac{1}{|\sin(\pi \omega_{lm})|}$:
\[
\left | e^{i \omega_{lm} v} \frac{1}{2N+1} \sum_{n=-N}^N e^{i 2\pi n \omega_{lm}} \right | = \frac{1}{2N+1} \left | e^{-i 2\pi N \omega_{lm}} \frac{1 - e^{i 2\pi (2N+1) \omega_{lm}}}{1 - e^{i 2\pi \omega_{lm}}} \right | \le \frac{1}{(2N+1) |\sin(\pi \omega_{lm})|}
\] 
 and this term  vanishes as $N \to \infty$. Because $CC(t)$ is defined as the operator-norm closure of the convex hull, the limit operator belongs to the accessible set:$$M_\infty(v, t) = \lim_{N \to \infty} S_N(v, t) \in CC(t)$$and its non-zero entries correspond exclusively to integer spectral gaps.

\textbf{Step 2: Fejér-Riesz elimination of integer harmonics.}
The problem is now reduced to the torus $\mathbb{T} = \mathbb{R} / 2\pi\mathbb{Z}$. Taking advantage of the measure formulation \eqref{EQ_Convexification_measure}, we seek a probability measure $\mu(v)$ on $\mathbb{T}$ whose Fourier coefficients satisfy $c_k = \int_{\mathbb{T}} e^{ikv} \mathrm{d}\mu(v) = 0$ for all $|k| \ge p$, while maximizing the modulus of the fundamental harmonic $|c_1|$.

Since $c_k = 0$ for $|k| \ge p$, the measure $\mu$ must admit a density which is a non-negative trigonometric polynomial of degree at most $N = p-1$:
\begin{equation}
    \mathrm{d}\mu(v) = \left( 1 + \sum_{k=1}^{p-1} c_k e^{ikv} + \bar{c}_k e^{-ikv} \right) \frac{\mathrm{d}v}{2\pi} \ge 0.
\end{equation}
Recall now the classical Fejér-Riesz theorem. 
\begin{theorem}[Fejér-Riesz theorem, \cite{fejer1915uber} and \cite{riesz1915uber}]\label{THE_Fejer-Riesz}
Let $P(\theta)$ be a trigonometric polynomial of degree $n$, given by
\begin{equation}
    P(\theta) = \sum_{k=-n}^{n} c_k e^{ik\theta},
\end{equation}
with complex coefficients $c_k$ satisfying $c_{-k} = \overline{c_k}$. 
If $P(\theta) \ge 0$ for all real $\theta$, then there exists an algebraic polynomial $A(z)$ of degree at most $n$, 
\begin{equation}
    A(z) = \sum_{k=0}^{n} a_k z^k,
\end{equation}
with complex coefficients $a_k$, such that for all $\theta \in \mathbb{R}$,
\begin{equation}
    P(\theta) = \left| A(e^{i\theta}) \right|^2.
\end{equation}
Furthermore, the roots of $A(z)$ can always be chosen such that they all lie within or on the unit circle ($|z| \le 1$).
\end{theorem}

By the Fejér-Riesz theorem,  
 since $\mu$ is a non-negative measure on $\mathbb{T}$ with a Fourier support  bounded by $p-1$, its density can be factorized as the squared modulus of an algebraic polynomial of degree at most $p-1$:
\begin{equation}
    \mathrm{d}\mu(v) = \left| \sum_{n=0}^{p-1} a_n e^{inv} \right|^2 \frac{\mathrm{d}v}{2\pi}, \quad a_n \in \mathbb{C}.
\end{equation}

The probability conservation constraint $c_0 = 1$ enforces a unit Euclidean norm on the coefficient vector $a \in \mathbb{C}^p$: $\sum_{n=0}^{p-1} |a_n|^2 = 1$. The fundamental coefficient expands into the bilinear form $c_1 = \sum_{n=0}^{p-2} a_n \overline{a_{n+1}}$. 

\textbf{Step 3:  Algebraic Derivation of the Spectral Norm}

By introducing the shift matrix $S \in \mathcal{M}_{p}(\mathbb{R})$ defined by $S_{m, m+1} = 1$ and $0$ elsewhere, the coefficient $c_1$ is expressed as a standard sesquilinear form:
$$c_1 = \sum_{n=0}^{p-2} a_n \overline{a_{n+1}}= a^\dagger S^\dagger a$$
where $a^\dagger$ is the conjugate transpose of vector $a$.

The objective is to maximize the modulus $|c_1|$ subject to the unit norm constraint $\|a\|_2^2 = a^\dagger a = 1$.

Applying a diagonal unitary phase transformation $D(\alpha) = \text{diag}(1, e^{i\alpha}, e^{i2\alpha},
\dots, e^{i(p-1)\alpha})$ to a vector $a$ generates a new vector $\tilde{a} = D(\alpha)a$. 
The norm is 
invariant ($\|\tilde{a}\|_2 = \|a\|_2$), and the sesquilinear form is modified by a global 
phase:
$$\tilde{a}^\dagger S^\dagger \tilde{a} = e^{-i\alpha} a^\dagger S^\dagger a = e^{-i\alpha} 
c_1$$
There always exists an angle $\alpha$ such that the resulting coefficient is  real and 
positive. Therefore, optimizing the modulus over $\mathbb{C}^p$ is geometrically equivalent 
to maximizing its real part:
$$\max_{\|a\|_2=1} |a^\dagger S a| = \max_{\|a\|_2=1} \text{Re}(a^\dagger S^\dagger a)$$
The real part of a complex quadratic form is extracted via the Hermitian component of its 
generating matrix:
$$\text{Re}(a^\dagger S^\dagger a) = a^\dagger \left( \frac{S + S^\dagger}{2} \right) a$$

Define the matrix $T = \frac{1}{2}(S + S^\dagger)$. By construction, $T$ is a real symmetric 
matrix, with zeros on the main diagonal and values $1/2$ on the first sub-diagonal and super-
diagonal ($T_{m, m+1} = T_{m+1, m} = 1/2$).

The constrained problem reduces to the Rayleigh quotient:

$$\max_{\|a\|_2=1} a^\dagger T a$$

According to the Rayleigh-Ritz min-max theorem, the supremum of this bilinear form on the unit sphere 
corresponds to the largest eigenvalue of the Hermitian matrix $T$, which equals here its spectral norm 
$\|T\|_2$ due to the symmetry of the spectrum around zero.

 The maximum modulus is hence bounded by the spectral radius (see Lemma \ref{LEM_eigenvectors_shift_matrix} in Appendix \ref{SEC_Appendix}):
\begin{equation}
    \max |c_1| = \|T\|_2 = \cos\left( \frac{\pi}{p+1} \right),
\end{equation}
with the normalized optimal eigenvector yielding $a_n = \sqrt{\frac{2}{p+1}} \sin\left(\frac{(n+1)\pi}{p+1}\right)$.

Let $\mu_0$ be the optimal measure yielding the strictly positive maximal first Fourier 
coefficient $c_1^{(0)} = \cos\left( \frac{\pi}{p+1} \right)$. By construction, the $(j,k)$ 
entry of the pre-conditioned generator accumulates phase linearly with $v$, evaluating to $
(M_\infty(v, t))_{jk} = |A_{jk}| e^{i(v + \xi_{jk}(t))}$, where $\xi_{jk}(t)$ is the 
dynamical phase at time $t$. 

We define the shifted measure $\mu_{opt}(v) = \mu_0(v - \phi)$. A simple change of variables $w = v - \phi$ formalizes the translation of the first moment:
$$
\int_{\mathbb{T}} e^{iv} \mathrm{d}\mu_{opt}(v) = \int_{\mathbb{T}} e^{i(w+\phi)} \mathrm{d}\mu_0(w) = e^{i\phi} c_1^{(0)}.
$$
Integrating the generator $M_\infty(v, t)$ against $\mu_{opt}$ evaluates the target entry to:
$$
(M^{(j,k)}_{part})_{jk} = |A_{jk}| e^{i\xi_{jk}(t)} \int_{\mathbb{T}} e^{iv} \mathrm{d}\mu_{opt}(v) = |A_{jk}| c_1^{(0)} e^{i(\xi_{jk}(t) + \phi)}.
$$
Setting the continuous shift parameter to $\phi = \theta - \xi_{jk}(t)$ aligns the argument of the entry with the target phase $\theta$, yielding $e^{i\theta} |A_{jk}| \cos\left( \frac{\pi}{p+1} \right)$. As established, all non-target edges evaluate to zero, concluding the proof.
\end{proof}

Recall that we have factorized the generalized solutions (Definition \ref{DEF_Generalized_solution}) $Z(t)=R(t) Y^{\Omega}_{t,0}$ 
following \eqref{eq:decomp}:
\[Z(t) \;=\; \underbrace{\exp\!\bigl(v(t)\,B + t\,A_D\bigr)}_{S_v(t)}
  \;\cdot\; \underbrace{W^{\widetilde{\Omega}}(t)}_{W},
\]
where $v(t) = \int_0^t u(s)\,ds$ for some $u\in L^1_{\mathrm{loc}}$
(unbounded) and $W^{\widetilde{\Omega}}$ satisfies
\[
  \dot{W}^{\widetilde{\Omega}}(t) \;=\; 
  \widetilde{\Omega}(t)\,W^{\widetilde{\Omega}}(t),
  \qquad \widetilde{\Omega}(t)\in CC(t)\subset\mathfrak{su}(d),
\]
with $CC(t) = \overline{\operatorname{conv}\{e^{-wB-tA_D}A_{ND}e^{wB+ tA_D}:w\in\mathbb{R}\}}$
compact and $\|\widetilde{\Omega}(t)\|_{op}\leq L_{CC}$.

\begin{corollary}\label{COR_STLC_W}
Under Assumptions \ref{ASS_control} and \ref{ASS_coupling_graph}, and for $d \ge 3$, for every $\varepsilon > 0$, the set
\begin{equation}
  \mathcal{A}^W_{\leq \varepsilon}
  := \bigcup_{\substack{\widetilde{\Omega}\in L^{\infty}([0,\varepsilon],\,\mathfrak{su}(d))\\
      \widetilde{\Omega}(t)\in CC(t)}}
     \{ W^{\widetilde{\Omega}}(t) : 0 \le t \le \varepsilon \}
\end{equation}
contains a neighborhood of $I_d$ in $SU(d)$.
\end{corollary}

\begin{proof}
We show that the identity $I_d$ is an interior point of the reachable set $\mathcal{A}^W_{\leq \varepsilon}$ for every $\varepsilon > 0$. The proof proceeds in two steps: first, we establish that the convexified control set $CC(t)$ contains a symmetric set whose Lie algebra is $\mathfrak{su}(d)$ (for $d \ge 3$), and second, we use a classical argument based on Chow theorem for symmetric systems to prove  small-time local controllability (STLC).

\medskip
\noindent\textbf{Step 1: Full-rank generation of the Lie algebra.}
By Proposition \ref{PRO_partial_edge_isolation}, under Assumption \ref{ASS_coupling_graph}, the set $CC(t)$ contains, for every $t$ and for each edge $e = (j,k)$ of a spanning tree $\mathbf{T}$ of the coupling graph, the family of matrices
\begin{equation}
  M^{(j,k)}_{part}(\theta) = \gamma_{jk} \left( e^{i\theta} E_{jk} - e^{-i\theta} E_{kj} \right) \in CC(t), \quad \theta \in [0, 2\pi),
\end{equation}
where $\gamma_{jk} = |A_{jk}| \cos\left(\frac{\pi}{p_{jk}+1}\right) > 0$, with $p_{jk}:=M_{jk}$ the harmonic-freeness order of the edge $(j,k)$ defined in Proposition~\ref{PRO_Result_Constructive_State_to_State}. Let $\mathcal{S}$ be the union of these families over all edges of $\mathbf{T}$. The set $\mathcal{S}$ is symmetric with respect to the origin, since $M^{(j,k)}_{part}(\theta + \pi) = -M^{(j,k)}_{part}(\theta)$.

We claim that the Lie algebra generated by $\mathcal{S}$ is $\mathfrak{su}(d)$, provided $d \ge 3$. To see this, write $G_{jk} = \gamma_{jk} E_{jk}$. The generators of the edge $(j,k)$ are $S_{jk} = G_{jk} - \overline{G}_{jk}$ and $C_{jk} = i(G_{jk} + \overline{G}_{jk})$. For any two adjacent edges $(j,k)$ and $(k,\ell)$ in the tree (with $j \neq \ell$), the Lie brackets are:
\begin{align}
  [S_{jk}, S_{k\ell}] &= \gamma_{jk} \gamma_{k\ell} (E_{j\ell} - E_{\ell j}) =: S^{j\ell}_{eff}, \\
  [S_{jk}, C_{k\ell}] &= i \gamma_{jk} \gamma_{k\ell} (E_{j\ell} + E_{\ell j}) =: i C^{j\ell}_{eff}.
\end{align}
These brackets generate the effective Givens generators for the pair $(j,\ell)$. By iterating this process along the unique paths in the tree $\mathbf{T}$, we obtain $S^{ab}_{eff}$ and $C^{ab}_{eff}$ for all pairs $(a,b)$ of vertices in the tree. The matrix units $E_{ab}$ and $E_{ba}$ are then recovered as linear combinations of $S^{ab}_{eff}$ and $C^{ab}_{eff}$. Finally, the Cartan subalgebra is generated by the commutators $[E_{ab}, E_{ba}] = E_{aa} - E_{bb}$, which span the $(d-1)$-dimensional space of traceless diagonal matrices. Thus, $\mathrm{Lie}(\mathcal{S}) = \mathfrak{su}(d)$.

For $d=2$, the spanning tree has a single edge, so the bracket iteration above produces no non-adjacent effective generators; that iteration is the part of the argument that requires $d\ge3$. (The generators do not commute: $S_{12}=\gamma(E_{12}-E_{21})$, $C_{12}=i\gamma(E_{12}+E_{21})$ satisfy $[S_{12},C_{12}]=2i\gamma^2(E_{11}-E_{22})$, and $\{S_{12},C_{12},E_{11}-E_{22}\}$ is a basis of $\mathfrak{su}(2)$; the STLC conclusion thus also holds for $d=2$, by a direct argument.)

\medskip
\noindent\textbf{Step 2: STLC via Chow theorem.}
We define $DD=\overline{\operatorname{conv}(\mathcal{S})}$.
Let $V = \mathrm{span}(DD)$. Since $DD$ is a compact, symmetric, and convex subset of $\mathfrak{su}(d)$, it contains the origin in its interior relative to $V$. Consequently, there exists a Euclidean ball $\mathcal{B} \subset DD$ centered at the origin in $V$.

Consider the restricted control system $\dot{x} = \Omega x$ with $\Omega \in \mathcal{B}$. Because $\mathcal{B} \subseteq DD$, the reachable set for the restricted system is contained in the reachable set of the original system $\mathcal{A}^W_{\leq T}$.

The set $\mathcal{B}$ defines a standard sub-Riemannian structure on $SU(d)$ via right translations. For this restricted driftless system, the minimal time function to reach a point from $I_d$ coincides with the Carnot-Carathéodory distance $d_{SR}$ associated with $\mathcal{B}$. Therefore, the reachable set in time $\leq \varepsilon$ using controls in $\mathcal{B}$ is the closed sub-Riemannian ball of radius $\varepsilon$ centered at $I_d$, denoted by $B_{SR}(I_d, \varepsilon)$.

By hypothesis, $\mathrm{Lie}(\mathcal{S}) = \mathfrak{su}(d)$, which trivially implies $\mathrm{Lie}(\mathcal{B}) = \mathfrak{su}(d)$ since $V = \mathrm{span}(\mathcal{B})$. This ensures that the associated distribution satisfies the Lie Algebra Rank Condition (LARC) everywhere on $SU(d)$. By the Chow-Rashevskii theorem, the topology induced by the Carnot-Carathéodory distance $d_{SR}$ is equivalent to the standard manifold topology on $SU(d)$.

Consequently, for any $\varepsilon > 0$, the sub-Riemannian ball $B_{SR}(I_d, \varepsilon)$ contains a neighborhood of $I_d$ in the manifold topology. The result follows since  $B_{SR}(I_d, \varepsilon) \subseteq \mathcal{A}^W_{\leq \varepsilon}$. 
\end{proof}

\subsection{Elementary Phased Givens Rotations}

Before addressing global routing on the graph, we derive time bounds for the elementary operations 
of the convexified set $CC(t)$. The basic operation is the phased Givens rotation on the 
two-dimensional subspace $\text{Span}(\phi_j^B, \phi_k^B)$. The general theory is standard
(e.g., \cite{Nielsen_Chuang_2010}, Chapter 4; \cite{Williams2011}, Chapter 2); however, the explicit 
parameters of Lemma \ref{LEM_Effective_Givens} are not readily extracted from these references and
are proved below, as the relative phases of Section~\ref{SEC_Upper_bound_StS} must be fixed from them.

  Let $(E_{mn})_{1 \leq m,n \leq d}$ denote the canonical matrix units with respect to $(\phi_1^B, \dots ,\phi_d^B)$, the eigenbasis of matrix $B$ , such that $E_{mn}$ maps $\phi_n^B$ to $\phi_m^B$ and acts trivially on the rest of the basis. 
\begin{definition}[Phased Givens Rotation]
A phased Givens rotation of angle $\theta \ge 0$ and relative phase $\phi \in [0, 2\pi)$ operating on the two-dimensional subspace $\text{Span}(\phi_j^B, \phi_k^B)$ is the unitary transformation defined by:
\begin{equation}
    R_{jk}(\theta, \phi) = \exp\left( \theta \left( e^{i\phi} E_{jk} - e^{-i\phi} E_{kj} \right) \right).
\end{equation}
\end{definition}

\begin{lemma}[Effective parameters for Givens Rotations]\label{LEM_Effective_Givens}
Let $\psi=(\psi_1,\psi_2)^T$ and $\tilde{\psi}=(\tilde{\psi}_1,\tilde{\psi}_2)^T$ be vectors in the unit sphere of $\mathbb{C}^2$. Then, choosing $\theta_1= \arg(\psi_1) - \arg(\tilde{\psi}_1)$, $\theta_2=\arg(\psi_2) - \arg(\tilde{\psi}_2)$,
\[\rho = \arccos(|\psi_1||\tilde{\psi}_1| + |\psi_2||\tilde{\psi}_2|) = 
\left| \arctan\left(\frac{|\psi_2|}{|\psi_1|}\right) - \arctan\left(\frac{|\tilde{\psi}_2|}{|\tilde{\psi}_1|}\right) \right|
\]
and $\phi=\arg\left( (|\tilde{\psi}_1||\psi_2| - |\tilde{\psi}_2||\psi_1|) e^{i(\arg(\psi_1) - \arg(\psi_2))} \right)$, one has
$$\exp \begin{pmatrix} 0 & \rho e^{i \phi} \\ -\rho e^{-i \phi} & 0 \end{pmatrix} \begin{pmatrix} \psi_1 \\ \psi_2 \end{pmatrix} = \begin{pmatrix} e^{i \theta_1} \tilde{\psi}_1 \\ e^{i \theta_2} \tilde{\psi}_2 \end{pmatrix}.$$
\end{lemma}

\begin{remark} \label{REM_convention_composantes_nulles}
In Lemma \ref{LEM_Effective_Givens}, we used the convention $\arctan(+\infty)=\frac{\pi}{2}$. 
If $|\tilde{\psi}_1||\psi_2| - |\tilde{\psi}_2||\psi_1| \ge 0$, then $\phi = \arg(\psi_1) - \arg(\psi_2)$. Otherwise, $\phi = \arg(\psi_1) - \arg(\psi_2) + \pi$. The expressions of $\theta_1$, $\theta_2$, $\phi$ remain valid for any arbitrary choice of arguments for zero components of $\psi$ or $\tilde{\psi}$, provided the arbitrary value assigned to the argument of any specific zero component remains identical throughout the calculation of all parameters ($\theta_1$, $\theta_2$, $\phi$) in which it appears.
\end{remark}

\begin{proof}[Proof of Lemma \ref{LEM_Effective_Givens}]
Let $R(\rho, \phi) = \exp \begin{pmatrix} 0 & \rho e^{i \phi} \\ -\rho e^{-i \phi} & 0 \end{pmatrix}$.
By diagonalizing the anti-Hermitian generator or using its Taylor series expansion, the matrix exponential resolves to:

$$R(\rho, \phi) = \begin{pmatrix} \cos \rho & \sin \rho ~ e^{i \phi} \\ -\sin \rho ~ e^{-i \phi} & \cos \rho \end{pmatrix}$$

The constraint is $R(\rho, \phi) \psi = \Psi$, where the target vector is defined as $\Psi = (e^{i\theta_1} \tilde{\psi}_1, e^{i\theta_2} \tilde{\psi}_2)^T$.
Define the polar representations $\psi_k = r_k e^{i\gamma_k}$ and 
$\tilde{\psi}_k = \tilde{r}_k e^{i\delta_k}$ for $k \in \{1, 2\}$.
Setting $\theta_k = \gamma_k - \delta_k$ aligns the phases of the 
target vector with the initial vector, such that 
$\Psi_k = \tilde{r}_k e^{i\gamma_k}$. The operational system of 
equations becomes:

$$\begin{cases} \cos \rho \cdot r_1 e^{i\gamma_1} + \sin \rho e^{i\phi} \cdot r_2 e^{i\gamma_2} = \tilde{r}_1 e^{i\gamma_1} \\ -\sin \rho e^{-i\phi} \cdot r_1 e^{i\gamma_1} + \cos \rho \cdot r_2 e^{i\gamma_2} = \tilde{r}_2 e^{i\gamma_2} \end{cases}$$

Assuming $r_1, r_2 \neq 0$ (the limits hold via continuity for the boundary cases), isolate $\sin \rho e^{i\phi}$ from the first equation and $\sin \rho e^{-i\phi}$ from the second:

$$\sin \rho \cdot e^{i\phi} = \frac{\tilde{r}_1 - \cos \rho \cdot r_1}{r_2} e^{i(\gamma_1 - \gamma_2)}$$

$$\sin \rho \cdot e^{-i\phi} = \frac{\cos \rho \cdot r_2 - \tilde{r}_2}{r_1} e^{-i(\gamma_1 - \gamma_2)}$$

Taking the complex conjugate of the second relation yields a second expression for $\sin \rho e^{i\phi}$. Equating them gives:

$$\frac{\tilde{r}_1 - \cos \rho \cdot r_1}{r_2} = \frac{\cos \rho \cdot r_2 - \tilde{r}_2}{r_1}$$

Rearranging this equality and using $r_1^2 + r_2^2 = 1$ leads to:

$$\cos \rho = r_1 \tilde{r}_1 + r_2 \tilde{r}_2 = |\psi_1||\tilde{\psi}_1| + |\psi_2||\tilde{\psi}_2|$$

Since $\|\psi\| = \|\tilde{\psi}\| = 1$, the Cauchy-Schwarz inequality ensures $0 \le r_1 \tilde{r}_1 + r_2 \tilde{r}_2 \le 1$. Thus, 
$\rho = \arccos(|\psi_1||\tilde{\psi}_1| + |\psi_2||\tilde{\psi}_2|)$ 
is well defined and restricted to $[0, \pi/2]$.

To isolate the off-diagonal parameter, substitute $\cos \rho$ back 
into the first isolated expression:

$$\sin \rho e^{i\phi} = (\tilde{r}_1 r_2 - \tilde{r}_2 r_1) e^{i(\gamma_1 - \gamma_2)}$$

Given that $\sin \rho \ge 0$ for $\rho \in [0, \pi/2]$, the complex 
argument $\phi$ is determined by the right-hand side of this equation.

To obtain the alternative expression of $\rho$ (in terms of arctan), 
note that since the vectors $\psi$ and $\tilde{\psi}$ belong to the 
unit sphere of $\mathbb{C}^2$, the moduli of their components satisfy 
$|\psi_1|^2 + |\psi_2|^2 = 1$ and $|\tilde{\psi}_1|^2 + 
|\tilde{\psi}_2|^2 = 1$.

These moduli can be parameterized using angles $\alpha$ and $\beta$ 
in the interval $[0, \pi/2]$, using the the two‑argument arctangent function $\operatorname{atan2}$:
$|\psi_1| = \cos \alpha$ and $|\psi_2| = \sin \alpha$, which implies $\alpha = \operatorname{atan2}(|\psi_2|, |\psi_1|)$
$|\tilde{\psi}_1| = \cos \beta$ and $|\tilde{\psi}_2| = \sin \beta$, 
which implies $\beta = \operatorname{atan2}(|\tilde{\psi}_2|, |\tilde{\psi}_1|)$. 

Recalling the established expression for $\cos \rho$:

$$\cos \rho = |\psi_1||\tilde{\psi}_1| + |\psi_2||\tilde{\psi}_2|$$

Substituting the trigonometric parametrization yields:

$$\cos \rho = \cos \alpha \cos \beta + \sin \alpha \sin \beta = \cos(\alpha - \beta)$$

By restricting $\rho$ to the interval $[0, \pi/2]$ to parameterize an effective Givens rotation, and given that $\alpha, \beta \in [0, \pi/2]$, the angular distance is given by the absolute value of their difference:

$$\rho = |\alpha - \beta|$$

Replacing $\alpha$ and $\beta$ with their expressions in terms of $\arctan$ yields the alternative expression of $\rho$.
\end{proof}

\begin{corollary}[Modulus routing and phase locking] \label{COR_SU2_phase_locking}
Let $\psi=(\psi_1,\psi_2)^T$ be a vector in the unit sphere of $\mathbb{C}^2$. Then, for every $\theta_1$ in $\mathbb{R}$, choosing $\phi_1 = \arg(\psi_1) - \arg(\psi_2) + \pi$, $\phi_2 = \theta_1 - \arg(\psi_2)$ (with the same conventions as in Remark \ref{REM_convention_composantes_nulles} for the argument of zero components), $\rho_1 = \arcsin(|\psi_1|)$ and $\rho_2 = \frac{\pi}{2}$ , one has

$$\exp \begin{pmatrix} 0 & \rho_2 e^{i \phi_2} \\ -\rho_2 e^{-i \phi_2} & 0 \end{pmatrix} \exp \begin{pmatrix} 0 & \rho_1 e^{i \phi_1} \\ -\rho_1 e^{-i \phi_1} & 0 \end{pmatrix} \begin{pmatrix} \psi_1 \\ \psi_2 \end{pmatrix} =  \begin{pmatrix}  e^{i \theta_1} \\ 0\end{pmatrix}.$$

Notice that the sum $\rho_1 + \rho_2$ belongs to $[0,\pi]$ and solely depends upon $|\psi_1|$.
\end{corollary}

\begin{proof}[Proof of Corollary \ref{COR_SU2_phase_locking}]
To map $\psi$ to $(e^{i\theta_1}, 0)^T$ for an arbitrary $\theta_1$ using two Givens rotations $R(\rho_2, \phi_2) R(\rho_1, \phi_1)$, the intermediate state must hold its entire norm in the second component. Let $R(\rho, \phi) = \begin{pmatrix} \cos\rho & \sin\rho e^{i\phi} \\ -\sin\rho e^{-i\phi} & \cos\rho \end{pmatrix}$.

Applying the first rotation $R(\rho_1, \phi_1)$ with $\rho_1 = \arcsin(|\psi_1|)$ yields $\cos\rho_1 = \sqrt{1-|\psi_1|^2} = |\psi_2|$ and $\sin\rho_1 = |\psi_1|$. With $\phi_1 = \arg(\psi_1) - \arg(\psi_2) + \pi$, the phase term becomes $e^{i\phi_1} = -e^{i(\arg(\psi_1) - \arg(\psi_2))}$. The application to $\psi$ yields:

$$R(\rho_1, \phi_1) \psi = \begin{pmatrix} |\psi_2| & -|\psi_1| e^{i(\arg(\psi_1) - \arg(\psi_2))} \\ |\psi_1| e^{-i(\arg(\psi_1) - \arg(\psi_2))} & |\psi_2| \end{pmatrix} \begin{pmatrix} |\psi_1| e^{i\arg(\psi_1)} \\ |\psi_2| e^{i\arg(\psi_2)} \end{pmatrix}$$

Evaluating the components:

1. The first component is:

$$|\psi_2||\psi_1| e^{i\arg(\psi_1)} - |\psi_1||\psi_2| e^{i(\arg(\psi_1) - \arg(\psi_2) + \arg(\psi_2))} = 0$$

2. The second component is:

$$|\psi_1|^2 e^{-i(\arg(\psi_1) - \arg(\psi_2))} e^{i\arg(\psi_1)} + |\psi_2|^2 e^{i\arg(\psi_2)} = (|\psi_1|^2 + |\psi_2|^2) e^{i\arg(\psi_2)} = e^{i\arg(\psi_2)}$$

The intermediate state is  $(0, e^{i\arg(\psi_2)})^T$.

Now, applying the second rotation $R(\rho_2, \phi_2)$ with $\rho_2 = \frac{\pi}{2}$ yields $\cos\rho_2 = 0$ and $\sin\rho_2 = 1$. The transformation becomes:

$$R(\rho_2, \phi_2) \begin{pmatrix} 0 \\ e^{i\arg(\psi_2)} \end{pmatrix} = \begin{pmatrix} 0 & e^{i\phi_2} \\ -e^{-i\phi_2} & 0 \end{pmatrix} \begin{pmatrix} 0 \\ e^{i\arg(\psi_2)} \end{pmatrix} = \begin{pmatrix} e^{i(\phi_2 + \arg(\psi_2))} \\ 0 \end{pmatrix}$$

By choosing $\phi_2 = \theta_1 - \arg(\psi_2)$, the phase of the first component evaluates  to \[\exp \left [i(\theta_1 - \arg(\psi_2) + \arg(\psi_2)) \right ] = e^{i\theta_1}.\] This establishes the mapping $R(\rho_2, \phi_2) R(\rho_1, \phi_1) \psi = (e^{i\theta_1}, 0)^T$.

Regarding the constraint on the sum of the parameters:

$$\rho_1 + \rho_2 = \arcsin(|\psi_1|) + \frac{\pi}{2}$$

Since $\psi$ is restricted to the unit sphere of $\mathbb{C}^2$, its modulus $|\psi_1| \in [0, 1]$, which implies $\arcsin(|\psi_1|) \in [0, \frac{\pi}{2}]$. Therefore, the sum $\rho_1 + \rho_2$ maps directly to the interval $[\frac{\pi}{2}, \pi] \subset [0, \pi]$ and is solely determined by $|\psi_1|$.
\end{proof}

\begin{proposition}[Time Bound for Elementary Rotations]\label{PRO_Time_bound_Gevens}
Let $(j,k)$ be an edge in the coupling graph, associated with the entry $A_{jk}$ of the non-diagonal drift $A_{ND}$, and let $p \ge 2$ be the order of the lowest integer harmonic resonance present in the graph (with $p=+\infty$ if the edge has no harmonic resonance). 

For every $\theta$ in $[0,\frac{\pi}{2}]$, for every $\phi$ in $[0,2\pi]$, there exist $\tau_{jk}(\theta)$
\begin{equation}
    \tau_{jk}(\theta)= \frac{\theta}{|A_{jk}| \cos\left( \frac{\pi}{p+1} \right)}
\end{equation}
 and  a sequence of continuous orbits $W_n:[0,\tau_{jk}(\theta)]\to SU(d)$ of the effective interaction dynamics \eqref{eq:W_dynamics}  such that 
$W_n(\tau_{jk}(\theta))$ converges to $R_{jk}(\theta, \phi)$ .
\end{proposition}

\begin{proof}
By relaxation theory (Theorem 8.7 in \cite{Agrachev2004}), the accessible set of the original system
is dense in that of the convexified dynamics: a unitary propagator generated by a constant element
$M \in CC(t)$ over a time $t$ is approximated uniformly by admissible control sequences in the same
time, so that $t$ is an upper bound for the infimum synthesis time $\tau_{jk}$.

From the Partial Edge Isolation (Proposition \ref{PRO_partial_edge_isolation}), we select the target complex phase of the isolated entry to be $\phi$. There exists an available constant generator $M^{(j,k)}_{part} \in CC(t)$ that annihilates all non-target edges, giving:
\begin{equation}
    M^{(j,k)}_{part} = |A_{jk}| \cos\left( \frac{\pi}{p+1} \right) \left( e^{i\phi} E_{jk} - e^{-i\phi} E_{kj} \right).
\end{equation}
Integrating the interaction dynamics $\dot{z} = M^{(j,k)}_{part} z$ over a continuous duration $t$ generates the unitary propagator $\exp(t M^{(j,k)}_{part})$. Identifying this propagator with the target operation $R_{jk}(\theta, \phi)$ imposes the algebraic constraint:
\begin{equation}
    t |A_{jk}| \cos\left( \frac{\pi}{p+1} \right) = \theta.
\end{equation}
Solving for $t$ yields the required interaction time $\tau_{jk}$.
\end{proof}

We define the characteristic time $C_e$ of an edge $e = (j,k)$. 
Since phase adjustments are handled later (via Kronecker torus density or composite return-trip 
sequences, see Sections~\ref{SEC_Upper_bound_StS} and~\ref{SEC_Contructive_bounds_Propagators}), 
the basic operation is to transfer population across the edge, which requires a rotation angle 
of at most $\theta = \pi/2$. Using the partial isolation generator, the worst-case time cost 
of this traversal is:
\begin{equation}\label{EQ_valeur_C_e}
    C_e = \frac{\pi}{2 |A_e| \cos\left( \frac{\pi}{p_e + 1} \right)} = \frac{\pi \alpha_e}{2 |A_e|},
\end{equation}
where $\alpha_e$ is the averaging efficiency correction factor introduced 
in Proposition \ref{PRO_Result_Constructive_State_to_State}.

\subsection{Constructive Bounds for State-to-State Transfer}\label{SEC_Upper_bound_StS}

The next result (Proposition \ref{PRO_au_plus_une_fois}) states that routing mass along a 
spanning tree with elementary Givens rotations can be achieved using each edge at most 
once. 

\begin{proposition}\label{PRO_au_plus_une_fois}
Let $T = (V, E)$ be a spanning tree of a graph where $d = |V| $. Let $p, q \in \mathbb{R}_{\geq 0}^d$ be 
the initial and target population vectors respectively, satisfying 
$\sum_{i \in V} p_i = \sum_{i \in V} q_i = 1$. Let the admissible operations be localized mass transfers 
along the edges of $T$, where an operation on edge $e = \{u, v\}$ parameterized by $\delta \in \mathbb{R}$ 
maps $p$ to $p'$ such that $p'_u = p_u - \delta$, $p'_v = p_v + \delta$, and $p'_w = p_w$ for all 
$w \in V \setminus \{u, v\}$, subject to the non-negativity constraint $p' \geq 0$. There exists a 
deterministic sequence of at most $|E|$ operations (at most one per edge) that transforms $p$ into $q$.
\end{proposition}

\begin{proof}
The proof proceeds by strong induction on the number of edges $|E|$ of the tree $T$. For $|E| = 0$, $T$ 
consists of a single vertex, $p = q$, and the sequence of operations is empty. Assume the proposition holds 
for any tree with fewer than $|E|$ edges.

For any edge $e = \{u, v\} \in E$, its removal partitions $T$ into two disjoint subtrees, $T_u$ and $T_v$, 
containing $u$ and $v$ respectively. By global mass conservation, the net mass transfer required from $T_u$ 
to $T_v$ across $e$ is determined by:
\begin{equation}
f_{u \to v} = \sum_{i \in V(T_v)} (q_i - p_i)
\end{equation}
Note that $f_{v \to u} = -f_{u \to v}$. We construct a directed flow graph $D_F = (V, E_F)$ by orienting 
each edge $\{u, v\} \in E$ as a directed edge $(u, v) \in E_F$ if and only if $f_{u \to v} > 0$. The weight 
of $(u, v)$ is assigned as $W(u, v) = f_{u \to v}$. Edges with zero required net flow 
are omitted from $E_F$.

Since the underlying graph $T$ is a tree, $D_F$ is a directed acyclic graph (a polytree). If 
$E_F = \emptyset$, then $p = q$ and no operations are required. If $E_F \neq \emptyset$, $D_F$ contains at 
least one source vertex $v^* \in V$ such that its in-degree $\text{deg}^-(v^*) = 0$ and its out-degree 
$\text{deg}^+(v^*) \geq 1$ (see Proposition 1.4.2 in \cite{Bang-Jensen2002}).

By the definition of the net flow requirements, the local mass balance at any vertex $i \in V$ satisfies:
\begin{equation}
p_i - q_i = \sum_{(i, j) \in E_F} W(i, j) - \sum_{(k, i) \in E_F} W(k, i)
\end{equation}
For the source vertex $v^*$, since $\text{deg}^-(v^*) = 0$, the second summation vanishes, yielding:
\begin{equation}
p_{v^*} = q_{v^*} + \sum_{(v^*, w) \in E_F} W(v^*, w)
\end{equation}
Given that $q_{v^*} \geq 0$ and $W(v^*, w) > 0$ for all $(v^*, w) \in E_F$, it follows that $p_{v^*} \geq W(v^*, w)$ for any chosen outgoing edge $(v^*, w)$.

Select an arbitrary outgoing edge $e^* = (v^*, w) \in E_F$. We execute a localized mass transfer on $e^*$ with parameter $\delta = W(v^*, w)$. The updated population vector $p'$ is given by:
\begin{align}
p'_{v^*} &= p_{v^*} - W(v^*, w) \geq 0 \\
p'_w &= p_w + W(v^*, w)
\end{align}
and $p'_i = p_i$ for all $i \in V \setminus {v^*, w}$.

The edge $e^*$ is then removed from $T$, partitioning the graph into two disjoint subtrees, $T_1$ and $T_2$. Because $e^*$ was the unique bridge between $T_1$ and $T_2$, and the mass transferred matches the required net flow across $e^*$, the restriction of $p'$ and $q$ to each subtree satisfies internal mass conservation:
\begin{equation}
\sum_{i \in V(T_1)} p'_i = \sum_{i \in V(T_1)} q_i \quad \text{and} \quad \sum_{i \in V(T_2)} p'_i =
\sum_{i \in V(T_2)} q_i
\end{equation}
Furthermore, for any remaining edge $e \in E \setminus \{e^*\}$, the partition of $T \setminus \{e^*\}$ 
does not alter the vertex sets of the components isolated by $e$ within $T_1$ or $T_2$. Thus, the required 
net flows $f'_{u \to v}$ computed with respect to $p'$ on the subtrees $T_1$ and $T_2$ are equal to 
the original values $f_{u \to v}$ for all remaining edges. The induction hypothesis applies independently 
to $T_1$ and $T_2$, ensuring the existence of a valid sequence of operations on the remaining edges.
Concatenating the transfer on $e^*$ with these two edge-disjoint sequences yields a valid sequence of
at most $1+|E(T_1)|+|E(T_2)|=|E|$ operations, at most one per edge, transforming $p$ into $q$. This
concludes the proof.
\end{proof}

\begin{lemma}[Time-Invariance of Inverse Givens Rotations]\label{LEM_Givens_inverse}
Let $R_{jk}(\theta, \phi)$ be a phased Givens rotation of angle $\theta \in [0, \frac{\pi}{2}]$ and phase $\phi \in [0, 2\pi)$ acting on the subspace $\text{Span}(\phi_j^B, \phi_k^B)$. If $R_{jk}(\theta, \phi)$ can be generated by the relaxed dynamics in an interaction time $\tau_{jk}(\theta)$, then its unitary inverse $R_{jk}(\theta, \phi)^{-1}$ can be generated within the identical interaction time $\tau_{jk}(\theta)$.
\end{lemma}

\begin{proof}
By definition of the phased Givens rotation, the unitary inverse 
satisfies:
\begin{equation}
    R_{jk}(\theta, \phi)^{-1} = R_{jk}(-\theta, \phi) = R_{jk}(\theta, \phi + \pi).
\end{equation}
According to Proposition \ref{PRO_Time_bound_Gevens}, the interaction time required to 
synthesize a rotation depends linearly on the population rotation 
angle $\theta$, evaluating to:
\begin{equation}
    \tau_{jk}(\theta) = \frac{\theta}{|A_{jk}| \cos\left(\frac{\pi}{p+1}\right)}.
\end{equation}
Since shifting the operational phase parameter from $\phi$ to $\phi + \pi$ leaves the modulus $\theta$ invariant, the temporal cost to synthesize the inverse operator is unchanged.
\end{proof}

\begin{proposition}[State-to-State Transfer Bound with Rational Rank]
\label{PRO_Result_Constructive_State_to_State_detailed}
Let $T$ be a spanning tree of the coupling graph. Let $\mathcal{E}_{free}$ be a connected subtree of $T$ such that the spectral gaps associated with its edges are linearly independent over $\mathbb{Q}$. 

Let $\mathcal{E}_{cons}$ be the set of the remaining edges in $T$. The minimal state-to-state transfer time $\tau_{\mathrm{SS}}$ is bounded from above by:
\begin{equation}
    \tau_{\mathrm{SS}}^{\max} \le \frac{\pi}{2} \sum_{e \in \mathcal{E}_{free}} \frac{\alpha_e}{|A_e|} + \pi \sum_{e \in \mathcal{E}_{cons}} \frac{\alpha_e}{|A_e|}.
\end{equation}
The numbers $\alpha_e$ share the definitions provided in Proposition \ref{PRO_Result_Constructive_State_to_State}.
\end{proposition}

\begin{proof}
The trajectory is constructed in the extended system combining the relaxed dynamics (controls in $CC(t)$)
and impulsive controls (instantaneous jumps $\exp(vB)$ generated by $B$), so that no error term is
propagated through the routing.
The degenerate case $\mathcal{E}_{free}=\emptyset$ (i.e.\ $\eta=\pi$) reduces to a transfer performed entirely on
$\mathcal{E}_{cons}=T$: Steps 2 and 3 below are then vacuous and the bound is the Step-1 quantity
$\tau_{cons}=\pi\sum_{e\in T}\alpha_e/|A_e|$. We henceforth assume $\mathcal{E}_{free}\neq\emptyset$.

\paragraph{Step 1: Local Gather and Scatter on $\mathcal{E}_{cons}$.}
The removal of $\mathcal{E}_{free}$ partitions the edges of $\mathcal{E}_{cons}$ into a forest of disjoint 
subtrees $\{T_m\}$, each intersecting $\mathcal{E}_{free}$ at a local root $u_m$. The population transfers 
(Gather towards $u_m$, then Scatter from $u_m$) operate on these subtrees. Within the relaxed system, 
Proposition \ref{PRO_partial_edge_isolation} provides access to constant generators $M^{(j,k)}_{part}$ 
allowing the synthesis of phased Givens rotations. The macroscopic physical time required to traverse 
each edge of $\mathcal{E}_{cons}$ twice (once in the Gather phase and once in the reverse Scatter phase,
whose identical cost is given by Lemma~\ref{LEM_Givens_inverse}) evaluates to:
\begin{equation}
    \tau_{cons} = 2 \sum_{e \in \mathcal{E}_{cons}} C_e = \pi \sum_{e \in \mathcal{E}_{cons}} 
    \frac{\alpha_e}{|A_e|}.
\end{equation}

\paragraph{Step 2: Impulsive Phase Alignment on the Core Tree.}
Following Step 1, the probability mass is localized on the vertices of the connected tree 
$\mathcal{E}_{free}$. The extended system allows the application of impulsive controls (Dirac measures) 
generating instantaneous jumps of the form $\exp(s B)$ for $s \in \mathbb{R}$. Because the spectral gaps 
associated with the edges of $\mathcal{E}_{free}$ are $\mathbb{Q}$-linearly independent, 
Kronecker's theorem states that the closure of the one-parameter subgroup $\exp(\mathbb{R}B)$ coincides 
with the maximal torus on the relative phases of $\mathcal{E}_{free}$. Within this extended framework, the 
state can transition to any target relative phase configuration on $\mathcal{E}_{free}$ with a 
 zero physical time cost ($T_{imp} = 0$).

\paragraph{Step 3: Core Tree Routing.}
With the masses localized and phased on $\mathcal{E}_{free}$, the populations are routed across the 
edges of $\mathcal{E}_{free}$ using the generators of the relaxed dynamics. The physical time cost to 
execute  one Givens rotation per edge in $\mathcal{E}_{free}$ is:
\begin{equation}
    \tau_{free} = \sum_{e \in \mathcal{E}_{free}} C_e = \frac{\pi}{2} \sum_{e \in \mathcal{E}_{free}} 
    \frac{\alpha_e}{|A_e|}.
\end{equation}

\paragraph{Step 4: Conclusion via Strong Limits.}
The concatenation defines an exact trajectory of the extended system reaching the target state $\psi_f$ 
in time $\tau = \tau_{cons} + \tau_{free}$, hence $\tau^{\max}_{\mathrm{ASS}} \le \tau$ by definition of the
relaxed attainable set; Lemma~\ref{LEM_approx_cont} gives $\tau^{\max}_{\mathrm{SS}} = \tau^{\max}_{\mathrm{ASS}} \le \tau$.
\end{proof}

\begin{proof}[Proof of Proposition \ref{PRO_Result_Constructive_State_to_State}]
The statement follows from Proposition \ref{PRO_Result_Constructive_State_to_State_detailed} with $\mathcal{E}_{free}=T$ for the case $\eta=\frac{\pi}{2}$ and with $\mathcal{E}_{free}=\emptyset$ for $\eta = \pi$. 
\end{proof}

\paragraph{Example.} To illustrate the mechanics of Proposition \ref{PRO_Result_Constructive_State_to_State_detailed}, we consider a state-to-state transfer in dimension $d=5$. We operate in the eigenbasis of the control matrix $B$. The active spanning tree of the coupling graph consists of edges $(1,2)$, $(2,3)$, $(3,4)$, and $(4,5)$. The objective is to steer the initial state $\psi_0=(1/2,1/2,0,1/2,1/2)^T$ to the target state $\psi_f=(-1/2,1/2,0,i/2,1/2)^T$.

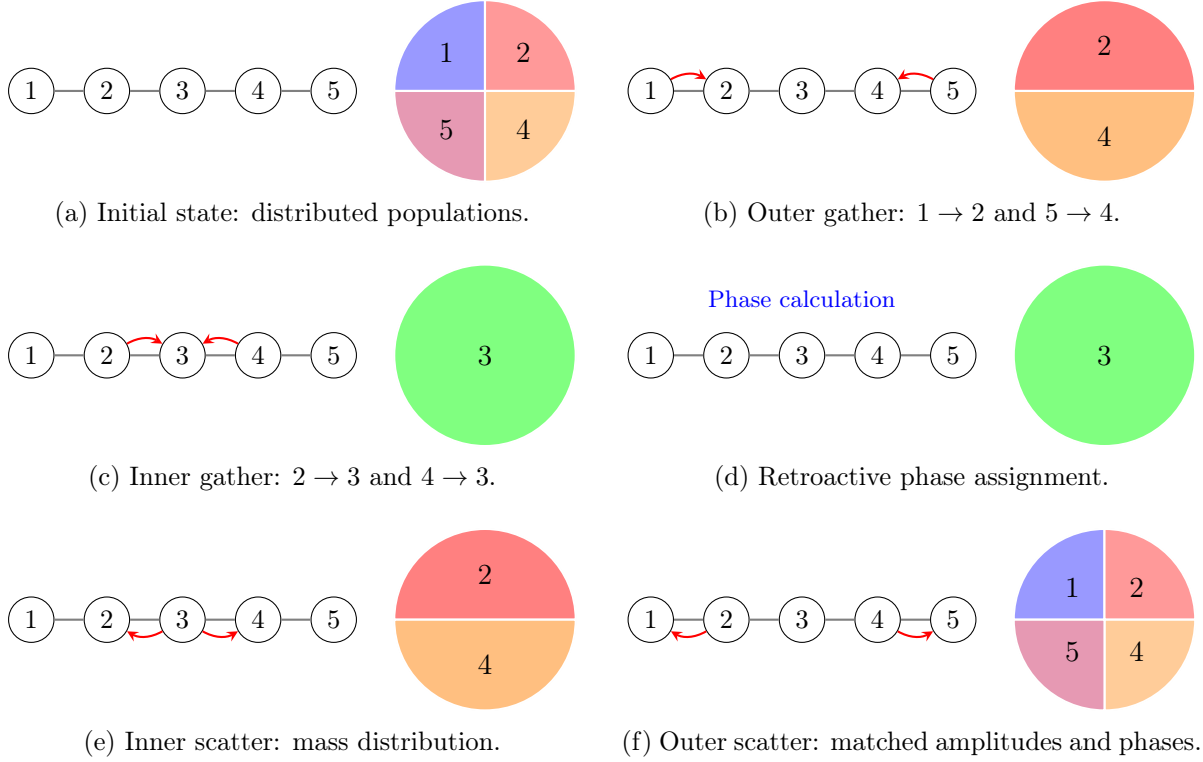
\begin{figure}[htbp]
\centering
\tikzset{
    node style/.style={circle, draw, minimum size=0.6cm, inner sep=1pt, font=\small},
    active edge/.style={->, thick, red, >=stealth},
    inactive edge/.style={-, gray, thick},
    pie slice/.style={draw=white, thick}
}

\newcommand{\drawroutingstep}[3]{
    \begin{subfigure}[t]{0.48\textwidth} 
        \centering
        \begin{tikzpicture}
            \begin{scope}[xshift=-2cm, yshift=0cm]
                \foreach \i in {1,...,5} {
                    \node[node style] (N\i) at (\i-3, 0) {\i};
                }
                \draw[inactive edge] (N1) -- (N2) (N2) -- (N3) (N3) -- (N4) (N4) -- (N5);
                #1 
            \end{scope}
            
            \begin{scope}[xshift=2cm, yshift=0cm, scale=1.2]
                \draw[gray!20, fill=gray!10] (0,0) circle (1); 
                #2 
            \end{scope}
        \end{tikzpicture}
        \caption{#3}
    \end{subfigure}
}

\drawroutingstep{
}
{
    \fill[blue!40, pie slice] (0,0) -- (90:1) arc (90:180:1) -- cycle; \node at (135:0.6) {1};
    \fill[red!40, pie slice] (0,0) -- (0:1) arc (0:90:1) -- cycle; \node at (45:0.6) {2};
    \fill[orange!40, pie slice] (0,0) -- (270:1) arc (270:360:1) -- cycle; \node at (315:0.6) {4};
    \fill[purple!40, pie slice] (0,0) -- (180:1) arc (180:270:1) -- cycle; \node at (225:0.6) {5};
}{Initial state: distributed populations.}
\hfill
\drawroutingstep{
    \draw[active edge] (N1) to[bend left=30] (N2);
    \draw[active edge] (N5) to[bend right=30] (N4);
}
{
    \fill[red!50, pie slice] (0,0) -- (0:1) arc (0:180:1) -- cycle; \node at (90:0.5) {2};
    \fill[orange!50, pie slice] (0,0) -- (180:1) arc (180:360:1) -- cycle; \node at (270:0.5) {4};
}{Outer gather: $1 \to 2$ and $5 \to 4$.}

\vspace{0.4cm}

\drawroutingstep{
    \draw[active edge] (N2) to[bend left=30] (N3);
    \draw[active edge] (N4) to[bend right=30] (N3);
}
{
    \fill[green!50, pie slice] (0,0) circle (1); \node at (0,0) {3};
}{Inner gather: $2 \to 3$ and $4 \to 3$.}
\hfill
\drawroutingstep{
    \node[above=0.2cm of N3, font=\footnotesize, color=blue] {Phase calculation};
}
{
    \fill[green!50, pie slice] (0,0) circle (1); \node at (0,0) {3};
}{Retroactive phase assignment.}

\vspace{0.4cm}

\drawroutingstep{
    \draw[active edge] (N3) to[bend left=30] (N2);
    \draw[active edge] (N3) to[bend right=30] (N4);
}
{
    \fill[red!50, pie slice] (0,0) -- (0:1) arc (0:180:1) -- cycle; \node at (90:0.5) {2};
    \fill[orange!50, pie slice] (0,0) -- (180:1) arc (180:360:1) -- cycle; \node at (270:0.5) {4};
}{Inner scatter: mass distribution.}
\hfill
\drawroutingstep{
    \draw[active edge] (N2) to[bend left=30] (N1);
    \draw[active edge] (N4) to[bend right=30] (N5);
}
{
    \fill[blue!40, pie slice] (0,0) -- (90:1) arc (90:180:1) -- cycle; \node at (135:0.5) {1};
    \fill[red!40, pie slice] (0,0) -- (0:1) arc (0:90:1) -- cycle; \node at (45:0.5) {2 };
    \fill[orange!40, pie slice] (0,0) -- (270:1) arc (270:360:1) -- cycle; \node at (315:0.5) {4};
    \fill[purple!40, pie slice] (0,0) -- (180:1) arc (180:270:1) -- cycle; \node at (225:0.5) {5};
}{Outer scatter: matched amplitudes and phases.}

\caption{State-to-state transfer on a fully resonant coupling graph. \textbf{(a)} The initial state. \textbf{(b)-(c)} The Gather sequence progressively concentrates the probability mass onto the topological root (node 3), operating independently of the initial internal phases. \textbf{(d)} Because the graph is fully resonant, no physical impulsive drift is executed; the required control phases for the subsequent steps are determined analytically by backward propagation. \textbf{(e)} The Scatter sequence distributes the populations back to the nodes. \textbf{(f)} The retroactive assignment of control phases ensures the final state matches the target complex amplitudes.}
\label{FIG_state_to_state_resonant}
\end{figure}

We first consider a fully resonant coupling graph, a scenario occurring if all spectral gaps $b_k - b_j$ are integers. As detailed in Figure \ref{FIG_state_to_state_resonant}, the routing requires a global Gather-Scatter sequence spanning the entire tree. Because the system lacks rational independence, no physical impulsive drift is executed; the required phase adjustments are entirely handled analytically via retroactive assignment of the control parameters during the Scatter phase. The population distributions at each step are represented by pie charts, where the area of each sector is proportional to the squared modulus of the corresponding coordinate.

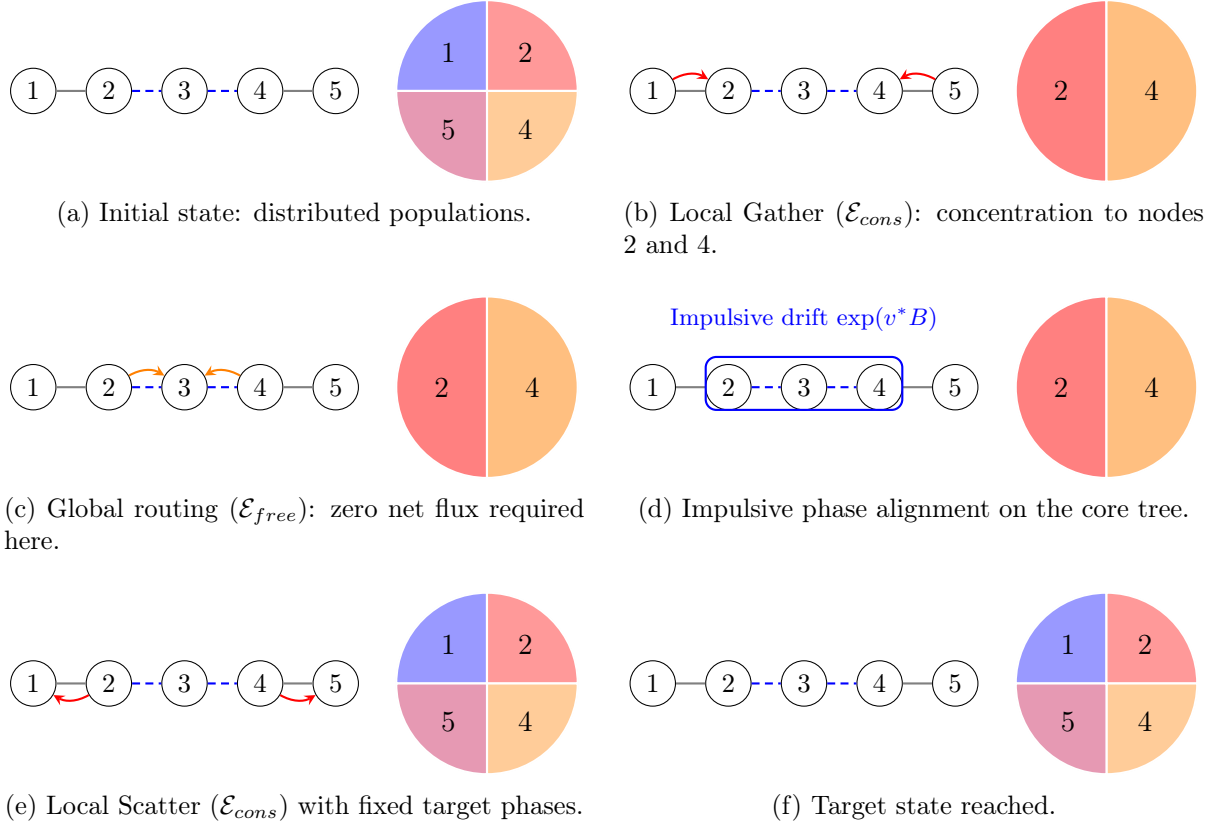
\begin{figure}[htbp]
\centering
\tikzset{
    node style/.style={circle, draw, minimum size=0.6cm, inner sep=1pt, font=\small},
    cons edge/.style={-, gray, thick},
    free edge/.style={-, blue, thick, densely dashed}, 
    active cons/.style={->, red, thick, >=stealth},
    active free/.style={->, orange, thick, >=stealth},
    pie slice/.style={draw=white, thick}
}

\newcommand{\drawroutingstep}[3]{
    \begin{subfigure}[t]{0.48\textwidth}
        \centering
        \begin{tikzpicture}
            \begin{scope}[xshift=-2cm, yshift=0cm]
                \foreach \i in {1,...,5} {
                    \node[node style] (N\i) at (\i-3, 0) {\i};
                }
                \draw[cons edge] (N1) -- (N2);
                \draw[free edge] (N2) -- (N3) (N3) -- (N4);
                \draw[cons edge] (N4) -- (N5);
                #1 
            \end{scope}
            
            \begin{scope}[xshift=2cm, yshift=0cm, scale=1.2]
                \draw[gray!20, fill=gray!10] (0,0) circle (1); 
                #2 
            \end{scope}
        \end{tikzpicture}
        \caption{#3}
    \end{subfigure}
}

\drawroutingstep{
}
{
    \fill[blue!40, pie slice] (0,0) -- (90:1) arc (90:180:1) -- cycle; \node at (135:0.6) {1};
    \fill[red!40, pie slice] (0,0) -- (0:1) arc (0:90:1) -- cycle; \node at (45:0.6) {2};
    \fill[orange!40, pie slice] (0,0) -- (270:1) arc (270:360:1) -- cycle; \node at (315:0.6) {4};
    \fill[purple!40, pie slice] (0,0) -- (180:1) arc (180:270:1) -- cycle; \node at (225:0.6) {5};
}{Initial state: distributed populations.}
\hfill
\drawroutingstep{
    \draw[active cons] (N1) to[bend left=30] (N2);
    \draw[active cons] (N5) to[bend right=30] (N4);
}
{
    \fill[red!50, pie slice] (0,0) -- (90:1) arc (90:270:1) -- cycle; \node at (180:0.5) {2};
    \fill[orange!50, pie slice] (0,0) -- (270:1) arc (270:450:1) -- cycle; \node at (0:0.5) {4};
}{Local Gather ($\mathcal{E}_{cons}$): concentration to nodes 2 and 4.}

\vspace{0.4cm}

\drawroutingstep{
    \draw[active free] (N2) to[bend left=30] (N3);
    \draw[active free] (N4) to[bend right=30] (N3);
}
{
    \fill[red!50, pie slice] (0,0) -- (90:1) arc (90:270:1) -- cycle; \node at (180:0.5) {2};
    \fill[orange!50, pie slice] (0,0) -- (270:1) arc (270:450:1) -- cycle; \node at (0:0.5) {4};
}{Global routing ($\mathcal{E}_{free}$): zero net flux required here.}
\hfill
\drawroutingstep{
    \node[above=0.3cm of N3, font=\footnotesize, color=blue] {Impulsive drift $\exp(v^* B)$};
    \draw[blue, thick, rounded corners] (-1.3, -0.3) rectangle (1.3, 0.4); 
}
{
    \fill[red!50, pie slice] (0,0) -- (90:1) arc (90:270:1) -- cycle; \node at (180:0.5) {2};
    \fill[orange!50, pie slice] (0,0) -- (270:1) arc (270:450:1) -- cycle; \node at (0:0.5) {4};
}{Impulsive phase alignment on the core tree.}

\vspace{0.4cm}

\drawroutingstep{
    \draw[active cons] (N2) to[bend left=30] (N1);
    \draw[active cons] (N4) to[bend right=30] (N5);
}
{
    \fill[blue!40, pie slice] (0,0) -- (90:1) arc (90:180:1) -- cycle; \node at (135:0.6) {1};
    \fill[red!40, pie slice] (0,0) -- (0:1) arc (0:90:1) -- cycle; \node at (45:0.6) {2};
    \fill[orange!40, pie slice] (0,0) -- (270:1) arc (270:360:1) -- cycle; \node at (315:0.6) {4};
    \fill[purple!40, pie slice] (0,0) -- (180:1) arc (180:270:1) -- cycle; \node at (225:0.6) {5};
}{Local Scatter ($\mathcal{E}_{cons}$) with fixed target phases.}
\hfill
\drawroutingstep{
}
{
    \fill[blue!40, pie slice] (0,0) -- (90:1) arc (90:180:1) -- cycle; \node at (135:0.6) {1};
    \fill[red!40, pie slice] (0,0) -- (0:1) arc (0:90:1) -- cycle; \node at (45:0.6) {2};
    \fill[orange!40, pie slice] (0,0) -- (270:1) arc (270:360:1) -- cycle; \node at (315:0.6) {4};
    \fill[purple!40, pie slice] (0,0) -- (180:1) arc (180:270:1) -- cycle; \node at (225:0.6) {5};
}{Target state reached.}

\caption{State-to-state transfer on a coupling graph with a rationally independent core. \textbf{(a)} The 
initial state. \textbf{(b)} The Local Gather sequence concentrates the probability mass onto the local 
roots (nodes 2 and 4) using  the rationally dependent edges $\mathcal{E}_{cons}$. \textbf{(c)} 
Populations are routed across the independent core tree $\mathcal{E}_{free}$. In this example, no routing 
is necessary.  \textbf{(d)} A physical impulsive drift $\exp(v^* B)$ is executed in zero macroscopic time. 
Exploiting the linear independence of the core edges, this aligns the macroscopic relative phases to their 
required target values. \textbf{(e)} The Local Scatter sequence redistributes the mass within the dependent 
subtrees. \textbf{(f)} The final target state is exactly synthesized.}
\label{FIG_state_to_state_non_resonant}
\end{figure}

Secondly, we consider the same graph topology but assume the central spectral gaps $b_3-b_2$ and $b_4-b_3$ 
are linearly independent over $\mathbb{Q}$ (for instance, evaluating to $\sqrt{2}$ and $\sqrt{3}$), while 
the peripheral gaps $b_2-b_1$ and $b_5-b_4$ are rationally dependent on them. This segregates the graph, 
defining an independent core tree $\mathcal{E}_{free} = \{(2,3), (3,4)\}$ and a dependent set 
$\mathcal{E}_{cons} = \{(1,2), (4,5)\}$. As illustrated in Figure \ref{FIG_state_to_state_non_resonant}, 
the Gather phase is now local, concentrating the populations onto the local roots (nodes 2 and 4). 
Because the macroscopic masses required on these nodes match the initial distributions, the global routing 
across $\mathcal{E}_{free}$ requires zero net flux. A physical impulsive drift is then executed to align 
the macroscopic relative phases on the core. Finally, a local Scatter phase distributes the amplitudes with 
the  internal target phases.

\subsection{Constructive Bounds for Unitary Synthesis (Propagators)}\label{SEC_Contructive_bounds_Propagators}

Synthesizing a full unitary operator $X_f \in G$ requires sequentially routing the orthonormal basis while matching $d^2-1$ relative phases. 

\begin{proposition}[Tree-Based Unitary Synthesis Bound]
\label{PRO_Result_Explicit_Propagator_precised}
Let $T$ be a spanning tree of the coupling graph. Let $\mathcal{E}_{slow}$ be a connected subtree whose corresponding spectral gaps are linearly independent over $\mathbb{Q}$.
The topological weight $w(e)$ is defined as $\min(n_1(e), n_2(e))$, where $n_1(e)$ and $n_2(e)$ are the sizes of the two connected components formed by removing $e$. 

Let $\mathcal{E}_{fast}$ contain the remaining edges of $T$. The minimal synthesis time $\tau_{\mathrm{US}}$ is bounded from above by:
\begin{equation}
    \tau_{\mathrm{US}}^{\max} \le \frac{\pi}{2} \sum_{e \in \mathcal{E}_{slow}} w(e) \frac{\alpha_e}{|A_e|} + \frac{\pi}{2} \sum_{e \in \mathcal{E}_{fast}} (w(e) + 1) \frac{\alpha_e}{|A_e|}.
\end{equation}
\end{proposition}

\begin{remark}
Proposition \ref{PRO_Result_Explicit_Propagator_precised} actually holds for any admissible bipartition
\(\mathcal{E}= \mathcal{E}_{\mathrm{slow}}\cup\mathcal{E}_{\mathrm{fast}}\) of the tree $T$
that satisfies the connectedness and Diophantine conditions.
Consequently, one may tailor the partition in order to sharpen the resulting
time‑estimate.

A convenient strategy is to assign the larger
penalty \(w(e)+1\) to the fast edges and keep the baseline penalty
\(w(e)\) for the slow edges.  Practically, the subtree \(\mathcal{E}_{\mathrm{slow}}\) is chosen so
as to maximize its unpenalized cost
\(
\sum_{e\in\mathcal{E}_{\mathrm{slow}}} \,\frac{\alpha_{e}}{|A_{e}|}.
\)
\end{remark}

\begin{proof}
The protocol is built in the relaxed (extended) system, as for Proposition~
\ref{PRO_Result_Constructive_State_to_State_detailed}. The mechanism, inspired by
\cite{PhysRevLett.73.58} (see also \cite{Clements:16}), routes the columns of $X_f^{-1}$
to the canonical eigenbasis of $B$ by a QR-like Givens decomposition, pruning the tree
from the leaves inward. We work in the relaxed sense; the transfer of the bound from
$\tau^{\max}_{\mathrm{AUS}}$ to $\tau^{\max}_{\mathrm{US}}$ is the content of the section
(Section~\ref{SUB_SEC_equivalence_worst_case}) and is recalled at the end.

\emph{Amplitude routing.}
Processing the tree leaf by leaf, each edge $e$ is used as a transit channel
$w(e)=\min(n_1(e),n_2(e))$ times, the $w(e)$-th sweep decoupling the leaf. Within the
relaxed system, Proposition~\ref{PRO_partial_edge_isolation} yields the exact synthesis
of phased Givens rotations, so the $w(e)$ sweeps are pure $SU(2)$ amplitude transits.

\emph{Phase management.}
A single impulsive drift $\exp(vB)$ is applied at the very end, producing a phase jump
$\delta_j$ on every state $j$. Because the gaps of $\mathcal{E}_{slow}$ are $\mathbb{Q}$-linearly
independent, this jump is free on the slow states: choosing $v$ sets their final phases to
the target at no extra cost. On the fast states the jump $\delta_j$ is instead imposed by
this same $v$; to compensate it, we keep track of a phase vector and, when a fast leaf $j$
is decoupled, apply the pre-compensating phase $\theta_j$ through a full $U(2)$ lock so that
$\theta_j + \delta_j$ equals the required phase on $j$. Because $v$ (hence each jump
$\delta_j=b_j\,v$) is fixed in advance by the slow-state target phases, every pre-compensation
$\theta_j$ is computed before the routing sweeps are executed
(Corollary~\ref{COR_SU2_phase_locking}).
This extra $U(2)$ sweep is the $+1$ on the fast edges; the slow edges carry no such
penalty, their phase being settled by the terminal drift.

\emph{Time cost.}
The $w(e)$ $SU(2)$ transits cost $w(e)\,C_e$; the terminal $U(2)$ lock on a fast edge adds one
sweep. Hence
\begin{equation}
    \tau = \sum_{e \in \mathcal{E}_{slow}} w(e)\,C_e + \sum_{e \in \mathcal{E}_{fast}}
    (w(e)+1)\,C_e ,
\end{equation}
which, with $C_e = \frac{\pi \alpha_e}{2|A_e|}$, gives
\begin{equation}
    \tau = \frac{\pi}{2} \sum_{e \in \mathcal{E}_{slow}} w(e) \frac{\alpha_e}{|A_e|}
    + \frac{\pi}{2} \sum_{e \in \mathcal{E}_{fast}} (w(e)+1) \frac{\alpha_e}{|A_e|}.
\end{equation}
The relaxed trajectory attains the target in time $\tau$, hence $\tau^{\max}_{\mathrm{AUS}} \le \tau$
by definition of the relaxed attainable set; the worst-case relaxed-to-exact equivalence
(Lemma~\ref{LEM_approx_cont}) then gives $\tau^{\max}_{\mathrm{US}} \le \tau$.
\end{proof}

\paragraph{Example.} To detail the QR-like factorization introduced in Proposition \ref{PRO_Result_Explicit_Propagator_precised}, we consider the synthesis of a unitary 
operator $X_f$ on the same linear graph in dimension $d=5$. We assume the mixed 
spectral configuration from the previous example 
(Figure \ref{FIG_state_to_state_non_resonant}): the central gaps $b_3-b_2$ and $b_4-b_3$ 
are rationally independent over $\mathbb{Q}$, defining the core 
$\mathcal{E}_{slow} = \{(2,3), (3,4)\}$, while the peripheral gaps $b_2-b_1$ and 
$b_5-b_4$ are rationally dependent, defining $\mathcal{E}_{fast} = \{(1,2), (4,5)\}$. 

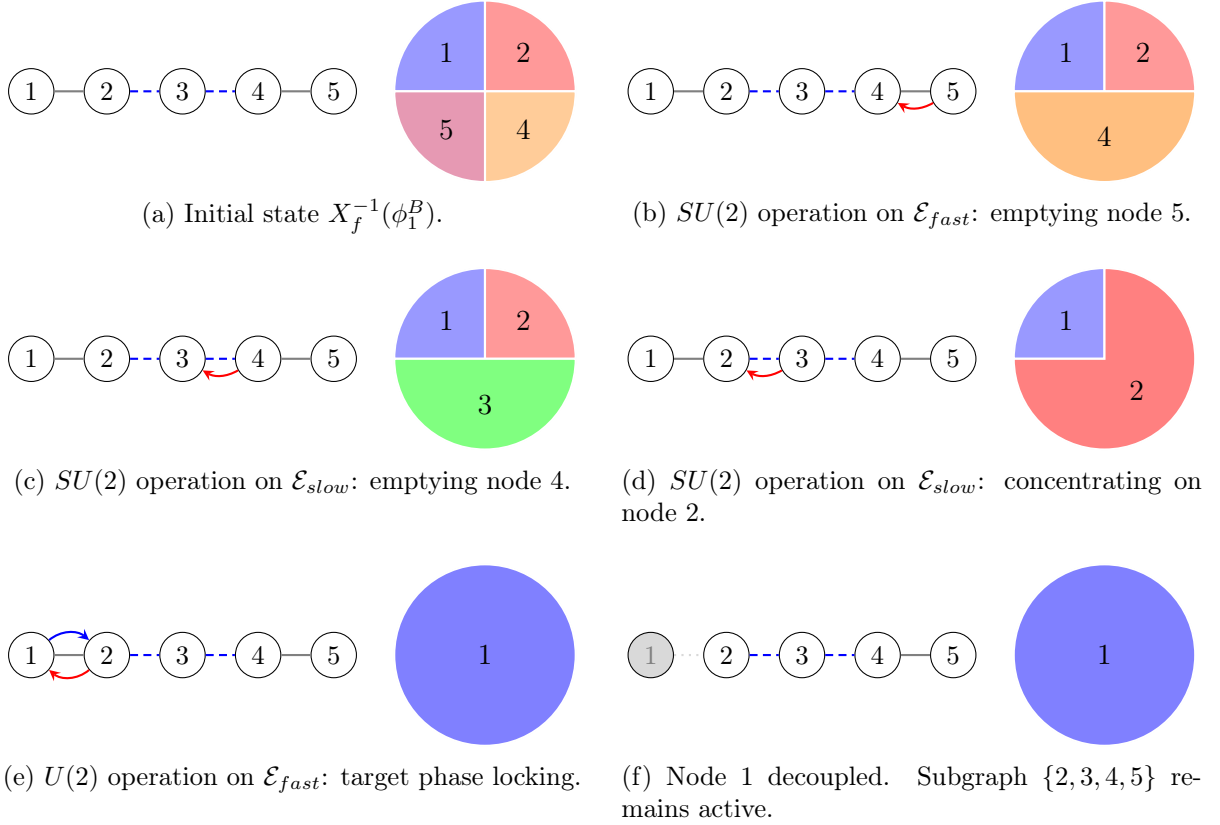
\begin{figure}[htbp] 
\centering
\tikzset{
    node style/.style={circle, draw, minimum size=0.6cm, inner sep=1pt, font=\small},
    decoupled node/.style={circle, draw, fill=gray!30, minimum size=0.6cm, inner sep=1pt, font=\small, text=gray},
    fast edge/.style={-, gray, thick}, 
    slow edge/.style={-, blue, thick, densely dashed}, 
    active edge/.style={->, red, thick, >=stealth},
    return edge/.style={->, blue, thick, >=stealth}, 
    decoupled edge/.style={-, gray!30, thick, dotted},
    pie slice/.style={draw=white, thick}
}

\newcommand{\drawroutingstep}[3]{
    \begin{subfigure}[t]{0.48\textwidth}
        \centering
        \begin{tikzpicture}
            \begin{scope}[xshift=-2cm, yshift=0cm]
                #1
            \end{scope}
            
            \begin{scope}[xshift=2cm, yshift=0cm, scale=1.2]
                \draw[gray!20, fill=gray!10] (0,0) circle (1); 
                #2
            \end{scope}
        \end{tikzpicture}
        \caption{#3}
    \end{subfigure}
}

\newcommand{\basegraph}{
    \foreach \i in {1,...,5} { \node[node style] (N\i) at (\i-3, 0) {\i}; }
    \draw[fast edge] (N1) -- (N2);
    \draw[slow edge] (N2) -- (N3) (N3) -- (N4);
    \draw[fast edge] (N4) -- (N5);
}


\drawroutingstep{
    \basegraph
}
{
    \fill[blue!40, pie slice] (0,0) -- (90:1) arc (90:180:1) -- cycle; \node at (135:0.6) {1};
    \fill[red!40, pie slice] (0,0) -- (0:1) arc (0:90:1) -- cycle; \node at (45:0.6) {2};
    \fill[orange!40, pie slice] (0,0) -- (270:1) arc (270:360:1) -- cycle; \node at (315:0.6) {4};
    \fill[purple!40, pie slice] (0,0) -- (180:1) arc (180:270:1) -- cycle; \node at (225:0.6) {5};
}{Initial state $X_f^{-1}(\phi_1^B)$.}
\hfill
\drawroutingstep{
    \basegraph
    \draw[active edge] (N5) to[bend left=30] (N4);
}
{
    \fill[blue!40, pie slice] (0,0) -- (90:1) arc (90:180:1) -- cycle; \node at (135:0.6) {1};
    \fill[red!40, pie slice] (0,0) -- (0:1) arc (0:90:1) -- cycle; \node at (45:0.6) {2};
    \fill[orange!50, pie slice] (0,0) -- (180:1) arc (180:360:1) -- cycle; \node at (270:0.5) {4};
}{$SU(2)$ operation on $\mathcal{E}_{fast}$: emptying node 5.}

\vspace{0.4cm}


\drawroutingstep{
    \basegraph
    \draw[active edge] (N4) to[bend left=30] (N3);
}
{
    \fill[blue!40, pie slice] (0,0) -- (90:1) arc (90:180:1) -- cycle; \node at (135:0.6) {1};
    \fill[red!40, pie slice] (0,0) -- (0:1) arc (0:90:1) -- cycle; \node at (45:0.6) {2};
    \fill[green!50, pie slice] (0,0) -- (180:1) arc (180:360:1) -- cycle; \node at (270:0.5) {3};
}{$SU(2)$ operation on $\mathcal{E}_{slow}$: emptying node 4.}
\hfill
\drawroutingstep{
    \basegraph
    \draw[active edge] (N3) to[bend left=30] (N2);
}
{
    \fill[blue!40, pie slice] (0,0) -- (90:1) arc (90:180:1) -- cycle; \node at (135:0.6) {1};
    \fill[red!50, pie slice] (0,0) -- (180:1) arc (180:450:1) -- cycle; \node at (315:0.5) {2};
}{$SU(2)$ operation on $\mathcal{E}_{slow}$: concentrating on node 2.}

\vspace{0.4cm}


\drawroutingstep{
    \basegraph
    \draw[active edge] (N2) to[bend left=40] (N1); 
    \draw[return edge] (N1) to[bend left=40] node[above, font=\footnotesize, color=blue] {} (N2); 
}
{
    \fill[blue!50, pie slice] (0,0) circle (1); 
    \node at (0,0) {1};
}{$U(2)$ operation on $\mathcal{E}_{fast}$: target phase locking.}
\hfill
\drawroutingstep{
    \node[decoupled node] (N1) at (-2, 0) {1};
    \foreach \i in {2,...,5} { \node[node style] (N\i) at (\i-3, 0) {\i}; }
    \draw[decoupled edge] (N1) -- (N2);
    \draw[slow edge] (N2) -- (N3) (N3) -- (N4);
    \draw[fast edge] (N4) -- (N5);
}
{
    \fill[blue!50, pie slice] (0,0) circle (1); 
    \node at (0,0) {1};
}{Node 1 decoupled. Subgraph $\{2,3,4,5\}$ remains active.}

\caption{Sequential QR-like factorization steps for the first column. \textbf{(a)} Initial distributed state. \textbf{(b)-(d)} Amplitude routing via $SU(2)$ operations, sequentially emptying the leaves inward. \textbf{(e)} $U(2)$ operation on $\mathcal{E}_{fast}$. Unlike internal transits, this terminal Givens rotation on the $(1,2)$ subspace must simultaneously annihilate node 2 and lock the pre-assigned phase on node 1. This requires a full $U(2)$ synthesis (illustrated by the bidirectional arrows), representing an accumulated rotation angle of $\pi$. This $+\pi/2$ penalty over a standard transit fully consumes the $2C_e$ time budget allocated to dependent edges. \textbf{(f)} The target node is decoupled from the active network, matching its intended amplitude and pre-assigned phase. The exact target phase will be reached after the global impulsive alignment on the core network.} 
\label{FIG_explication_synthese_unitaire}
\end{figure}

Figure \ref{FIG_explication_synthese_unitaire} depicts the first macroscopic step of 
the synthesis: routing the column vector $X_f^{-1}(\phi_1^B)$ to the canonical basis 
state $\phi_1^B$ (node 1). During the internal routing sweeps (from node 5 down to node 2), the 
Givens rotations act as pure amplitude transits within $SU(2)$,  independent of 
whether the traversed edge belongs to $\mathcal{E}_{fast}$ or $\mathcal{E}_{slow}$. The 
arbitrary residual phases left on these intermediate nodes are irrelevant at this stage 
since they remain part of the active network. The time cost of synthesizing 
a full $U(2)$ operation—allowing a maximum accumulated rotation angle of $\pi$ instead 
of $\pi/2$—is paid on the terminal edge of the sweep, which here is 
$(1,2) \in \mathcal{E}_{fast}$. This specific operation is the only one tasked with 
simultaneously annihilating the coordinate on node 2 and locking the exact target phase 
on node 1 right before its permanent decoupling from the network (Corollary~\ref{COR_SU2_phase_locking}).

\begin{remark}[Connection to Weighted Median Transmission]\label{REM_Connection_Transmission}
The combinatorial sum appearing in the upper bound can be reformulated using classical tree invariants. Let $c \in V(T)$ be the unweighted topological centroid of the tree $T$, which by definition minimizes the maximum component size upon removal. It is a standard result in graph theory (see, e.g., \cite{zbMATH03328106}) that in a tree, the centroid coincides with the median graph, which minimizes the sum of distances to all other vertices.

For any edge $e \in E(T)$, the connected component not containing $c$ has a cardinality  equal to $\min(n_1(e), n_2(e))$. By defining the effective weighted distance between any vertex $v$ and the centroid $c$ as 
\[\text{dist}_W(v, c) = \sum_{e \in \mathcal{P}(v,c)} \frac{\alpha_e}{|A_e|},\] 
where $\mathcal{P}(v,c)$ denotes the unique path linking $v$ to $c$ in $T$, the bound satisfies the algebraic identity:
\begin{equation}
    \sum_{e \in T} \min(n_1(e), n_2(e)) \frac{\alpha_e}{|A_e|} = \sum_{v \in V(T)} \text{dist}_W(v, c).
\end{equation}
In graph theory, the right-hand side is the weighted transmission (or status) of the median node.
The identity rewrites the bound of Proposition~\ref{PRO_Result_Explicit_Propagator_precised}:
routing all basis states through the topological centroid of the tree makes the total
synthesis time bounded above by $\frac{\pi}{2}$ times the sum of the weighted path lengths
from all vertices to this centroid, plus the $\frac{\pi}{2}$ contribution of the $+1$
penalty on the fast edges.
\end{remark}

\begin{proof}[Proof of Proposition \ref{PRO_Result_Explicit_Propagator}]
We consider the worst-case scenario for Proposition \ref{PRO_Result_Explicit_Propagator_precised}, where all ratios  $\frac{\alpha_e}{|A_e|}$ are identically equal to their maximum for every  $e$ in $T$, and the rational rank of the spectral gaps is equal to one, implying  $\mathcal{E}_{slow}=\emptyset$. Following Remark \ref{REM_Connection_Transmission} and results from \cite{zbMATH03328106} and \cite{zbMATH03515496}, the sum
 $\sum_{e \in T} \min(n_1(e), n_2(e)) \frac{\alpha_e}{|A_e|} = \sum_{v \in V(T)} \text{dist}_W(v, c)$ 
is maximized for a chain tree, for which
$\sum_{e \in T} \min(n_1(e), n_2(e)) = \lfloor d^2/4 \rfloor$. 

We obtain
\[ \tau_{\mathrm{US}}^{\max}\leq \frac{\pi}{2} \left (\frac{d^2}{4} + (d-1) \right ) \max_{e \in T} \frac{\alpha_e}{|A_e|} \leq \pi \frac{d^2 +4 d}{8} \max_{e \in T} \frac{\alpha_e}{|A_e|}
\] 
which is the claimed bound.
\end{proof}

\subsection{Small time local controllability}\label{SUBSEC_STLC}

We now prove the small-time local controllability near the identity that
underlies the pointwise equality and continuity of
Proposition~\ref{PRO_egalite}.  Under Assumption~\ref{ASS_drift_absorption}
(drift absorption) the relaxed attainable set contains a neighbourhood of
the identity (Proposition~\ref{PRO_relaxed_Nbhd}); combining this with the
closure of the exact attainable set and the LARC condition yields the same
neighbourhood for the exact system (Proposition~\ref{PRO_exact_Nbhd}).
Finally, the drift absorption is itself generic
(Proposition~\ref{PRO_Generic_Attainable_Neighborhood}).

\begin{proposition}[Relaxed small-time controllability near the identity]\label{PRO_relaxed_Nbhd}
Assume $d\geq 3$ and Assumptions~\ref{ASS_control} (simple spectrum of $B$,
$\mathrm{Lie}(A,B)\supset\mathfrak{g}$), \ref{ASS_coupling_graph}
(coupling graph connected) and \ref{ASS_drift_absorption}.
Then, for every $\varepsilon>0$, the relaxed attainable set
$\mathcal{A}^{\mathrm{rel}}_{\leq\varepsilon}$ contains an operator-norm
neighborhood of the identity $I_d$.
\end{proposition}

\begin{proof}
Fix $\varepsilon>0$. We show that $\mathcal{A}^{\mathrm{rel}}_{\leq\varepsilon}$
contains a ball $B_{op}(I_d,\rho)$ for some $\rho>0$ depending only on
$\varepsilon$, $A$, $B$, and the coupling graph.

Recall that we have defined the generalized solutions of \eqref{EQ_main} (Definition
\ref{DEF_Generalized_solution}) write as $Z(t)=R(t) Y^{\Omega}_{t,0}$. 
Following \eqref{eq:decomp} and defining $W^{\tilde{\Omega}}(t)= \exp(-t A_D ) 
Y^{\Omega}_{t,0}$, $Z(t)$ also factorizes as $Z(t)= R(t) \exp( t A_D) 
W^{\tilde{\Omega}}(t)$.
 
By Corollary \ref{COR_STLC_W}, there exists $\rho>0$ such that
\begin{equation}\label{EQ_relax_Nbhd_etape_1}
B(I_d,\rho) \subset \bigcup_{0 \leq t \leq \varepsilon} \quad \bigcup_{{\Omega}\in L^\infty([0,t],CC)} \{\exp(-t A_D ) 
Y^{\Omega}_{t,0} \}.\end{equation}

By Assumption~\ref{ASS_drift_absorption}, for every $t$ in $[0,\varepsilon]$, for every $\Omega$ in $L^\infty([0,t],CC)$, 
\begin{equation} \label{EQ_relax_Nbhd_etape_2}
\exp(-t A_D ) 
Y^{\Omega}_{t,0}  \in \overline{\exp(\R B)} Y^{\Omega}_{t,0} \subset \mathcal{A}_{\leq t}^{\mathrm{rel}} \subset \mathcal{A}_{\leq \varepsilon}^{\mathrm{rel}}.
\end{equation}
 Combining \eqref{EQ_relax_Nbhd_etape_1} and \eqref{EQ_relax_Nbhd_etape_2}, we get
 $B(I_d,\rho) \subset\mathcal{A}_{\leq \varepsilon}^{\mathrm{rel}}$ which finishes the proof. 
\end{proof}

\begin{proposition}[Exact small-time controllability near the identity]\label{PRO_exact_Nbhd}
Assume $d \geq 3$ and Assumptions~\ref{ASS_control} (simple spectrum of $B$,
$\operatorname{Lie}(A,B)\supset\mathfrak{g}$), \ref{ASS_coupling_graph}
(coupling graph connected) and \ref{ASS_drift_absorption}.
Then, for every $\varepsilon>0$, the exact attainable set
$\mathcal{A}_{\leq\varepsilon}$ contains an operator-norm neighborhood of
the identity $I_d$.
\end{proposition}

\begin{proof}
We show: $\forall\,\varepsilon>0$, $\exists\,\rho(\varepsilon)>0$ such that
$B_{op}(I_d,\rho(\varepsilon))\subset\mathcal{A}_{\leq\varepsilon}$.
Fix $\varepsilon>0$. Throughout, the concatenation convention is the
left-invariant one: applying the control that reaches $h$ in time
$t_1$ followed by the control that reaches $c$ in time $t_2$ yields the
total propagator $c\cdot h$ (the later propagator multiplies on the left).

\emph{Step 1 (relaxed neighborhood).}
By Proposition~\ref{PRO_relaxed_Nbhd}, there exists $\delta_1>0$ such that
\begin{equation}\label{EQ_rel_Nbhd}
B_{op}(I_d,\delta_1)\;\subset\;\mathcal{A}^{\mathrm{rel}}_{\leq\varepsilon/2}.
\end{equation}

\emph{Step 2 (density of the exact in the relaxed).}
By Lemma~\ref{LEM_cloture_ensemble_attaignable},
$\mathcal{A}^{\mathrm{rel}}_{\leq\varepsilon/2}
=\overline{\mathcal{A}_{\leq\varepsilon/2}}$ (closure in $\|\cdot\|_{op}$).
Combining with \eqref{EQ_rel_Nbhd}:
\begin{equation}\label{EQ_density_in_ball}
 B_{op}(I_d,\delta_1)\;\subset\;\overline{\mathcal{A}_{\leq\varepsilon/2}},
\end{equation}
i.e., $\mathcal{A}_{\leq\varepsilon/2}$ is dense in $ B_{op}(I_d,\delta_1)$.

\emph{Step 3 (local controllability of the exact system).}
The conclusion is a classical consequence of the LARC condition (see Lemmas 6.1 and 6.3 in \cite{JURDJEVIC1972313}). We include the proof for completeness.

By Assumption~\ref{ASS_control}
($\operatorname{Lie}(A,B)\supset\mathfrak{g}$) and the Krener theorem, the initial condition $I_d$ belongs to the topological closure of the interior of the attainable
set $\mathcal{A}_{\leq \varepsilon/2}$. Let us note $V_1$ a non-empty open set contained in
$\mathcal{A}_{\leq \varepsilon/2} \cap B_{op}(I_d,\delta_1) $. Then the set $V=V_1^{-1}$ of the inverse (in the matrix sense) of the elements of $V_1$, is open and contained in  $ B_{op}(I_d,\delta_1)$. By density, there exists $X_a \in V \cap \mathcal{A}_{\leq \varepsilon/2}$.     The set $V_1 X_a$ is open because the right multiplication by $X_a$ is a homeomorphism. By construction, $I_d$ belongs to $V_1 X_a$ and 
$V_1 X_a \subset \mathcal{A}_{\leq \varepsilon}$. Hence $\mathcal{A}_{\leq \varepsilon}$ is a neighborhood of $I_d$ which concludes the proof. 
\end{proof}

\begin{proposition}[Generically satisfied drift absorption]\label{PRO_Generic_Attainable_Neighborhood}
Under Assumptions \ref{ASS_control} and \ref{ASS_coupling_graph}, if the
spectral gaps of $B$ have maximal rational rank (which
holds for a full-measure set of $B\in\mathfrak{g}$ with simple
spectrum), then Assumption~\ref{ASS_drift_absorption} holds:
\[
\exp(\mathbb{R}\,A_D)\;\subset\;\overline{\exp(\mathbb{R}\,B)}.
\]
\end{proposition}

\begin{proof}
Since $A_D$ is diagonal in the eigenbasis of $B$,
$\exp(tA_D)\in\mathcal{T}_B$ (the maximal diagonal torus in $G$) for all
$t\in\mathbb{R}$. By Kronecker's theorem, the rational independence of
the spectral gaps implies that the one-parameter subgroup
$\exp(\mathbb{R}B)=\{e^{i\theta B}:\theta\in\mathbb{R}\}$ is \emph{dense}
in $\mathcal{T}_B$. Hence
\[
\exp(\mathbb{R}\,A_D)\;\subset\;\mathcal{T}_B
\;=\;\overline{\exp(\mathbb{R}\,B)},
\]
which is \eqref{EQ_drift_absorption}.
\end{proof}

\begin{remark}
The drift absorption condition Assumption~\ref{ASS_drift_absorption} is generically satisfied. Notice that it does not require the spectral gaps of $B$ to have a rational maximal rank.  For instance, if $B$ has an integer spectrum (making $\exp(\mathbb{R}B)$ a 1D closed curve), Assumption~\ref{ASS_drift_absorption} still holds provided $A_D$ is linearly dependent on $B$ (e.g., $A_D = 0$ or $A_D \propto B$).
\end{remark}

\section{Application: Spectral Symmetry Breaking, example of the Planar Rotor}\label{SEC_Chiral}

\subsection{The Planar Rotor Model and Finite-Dimensional Truncation}
\label{SEC_model_planar}

The rotational dynamics of a quantum planar rotor controlled by an external linearly polarized electric field are governed by the time-dependent Schrödinger equation on the unit circle $\mathbb{S}^1$. In atomic units, the evolution of the wave function $\psi(\theta, t)$ is given by:
\begin{equation}
    i \partial_t \psi(\theta, t) = \left( -\Delta + u(t)\cos(\theta) \right) \psi(\theta, t), \quad \theta \in [-\pi, \pi],
\end{equation}
where $\Delta = \frac{\partial^2}{\partial \theta^2}$ is the Laplace-Beltrami operator on $\mathbb{S}^1$, and $u(t) \in \mathbb{R}$ acts as the scalar control input representing the time-varying amplitude of the external field. 

We project the system onto the spectral basis of the free Hamiltonian. The orthonormal eigenbasis of the Laplacian is defined by the Fourier modes $ \phi_k^A(\theta) = \frac{1}{\sqrt{2\pi}} e^{ik\theta}$ for $k \in \mathbb{Z}$, satisfying $-\Delta \phi_k^A = k^2 \phi_k^A$.

We perform a Galerkin truncation of the state space, restricting the dynamics to the lowest 
$2d+1$ momentum eigenstates, $\mathcal{H}_d = \text{span}\{\phi_k^A\}_{k=-d}^{d}$. We 
denote by $\Pi_d$ the orthogonal projection on the $D=2d+1$-dimensional space $
\mathcal{H}_d$. Applying this projection to the Schrödinger equation yields a finite-
dimensional   bilinear control system (with size: $D\times D$) of the form:
\begin{equation}
    \frac{d}{dt} \mathbf{c}(t) = \left( A + u(t)B \right) \mathbf{c}(t),
\end{equation}
where $\mathbf{c}(t) \in \mathbb{C}^{2d+1}$ is the vector of probability amplitudes. The matrices $A$ and $B$, representing respectively the unperturbed drift and the control generator, are anti-Hermitian. In this spectral basis, the drift matrix $A = -i \Pi_d (-\Delta) \Pi_d$ is diagonal:
\begin{equation}
    A_{j,k} = -i k^2 \delta_{j,k}, \quad j,k \in \{-d, \dots, d\}.
\end{equation}
The interaction potential $\cos(\theta) = \frac{1}{2}(e^{i\theta} + e^{-i\theta})$ acts as a shift operator in momentum space:
\[\cos(\theta) \phi_k^A=\cos(\theta) e^{i k \theta} = \frac{1}{2} e^{i(k+1)\theta} + \frac{1}{2}e^{i(k-1)\theta}=\frac{1}{2} \phi_{k+1}^A + \frac{1}{2} \phi_{k-1}^A.\] Consequently, the $D \times D$ control matrix $B = -i \Pi_d (\cos\theta) \Pi_d$ exhibits a pure tridiagonal structure with zero diagonal entries:
\begin{equation}
    B_{j,k} = -\frac{i}{2} \left( \delta_{j, k+1} + \delta_{j, k-1} \right).
\end{equation}

This computational framework implies a symmetric spectrum for $A$ and a pure Toeplitz structure for $B$:
$$A=-i \begin{pmatrix}
  d^2 & 0 & \cdots  & \cdots & \cdots & 0 \\
 0 	& (d-1)^2 	& \ddots 	&  	& 	& \vdots\\
\vdots & \ddots & \ddots & \ddots  & & \vdots \\
\vdots &  &\ddots   & 0 & \ddots & \vdots \\
\vdots &  &   & \ddots  &\ddots & 0\\
0 & \cdots & \cdots & \cdots  & 0 & d^2
\end{pmatrix}, \quad 
B=-\begin{pmatrix}
 0 & i/2 & 0  & \cdots & 0 \\
 i/2 & 0 & \ddots & \ddots & \vdots\\
0 & \ddots & \ddots & \ddots & 0 \\
\vdots & \ddots & \ddots & 0 & i/2 \\
0 & \cdots & 0 & i/2  & 0
\end{pmatrix}$$

However, the spatial reflection symmetry of the dipole interaction, invariant under $\theta \mapsto -\theta$, suggests an orthogonal decomposition of $L^2(\mathbb{S}^1)$ into symmetric and antisymmetric subspaces. By applying a passive unitary transformation to the real Fourier basis—composed of the normalized cosine states $\phi_k^c(\theta) \propto \cos(k\theta)$ for $k \in \{0, \dots, d\}$ and sine states $\phi_k^s(\theta) \propto \sin(k\theta)$ for $k \in \{1, \dots, d\}$—the matrices undergo a block-diagonalization. 

In this parity-adapted basis, the system explicitly decouples into two independent invariant subspaces, represented by the ordered direct sums $A = A_{\text{even}} \oplus A_{\text{odd}}$ and $B = B_{\text{even}} \oplus B_{\text{odd}}$. The unperturbed drift matrices remain diagonal:
$$ A_{\text{even}} = -i \, \mathrm{diag}(0, 1, 4, \dots, d^2), \quad A_{\text{odd}} = -i \, \mathrm{diag}(1, 4, \dots, d^2). $$
The control generator $B$ forbids transitions between the two parity sectors. Its non-zero intra-block elements are confined to the adjacent off-diagonals:
$$ B_{\text{even}} = -i \begin{pmatrix} 
0 & \frac{1}{\sqrt{2}} & 0 & \cdots & 0 \\ 
\frac{1}{\sqrt{2}} & 0 & \frac{1}{2} & \ddots & \vdots \\ 
0 & \frac{1}{2} & 0 & \ddots & 0 \\ 
\vdots & \ddots & \ddots & 0 & \frac{1}{2} \\ 
0 & \cdots & 0 & \frac{1}{2} & 0 
\end{pmatrix}, \quad 
B_{\text{odd}} = -i \begin{pmatrix} 
0 & \frac{1}{2} & 0 & \cdots & 0 \\ 
\frac{1}{2} & 0 & \frac{1}{2} & \ddots & \vdots \\ 
0 & \frac{1}{2} & 0 & \ddots & 0 \\ 
\vdots & \ddots & \ddots & 0 & \frac{1}{2} \\ 
0 & \cdots & 0 & \frac{1}{2} & 0 
\end{pmatrix} $$

While this parity-adapted basis makes the decoupling explicit, the altered 
boundary condition at $k=0$ (the $1/\sqrt{2}$ coupling term) breaks the uniform 
Toeplitz structure of the symmetric block. We therefore return to the uniform exponential 
basis $e^{ik\theta}$, which keeps $B$ Toeplitz — a fact used below to extract the 
closed-form eigenvalues and to solve the Diophantine resonance equations.

\subsection{Spectral Chiral Symmetry and Effective Couplings}
\label{sec:chiral_structure}

As established  in Appendix \ref{SEC_Tridiagonal}, the eigenvalues $\lambda_m$ and the components of the corresponding orthonormal eigenvectors $\phi_m^B$ of the $D\times D$ control matrix $B$ are given by the closed-form expressions:
\begin{align}
    \lambda_m &= -i \cos\left( \frac{m \pi}{D + 1} \right), \quad m \in \{1, \dots, D=2d+1\}, \label{eq:eigenvalues} \\
    \langle \phi_j^A , \phi_m^B \rangle &= \sqrt{\frac{2}{D+1}} \sin\left( \frac{m j \pi}{D+1} \right), \quad j \in \{1, \dots, D\}. \label{eq:eigenvectors}
\end{align}

The spectrum exhibits a  chiral symmetry. For any mode $m \in \{1, \dots, d\}$, the conjugate mode $\bar{m} = D + 1 - m=2d+2 -m $ satisfies:
\begin{equation}
    \lambda_{\bar{m}} = -i \cos\left( \frac{(D+1-m) \pi}{D + 1} \right) = -i \cos\left( \pi - \frac{m \pi}{D + 1} \right) = -\lambda_m.
\end{equation}
This chiral involution $m \leftrightarrow \bar{m}$ partitions the spectrum into two symmetric branches (a positive one and a negative one), separated by the zero mode $m = (D+1)/2 = d+1$. 

Consequently, any transition between adjacent momentum states $(m, m+1)$ has a degenerate "chiral clone" transition $(\bar{m}-1, \bar{m}) = (D-m, D+1-m)$. Lemma~\ref{LEM_Global_Non_Resonance} below shows that this 2-fold chiral degeneracy is the only exact resonance among adjacent spectral gaps in the odd subspace (the subspace used in the proof of Proposition~\ref{PRO_rotateur_introduction}); for specific values of $D$, parasitic resonances occur in the even subspace, and the lemma characterizes these exceptions arithmetically.

Classifying and ruling out parasitic resonances reduces to solving trigonometric Diophantine equations over rational multiples of $\pi$; the theorem of Conway and Jones (Theorem~\ref{TH_liste_Conway_Jones}) classifies the rational linear relations involving up to four such cosines.

\begin{theorem}[Theorem 7 in \cite{Conway1976}] \label{TH_liste_Conway_Jones}
Suppose we have at most four distinct rational multiples of $\pi$ lying strictly between $0$ and $\pi/2$ for which some rational linear combination of their cosines is rational but no proper subset has this property. Then the appropriate linear combination is proportional to one from the following list:
\begin{align}
    &\cos\left(\frac{\pi}{3}\right) = \frac{1}{2} \label{eq:CJ_1} \\
    -&\cos(\phi) + \cos\left(\frac{\pi}{3}-\phi\right) + \cos\left(\frac{\pi}{3}+\phi\right) = 0 \quad (0<\phi <\pi/6) \label{eq:CJ_param} \\
    &\cos\left(\frac{\pi}{5}\right) - \cos\left(\frac{2\pi}{5}\right) = \frac{1}{2} \label{eq:CJ_2+2} \\
    &\cos\left(\frac{\pi}{7}\right) - \cos\left(\frac{2\pi}{7}\right) + \cos\left(\frac{3\pi}{7}\right) = \frac{1}{2} \label{eq:CJ_42} \\
    &\cos\left(\frac{\pi}{5}\right) - \cos\left(\frac{\pi}{15}\right) + \cos\left(\frac{4\pi}{15}\right) = \frac{1}{2} \label{eq:CJ_30A} \\
    -&\cos\left(\frac{2\pi}{5}\right) + \cos\left(\frac{2\pi}{15}\right) - \cos\left(\frac{7\pi}{15}\right) = \frac{1}{2} \label{eq:CJ_30B} \\
    &\cos\left(\frac{\pi}{7}\right) + \cos\left(\frac{3\pi}{7}\right) - \cos\left(\frac{\pi}{21}\right) + \cos\left(\frac{8\pi}{21}\right) = \frac{1}{2} \label{eq:CJ_4_1} \\
    &\cos\left(\frac{\pi}{7}\right) - \cos\left(\frac{2\pi}{7}\right) + \cos\left(\frac{2\pi}{21}\right) - \cos\left(\frac{5\pi}{21}\right) = \frac{1}{2} \label{eq:CJ_4_2} \\
    -&\cos\left(\frac{2\pi}{7}\right) + \cos\left(\frac{3\pi}{7}\right) + \cos\left(\frac{4\pi}{21}\right) + \cos\left(\frac{10\pi}{21}\right) = \frac{1}{2} \label{eq:CJ_4_3} \\
    -&\cos\left(\frac{\pi}{15}\right) + \cos\left(\frac{2\pi}{15}\right) + \cos\left(\frac{4\pi}{15}\right) - \cos\left(\frac{7\pi}{15}\right) = \frac{1}{2} \label{eq:CJ_4_4}
\end{align}
\end{theorem}

\begin{lemma}[Arithmetic Non-Resonance of the Parity-Preserving Ladder]
\label{LEM_Global_Non_Resonance}
Let $N=D+1=2d+2$, and let $c_s=\cos(s\pi/N)$. Consider the parity-preserving transition graph, i.e. transitions $(j,k)$ with $j\equiv k\pmod 2$.

\begin{enumerate}
    \item \textbf{Inter-branch resonances.}
    The chiral transitions $(m,N-m)$ are unique up to the trivial identity, for every even $N$.

    \item \textbf{Vertical intra-branch resonances.}
    The vertical transitions $(m,m+2)$ are unique up to chiral symmetry, except for the following exact resonances:
    \begin{itemize}
        \item If $N=6q$ with $q\ge 3$, the four transitions $(3q-2, 3q)$, $(3q, 3q+2)$, $(q-2, q+2)$, and $(5q-2, 5q+2)$ are  co-resonant with frequency gap $\Delta = \sin(\pi/3q)$.
        \item If $N=30$,  the four transitions $(8,10)$, $(20,22)$, $(2,6)$, and $(24,28)$ are  co-resonant with frequency gap  $c_8 - c_{10} = c_2 - c_6=\cos \left ( \frac{\pi}{15} \right ) - \cos \left ( \frac{\pi}{5} \right )$. The four transitions  $(13,15)$, $(15,17)$, $(3,7)$ and $(23,27)$ are  co-resonant with frequency gap $ \cos\left(\frac{\pi}{10}\right) - \cos\left(\frac{7\pi}{30}\right)$.
    \end{itemize}
    In both cases, the companion transition has index gap $4$.

    \item \textbf{Harmonic resonances.}
    If $\gcd(N,105)=1$, then no harmonic resonance
    \[
        p(c_j-c_k)=c_l-c_m,\qquad p\ge 2,
    \]
    exists between two non-degenerate transitions of the parity-preserving graph.
\end{enumerate}
\end{lemma}

\begin{remark}
Restricting the analysis to the specific subspace $(A_o, B_o)$, spanned by the even spectral indices $m \in \{2, 4, \dots, 2d\}$, the exact intra-branch resonances reduce to the following two cases:
\begin{enumerate}
    \item \textbf{Case $N = 12r$ ($r \ge 2$, $d = 6r-1$):} The quadruplet formed by the adjacent vertical edges $(6r-2, 6r)$ and $(6r, 6r+2)$, together with the companion transitions $(2r-2, 2r+2)$ and $(10r-2, 10r+2)$, is strictly co-resonant. Here, the central mode $6r = N/2$ is an even index, and therefore explicitly belongs to the ladder.
    \item \textbf{Case $N = 30$ ($d = 14$):} The quadruplet formed by the vertical edges $(8,10)$ and $(20,22)$, together with the companion transitions $(2,6)$ and $(24,28)$, is strictly co-resonant.
\end{enumerate}
\end{remark}

\begin{proof}[Proof of Lemma \ref{LEM_Global_Non_Resonance}]
For each index $s \in \{1,\dots,N-1\}$, define the folded angle $\alpha_s \in (0,\pi/2]$ by $\alpha_s = \theta_s$ if $s \le N/2$, and $\alpha_s = \pi - \theta_s$ otherwise. Consequently, $c_s = \pm \cos\alpha_s$. The initial logic and phase-folding mappings establish that any resonance equation reduces to a rational linear combination of cosines of distinct angles in $(0,\pi/2]$:
\begin{equation}
    \sum_{i} \eta_i \cos\alpha_i = 0,
\end{equation}
where the collected integer coefficients $\eta_i$ take values $\pm 1$, $\pm 2$ (e.g., from self-chiral terms $c_m - c_{N-m} = 2c_m$), or include the multiplier $p \ge 2$ in the harmonic case.

By Theorem \ref{TH_liste_Conway_Jones}, the set of irreducible four-cosine relations (equations \eqref{eq:CJ_4_1} to \eqref{eq:CJ_4_4}) evaluates only to $1/2$. A vanishing relation with heterogeneous coefficient magnitudes (such as $\{2, \pm 1\}$ or $\{p, \pm 1\}$ with $p \ge 2$) cannot be irreducible. The relation must be reducible into rational sub-families.

\textbf{1. Inter-branch resonances ($p=1$):}
After folding, the inter-branch condition $2c_m = c_l - c_k$ generates a coefficient $2$, leading to an equation of the form $2\cos\alpha_m - \cos\alpha_l \pm \cos\alpha_k = 0$. By Theorem \ref{TH_liste_Conway_Jones}, irreducible relations possess uniform coefficient magnitudes (up to global scaling). An equation with coefficients $\{2, -1, \pm 1\}$ is therefore reducible. There are two potential decomposition pathways:

\textit{Case A (Symmetric split):} The coefficient 2 is split, forming $(\cos\alpha_m - 
\cos\alpha_l) + (\cos\alpha_m \pm \cos\alpha_k) = 0$. For this to reduce into non-
trivial rational sub-relations, each bracket must evaluate to $\pm q \in \mathbb{Q}^*$. 
According to Theorem \ref{TH_liste_Conway_Jones} (Eq. \ref{eq:CJ_2+2}), the only non-
trivial rational evaluation for two cosines requires the angle pair $\{\pi/5, 2\pi/5\}
$. Both brackets must independently draw from this exact angular set, which 
algebraically forces $\alpha_l = \alpha_k$. This violates the topological assumption 
that the edges involve distinct spectral indices.

\textit{Case B (Bipartite split):} The relation reduces to 
$2\cos\alpha_m = \cos\alpha_l$ (or $\cos\alpha_k$). Since $\{2, -1\}$ is not an 
irreducible Conway-Jones family, it reduces further, demanding that 
$\cos\alpha_m \in \mathbb{Q}$ independently. By Niven's theorem, the only rational 
values for the cosine of a rational multiple of $\pi$ in the interval $(0, \pi/2]$ are 
$0$ and $1/2$. Testing $\cos\alpha_m = 1/2$ yields $2(1/2) = 1 = \cos\alpha_l$ hence  
$\alpha_l = 0$, which is excluded from the domain. Testing $\cos\alpha_m = 0$ yields 
$\alpha_m = \alpha_l = \pi/2$, violating the distinct angle requirement.

\textit{Conclusion:} The decomposition must isolate the term $2\cos\alpha_m$ to 
evaluate to a rational independently, alongside $\pm \cos\alpha_l \pm \cos\alpha_k \in 
\mathbb{Q}$. By Niven's theorem, $\cos\alpha_m \in \{0, 1/2\}$. If $\cos\alpha_m = 
1/2$, then $2\cos\alpha_m = 1$, requiring $\pm \cos\alpha_l \pm \cos\alpha_k = -1$. In 
$(0, \pi/2]$, no sum or difference of two cosines can equal $-1$ without invoking the 
angle $0$. Therefore, $\cos\alpha_m= 0$, implying $\alpha_m = \pi/2$ (the central mode 
$N/2$). The relation decomposes entirely into trivial/chiral cancellations. Thus, $(m, 
N-m)$ is unique (being self-chiral).

\textbf{2. Vertical intra-branch resonances ($p=1$):}
The equation $c_m - c_{m+2} = c_j - c_k$ decomposes according to the rational evaluations provided in Theorem \ref{TH_liste_Conway_Jones}.

\textit{Case 2+2:} Deriving from \eqref{eq:CJ_2+2}, the index difference is $2N/10$. 
Setting $2N/10 = 2$, hence $N=10$ generates the pairs $(2,4)$ and $(6,8)$. These are chiral clones, creating no parasitic resonance.

\textit{Case 3+1 (sum=0):} Derived from the parametric identity \eqref{eq:CJ_param}, this requires $N=6q$, yielding $c_{2q-a} + c_{2q+a} - c_a - c_{3q} = 0$. For a vertical edge, the index difference must be 2:

- If $2(q-a) = 2$ then $a = q-1$ and the secondary pair $(3q-1, 3q)$ has an index gap of 1, violating parity.

- If $q-a = 2$ then $a = q-2$, the pairs are $(3q-2, 3q)$ and $(q-2, q+2)$. This yields the first allowed exception.

\textit{Case 3+1 (sum=1/2):} Derived from the three-cosine sums evaluating to $1/2$:

- Family 30A \eqref{eq:CJ_30A}: Since $\frac{\pi}{15}=\frac{s\pi}{N}$ for some $s$ in $\N$, $N$ is a multiple of $15$. Since $N=2d+2$ is even by construction, this implies that $N$ is a multiple of 30. Hence  $N=30r$ for some $r \in \N$ and  \eqref{eq:CJ_30A} reads 
\[c_{6r}+c_{8r}-c_{2r}-c_{10r}=0.\]
 Setting the difference to 2 requires $r=1$, hence $N=30$, yielding the resonance $(8,10) \leftrightarrow (2,6)$ with $d=14$.

- Families 30B \eqref{eq:CJ_30B} and 42 \eqref{eq:CJ_42}: $N$ must be a multiple of $30$ or $42$ by the same reasoning as for the family 30A. Equation \eqref{eq:CJ_30B} 
\[ -c_{12 r} + c_{4r} -c_{10 r} -c_{14 r} =0\]
reads, after chiral folding $c_{14 r}=-c_{16 r}$:
\[ c_{4r} -c_{10 r} = c_{12 r}-c_{16 r}.\]
Equation \eqref{eq:CJ_42} reads
 \[c_{6r}-c_{12r}=c_{14r}-c_{18r}.\]
 These two remaining families yield a minimum index difference of at least $4r\ge 4$ across 
 all possible pairings. Consequently, these algebraic relations cannot be rearranged to 
 contain an index difference of $2$. They represent isolated resonances only between 
 higher-order transitions. Because they do not involve any vertical edge $(m, m+2)$, these 
 resonances remain decoupled from the active ladder and are not activated by the routing 
 protocol.

\textbf{3. Harmonic resonances ($p \ge 2$):}
Suppose $\gcd(N, 105)=1$.  After folding, the relation $p c_j - p c_k - c_l + c_m = 0$ must reduce. Because $\operatorname{lcm}(3,5,7,15,21)= 105$ and $\gcd(N, 105)=1$, the grid contains none of the special angles in Theorem \ref{TH_liste_Conway_Jones}. The entire relation decomposes into trivial terms, forcing $c_j = c_k$ or $c_l = c_m$. The spectral gaps vanish, neutralizing any harmonic resonance on the active edges.
\end{proof}

\begin{remark}
These two exceptions concern the \emph{reduced} ladder, obtained after the
central (centroid) mode has been removed. In Section~\ref{SEC_Chiral_StS}
the centroid-depopulation result of Section~\ref{SEC_Chiral_Topological_Reduction}
(Lemma~\ref{LEM_Tip_Depopulation}) is used to eliminate the central mode.
Once this mode is gone, the $N=6q$ exception is no longer active, since the
vertical edge $(3q-2,3q)$ it involves is incident to the eliminated mode.
The $N=30$ exception $(8,10)\leftrightarrow(2,6)$, by contrast, does not
involve the central mode: it survives the reduction and must still be handled
explicitly at the protocol level (Section~\ref{SEC_Chiral_StS}).
\end{remark}

\begin{proposition}\label{prop:effective_couplings}
Let $D = 2d+1$ be the dimension of the truncated state space $\mathcal{H}_d$. In the eigenbasis of the control generator $B$, the effective drift amplitudes $\tilde{A}_{m,n} = \langle \phi_m^B, A \phi_n^B \rangle$ for the unperturbed Laplacian $A = -i \, \mathrm{diag}(k^2)$  vanish  for any pair $(m,n)$ having different parities (i.e., $\tilde{A}_{m,n} = 0$ if $m-n$ is odd). 

For $m \neq n$ with $m-n$ even, the amplitudes are given by the closed-form expression:
\begin{equation} \label{eq:A_tilde_general}
\tilde{A}_{m,n} = \frac{-i}{2} \frac{\sin\left(\frac{m\pi}{D+1}\right) \sin\left(\frac{n\pi}{D+1}\right)}{\sin^2\left(\frac{(m-n)\pi}{2(D+1)}\right) \sin^2\left(\frac{(m+n)\pi}{2(D+1)}\right)}
\end{equation}
Consequently, the amplitude of the inter-block rung connecting an eigenstate $m$ to its chiral partner 
$\bar{m} = D+1-m$ is always allowed and reduces to:
\begin{equation} \label{eq:A_tilde_chiral}
\tilde{A}_{m, \bar{m}} = -\frac{i}{2} \tan^2\left(\frac{m\pi}{D+1}\right)
\end{equation}
\end{proposition}

\begin{remark}[Parity Selection Rule and Subspace Controllability]\label{REM_Subspace_decoupling}
 The parity operator $\mathcal{P}\psi(\theta)=\psi(-\theta)$ commutes with $A$, $B$, and $\Pi_d$, so 
 $[\mathcal{P},A]=[\mathcal{P},B]=0$. Since $B$ has a non-degenerate spectrum, its eigenvectors $\phi_m^B$ 
 carry a definite (alternating) parity, and the effective drift $\tilde 
 A_{m,n}=\langle\phi_m^B,A\phi_n^B\rangle$ vanishes whenever $m-n$ is odd. This proves the 
 block-diagonalization of Proposition \ref{prop:effective_couplings} without trigonometric sums.

 Consequently, $\mathrm{Lie}(A,B)$ does not span $\mathfrak{u}(D)$, and population transfer between the 
 even ($d{+}1$ modes) and odd ($d$ modes) sectors is forbidden. The system is, however, fully controllable 
 within each sector.

 To route both sectors under the single scalar control $u(t)$, we embed the two-component forest into a 
 pseudo-spanning tree $\widetilde T$ with a zero-weight virtual edge and a zero-flux constraint on 
 that edge. The control-theoretic justification is the block-diagonal structure of the interaction-picture 
 generator $M(v,t)=e^{-vB-tA_D}A_{ND}e^{vB+tA_D}=M_{\mathrm{even}}(v,t)\oplus M_{\mathrm{odd}}(v,t)$: the 
 shared parameter $v(t)$ provides a common time scale for the convexification but generates no inter-block 
 coupling. When a Givens rotation is synthesized on an edge of one sector, the off-diagonal entries of the 
 other sector average to zero (the efficiency of the isolation is governed by the lowest harmonic resonance, which may involve cross-sector gaps; the worst case is bounded by $\alpha^\ast \leq 2$), so the 
 passive sector undergoes only a diagonal phase drift, which is bookkept and compensated at the end of the 
 protocol.

\end{remark}

\begin{proof}[Proof of Proposition \ref{prop:effective_couplings}]
Let $X = [B, [B, A]]$ be the double commutator. For $m \neq n$, the effective couplings obey the
algebraic identity $\tilde{A}_{m,n} = \frac{\langle \phi_m^B, X \phi_n^B \rangle}{(\lambda_m - \lambda_n)^2}$.
In the mapped basis $j \in \{1, \dots, D\}$, the skew-adjoint matrices are $A_{jj} = -i(j-d-1)^2$ and $B_{jk} = -\frac{i}{2}(\delta_{j,k+1} + \delta_{j,k-1})$.

Direct evaluation of the double commutator yields a bulk action and localized boundary defects: $X = X_{bulk} + \Delta X$. 
The bulk matrix is $X_{bulk} = -i I + \frac{i}{2} P$, where $P_{jk} = \delta_{j,k+2} + \delta_{j,k-2}$ is the distance-2 adjacency operator. 
The sine eigenvector $\phi_n^B(j)=\sqrt{\frac{2}{D+1}}\sin(\theta_n j)$, $\theta_n=\frac{n\pi}{D+1}$, satisfies $P\phi_n^B = 2\cos(2\theta_n)\phi_n^B$ in the bulk, whence $X_{bulk}\phi_n^B = -2i\sin^2(\theta_n)\phi_n^B$ there. At the endpoints the distance-2 rule loses one neighbour, so $(P\phi_n^B)(1)=(P\phi_n^B)(D)=2\cos(2\theta_n)\phi_n^B + \phi_n^B$. For $m\neq n$ of the same parity, with $\phi_m^B(D)\phi_n^B(D)=\phi_m^B(1)\phi_n^B(1)$, the bulk therefore contributes
\[
  \langle\phi_m^B, X_{bulk}\phi_n^B\rangle
  = \tfrac{i}{2}\big[\phi_m^B(1)\phi_n^B(1)+\phi_m^B(D)\phi_n^B(D)\big]
  = i\,\phi_m^B(1)\phi_n^B(1) = \tfrac{2i}{D+1}\sin\theta_m\sin\theta_n,
\]
not zero: $\phi_n^B$ is an eigenvector of $X_{bulk}$ only away from the ends.

The off-diagonal amplitude is driven by these bulk endpoint corrections together with the finite-truncation boundary defects $\Delta X$. 
Evaluating the action of the double commutator at the boundaries $j=1$ and $j=D$ yields $X_{1,1} = -\frac{1}{2}(A_{1,1} - A_{2,2}) = i(d - 1/2)$ and symmetrically for $X_{D,D}$. 
Subtracting the bulk diagonal contribution $X_{bulk,1,1} = -i$ isolates the defect amplitude: 
\begin{equation}
    \Delta X_{1,1} = \Delta X_{D,D} = i(d - 1/2) - (-i) = i(d + 1/2) = i\frac{D}{2}.
\end{equation}
The defect operator therefore acts as $(\Delta X \phi_n^B)(1) = i\frac{D}{2} \phi_n^B(1)$ and $(\Delta X \phi_n^B)(D) = i\frac{D}{2} \phi_n^B(D)$. 

For $m$ and $n$ of the same parity, $\phi_m^B(D)\phi_n^B(D) = \phi_m^B(1)\phi_n^B(1)$. Taking the inner product and substituting the explicit normalization $\phi_n^B(1) = \sqrt{\frac{2}{D+1}} \sin(\theta_n)$ sums the defect contributions to:
\begin{equation}
    \langle \phi_m^B, \Delta X \phi_n^B \rangle = i D \, \phi_m^B(1) \phi_n^B(1) = 2i \left(\frac{D}{D+1}\right) \sin(\theta_m) \sin(\theta_n).
\end{equation}
Summing the bulk-endpoint and defect contributions gives, for $m\neq n$ of the same parity,
$\langle\phi_m^B, X\phi_n^B\rangle = \frac{2i}{D+1}\sin\theta_m\sin\theta_n + 2i\,\frac{D}{D+1}\sin\theta_m\sin\theta_n = 2i\,\sin\theta_m\sin\theta_n$.

Substituting the eigenvalues $\lambda_m = -i \cos(\theta_m)$ into the squared spectral gap yields:
\begin{equation}
    (\lambda_m - \lambda_n)^2 = - (\cos \theta_m - \cos \theta_n)^2 = -4 \sin^2\left(\frac{\theta_m - \theta_n}{2}\right) \sin^2\left(\frac{\theta_m + \theta_n}{2}\right).
\end{equation}

Dividing the off-diagonal amplitude $2i\sin\theta_m\sin\theta_n$ by the squared gap isolates the off-diagonal matrix elements, producing Equation \eqref{eq:A_tilde_general}. For the chiral coupling $n = D+1-m$, substituting $\theta_{\bar{m}} = \pi - \theta_m$ algebraically simplifies the result to Equation \eqref{eq:A_tilde_chiral}.
\end{proof}

\begin{corollary}[Block-Diagonal Structure of the Effective Drift]
\label{cor:block_diagonal_A}
Let $\mathcal{H}_{\text{even}}$ and $\mathcal{H}_{\text{odd}}$ be the invariant symmetric and antisymmetric subspaces of dimensions $d+1$ and $d$, respectively. In the parity-adapted eigenbasis of the control generator $B = B_{\text{even}} \oplus B_{\text{odd}}$, the effective drift matrix decomposes as $\tilde{A} = \tilde{A}^{\text{even}} \oplus \tilde{A}^{\text{odd}}$. 

The transition amplitudes within each invariant block are completely determined by restricting the global solution (Equation \eqref{eq:A_tilde_general}) to states of identical parity:
\begin{enumerate}
    \item \textbf{Symmetric subspace:} The entries of $\tilde{A}^{\text{even}}$ evaluate on the odd indices $m = 2r-1$ and $n = 2r'-1$, where $r, r' \in \{1, \dots, d+1\}$.
    \item \textbf{Antisymmetric subspace:} The entries of $\tilde{A}^{\text{odd}}$ evaluate on the even indices $m = 2s$ and $n = 2s'$, where $s, s' \in \{1, \dots, d\}$.
\end{enumerate}
For any $m \neq n$ satisfying the block parity constraint, the coupling strength is given by:
\begin{equation}
    \tilde{A}^{\text{block}}_{m,n} = -i \frac{\sin\left(\frac{m\pi}{D+1}\right) \sin\left(\frac{n\pi}{D+1}\right)}{2 \sin^2\left(\frac{(m-n)\pi}{2(D+1)}\right) \sin^2\left(\frac{(m+n)\pi}{2(D+1)}\right)}.
\end{equation}
\end{corollary}

\begin{proof}
The block-diagonalization ($\tilde{A}_{m,n} = 0$ for $m-n \equiv 1 \pmod 2$) is an immediate consequence of 
the parity selection rule established in Proposition \ref{prop:effective_couplings}. The explicit 
expressions for the non-zero intra-block elements $\tilde{A}^{\text{even}}$ and $\tilde{A}^{\text{odd}}$ 
require no further calculation; they are extracted by direct restriction of the closed-form expression 
\eqref{eq:A_tilde_general} from Proposition \ref{prop:effective_couplings} to the corresponding 
parity-matched index pairs.
\end{proof}

\begin{corollary}[Chiral Invariance of Effective Couplings]
\label{COR_Chiral_Coupling}
The effective drift amplitudes exhibit chiral symmetry. For any allowed transition between modes $m$ and $n$, the coupling strength is identical to the coupling strength between their respective chiral partners $\bar{m} = D+1-m$ and $\bar{n} = D+1-n$:
\begin{equation}
    \tilde{A}_{\bar{m}, \bar{n}} = \tilde{A}_{m,n}.
\end{equation}
\end{corollary}

\begin{proof}
The proof follows by direct substitution of the chiral indices into the  formula for $\tilde{A}_{m,n}$ given in Eq. \eqref{eq:A_tilde_general}. 
For the numerator, the reflection identity of the sine function yields:
\begin{equation}
    \sin\left(\frac{\bar{m}\pi}{D+1}\right) = \sin\left(\pi - \frac{m\pi}{D+1}\right) = \sin\left(\frac{m\pi}{D+1}\right),
\end{equation}
leaving the numerator invariant.
For the denominators, the sum and difference of the chiral indices evaluate to:
\begin{align}
    \bar{m} - \bar{n} &= (D+1-m) - (D+1-n) = -(m-n), \\
    \bar{m} + \bar{n} &= 2(D+1) - (m+n).
\end{align}
Since the sine is squared, the sign inversion in the difference term is absorbed: $\sin^2\left(\frac{-(m-n)\pi}{2(D+1)}\right) = \sin^2\left(\frac{(m-n)\pi}{2(D+1)}\right)$. 
For the sum term, a phase shift of $\pi$ occurs:
\begin{equation}
    \sin^2\left(\frac{(2D+2 - (m+n))\pi}{2(D+1)}\right) = \sin^2\left(\pi - \frac{(m+n)\pi}{2(D+1)}\right) = \sin^2\left(\frac{(m+n)\pi}{2(D+1)}\right).
\end{equation}
Every trigonometric term in the closed-form expression is individually invariant under the joint chiral involution $(m,n) \mapsto (\bar{m}, \bar{n})$, establishing the equality of the effective couplings.
\end{proof}

\begin{remark}[Diophantine Accidental Resonances on the Ladder]
The Conway-Jones argument of Lemma \ref{LEM_Global_Non_Resonance} rules out exact equalities between spectral gaps of the ladder, but it does not exclude integer harmonic relations $|\omega_{j'k'}|=p\,|\omega_{jk}|$ with $p\geq 2$. Whether such a relation involves an active protocol edge depends on the arithmetic of $D+1$.

 For $D=3$ ($D{+}1=4$), the only same-parity pair is $(1,3)$, with gap $\sqrt{2}$. No harmonic multiple exists within the graph, so $M_{1,3}=+\infty$ and $\alpha_{1,3}=1$.

 For $D=5$ ($D{+}1=6$), one checks that
 \[
 |b_1-b_5| \;=\; \sqrt{3} \;=\; 2\,|b_1-b_3| \;=\; 2\,|b_3-b_5|,
 \]
 where $(1,5)$ is the chiral rung and $(1,3),(3,5)$ are verticals. The pairs $(1,3)$ and $(3,5)$ are not chiral clones of $(1,5)$, so $M_{1,3}=M_{3,5}=2$ and $\alpha_{1,3}=\alpha_{3,5}=\left|\cos(\pi/3)\right|^{-1}=2$.

 Since $D+1$ can be highly composite for arbitrary $d$ (e.g., $D+1\in\{6,8,10,12,14,15,\dots\}$) and the resonant pair may or may not coincide with an active edge of the spanning tree, we do not track $\alpha_e$ edge-by-edge. All estimates of Section \ref{SUBSEC_Topology_Reduction} are stated with the uniform bound $\alpha^*:=\max_{e}\alpha_e\leq 2$, which is the worst case admissible by Proposition \ref{PRO_Result_Constructive_State_to_State_detailed}.
\end{remark}

 \begin{lemma}[Time constants of the protocol edges] \label{LEM_time_constants_C_cross_C_vert}
  Let $\delta=\pi/(2d+2)$ and $\alpha^*=\max_e\alpha_e\leq 2$. The routing protocols of the next sections \ref{SEC_Chiral_StS} and \ref{SEC_Chiral_unitary_synthesis} use three types of edges, each with a characteristic time $C_e=\pi\alpha_e/(2|A_e|)$ (equation \eqref{EQ_valeur_C_e}):
\begin{description}
 \item \textbf{Verticals.}  $(m,m+2)$ : the innermost verticals of the ladder.
 \item \textbf{Central rung.} $(m_c-2,\,m_c+2)$, where $m_c=d+1$ is the self-chiral mode : the dedicated chiral crossing of the ladder.
 \item \textbf{Adjacent rung.}  $(m_c-3,\,m_c+3)$ : the next chiral crossing out from the center.
\end{description}
 Define
 $$
 C_{vert}:=\max_{\text{verticals}} C_e,
 \qquad
 C_{cross}:=\max_{\text{central rung}} C_e,
 \qquad
 C_{cross}':=\max_{\text{adjacent rung}} C_e.
 $$
 Then:
\begin{align*}
C_{vert}\;\leq\;\frac{D+1}{D}\frac{3\pi^3\alpha^*}{8(d+1)^2},
 & \mbox{ for } d \geq 2\\
 C_{cross}\;\leq\;\frac{D+1}{D}\frac{16\pi\alpha^*}{(d+1)^2},
 & \mbox{ for } d \geq 3\\ 
 C_{cross}'\;\leq\;\frac{D+1}{D}\frac{36\pi\alpha^*}{(d+1)^2}
  & \mbox{ for } d \geq 5\\. 
\end{align*}
\end{lemma}

\begin{proof}
These are the edge times $C_e=\pi\alpha_e/(2|A_e|)$ of \eqref{EQ_valeur_C_e}, specialized to the three edge types.

 \textbf{Verticals.} From Proposition \ref{prop:effective_couplings},
 $$ |\tilde A_{m,m+2}|^{-1} = \frac{D+1}{D} 2\sin^2(\delta) 
 \left[ 1-\frac{\sin^2(\delta)}{\sin^2((m+1)\delta)} \right]^{-1}. $$
 The bracket is maximized at $m=1$, giving $\lvert\tilde A\rvert^{-1}\leq 3\sin^2(\delta)$ for $d\geq 2$ (since 
$\delta\leq\pi/6$). With $\sin(\delta)\leq\delta$ :
$$
C_{vert} \leq \frac{\pi\alpha^*}{2}\cdot \frac{D+1}{D} 3\delta^2 
\leq \frac{D+1}{D} \frac{3\pi^3\alpha^*}{8(d+1)^2}.
 $$

\textbf{Central rung.} The rung connects $m_c-2=d-1$ to $m_c+2=d+3$, at distance $k=2$ from the spectral center. From Proposition \ref{prop:effective_couplings},
$$|\tilde A_{d-1,d+3}| = \frac{D}{D+1} \frac{1}{2}\cot^2(2\delta).$$
 Using $\tan(x)\leq\frac{4}{\pi}x$ for $x\in[0,\pi/4]$, equivalently $\cot(x)\geq\frac{\pi}{4x}$, (with $x=2 \delta \leq \frac{\pi}{4}$ hence $d\geq 3$) :
 $\displaystyle  |\tilde A_{d-1,d+3}| \geq \frac{1}{2}\left(\frac{\pi}{8\delta}\right)^2 = \frac{\pi^2}{128\delta^2}
 = \frac{(d+1)^2}{32}$, hence 
 \[  C_{cross} \leq \frac{\pi\alpha^*}{2}\cdot 
 \frac{D+1}{D}\frac{32}{(d+1)^2} = \frac{D+1}{D} \frac{16\pi\alpha^*}{(d+1)^2}. \]

\textbf{Adjacent rung.} Same argument as for the central rung,  with $k=3$ (note that
we need $3 \delta \leq \pi/4$ hence $d\geq 5$):
$$|\tilde A_{d-2,d+4}| = \frac{D}{D+1} \frac{1}{2}\cot^2(3\delta) \geq \frac{D}{D+1} \frac{(d+1)^2}{72},$$
 hence  \[C_{cross}' \leq \frac{\pi\alpha^*}{2}\cdot \frac{D+1}{D} \frac{72}{(d+1)^2} = \frac{D+1}{D} \frac{36\pi\alpha^*}{(d+1)^2}. \]
 \end{proof}

\subsection{Topological Reduction of the Parity-Adapted Transition Graphs}\label{SEC_Chiral_Topological_Reduction}
\label{SUBSEC_Topology_Reduction}

 The routing network is determined by the effective couplings of Proposition \ref{prop:effective_couplings} 
 and the non-resonance conditions of Lemma \ref{LEM_Global_Non_Resonance}. The parity selection rule 
 (Remark \ref{REM_Subspace_decoupling}) partitions it into two disjoint invariant sub-graphs. Their topology
 depends on the parity of $d$: the self-chiral mode $m_c=d+1$ belongs to exactly one sector.

\begin{figure}[htbp]
    \centering
    \begin{tikzpicture}[
        vertex/.style={circle, draw, fill=black!10, inner sep=1pt, minimum size=18pt, font=\small},
        edge/.style={thick, black!70}
    ]

        \begin{scope}[shift={(0,0)}]
            \node[vertex] (L1) at (0,0) {2};
            \node[vertex] (L2) at (0,1.5) {4};
            \node[vertex] (L3) at (0,3) {6};
            \node[vertex] (L4) at (0,4.5) {8};

            \node[vertex] (R1) at (2.5,0) {16};
            \node[vertex] (R2) at (2.5,1.5) {14};
            \node[vertex] (R3) at (2.5,3) {12};
            \node[vertex] (R4) at (2.5,4.5) {10};

            \draw[edge] (L1) -- (L2) -- (L3) -- (L4);
            \draw[edge] (R1) -- (R2) -- (R3) -- (R4);

            \draw[edge] (L1) -- (R1);
            \draw[edge] (L2) -- (R2);
            \draw[edge] (L3) -- (R3);
            \draw[edge] (L4) -- (R4);

            \node[below=0.8cm of {(1.25,0)}] {\textbf{Perfect Ladder} };
        \end{scope}

        \begin{scope}[shift={(7,0)}]
            \node[vertex] (L1) at (0,0) {1};
            \node[vertex] (L2) at (0,1.5) {3};
            \node[vertex] (L3) at (0,3) {5};
            \node[vertex] (L4) at (0,4.5) {7};

            \node[vertex] (R1) at (2.5,0) {17};
            \node[vertex] (R2) at (2.5,1.5) {15};
            \node[vertex] (R3) at (2.5,3) {13};
            \node[vertex] (R4) at (2.5,4.5) {11};

            \node[vertex, fill=black!25] (C) at (1.25, 6) {9};

            \draw[edge] (L1) -- (L2) -- (L3) -- (L4) -- (C);
            \draw[edge] (R1) -- (R2) -- (R3) -- (R4) -- (C);

            \draw[edge] (L1) -- (R1);
            \draw[edge] (L2) -- (R2);
            \draw[edge] (L3) -- (R3);
            \draw[edge] (L4) -- (R4); 

            \node[below=0.8cm of {(1.25,0)}] {\textbf{Triangle-Tipped Ladder} };
        \end{scope}

    \end{tikzpicture}
    \caption{Topological distinction between the two invariant parity sub-graphs for a truncation degree $d=8$ ($D=17$). The odd-index sub-graph (right, 9 vertices) absorbs the self-chiral spectral center $m_c=9$, sealing the ladder with a rigid 3-cycle. The even-index sub-graph (left, 8 vertices) lacks the central node, forming a perfect bipartite ladder graph. Because $D+1=18 \equiv 0 \pmod 6$, Lemma \ref{LEM_Global_Non_Resonance} applies: the odd-index sub-graph hosts an exact parasitic resonance. Specifically, the active edge $(7,9)$ resonates exactly with the unplotted transition $(1,5)$, and their respective chiral clones $(9,11)$ and $(13,17)$ similarly resonate. The even-index sub-graph (corresponding to the odd subspace generated by matrices $A_o$ and $B_o$) remains immune to this degeneracy.}
\label{fig:ladder_topology_d8}
\end{figure}
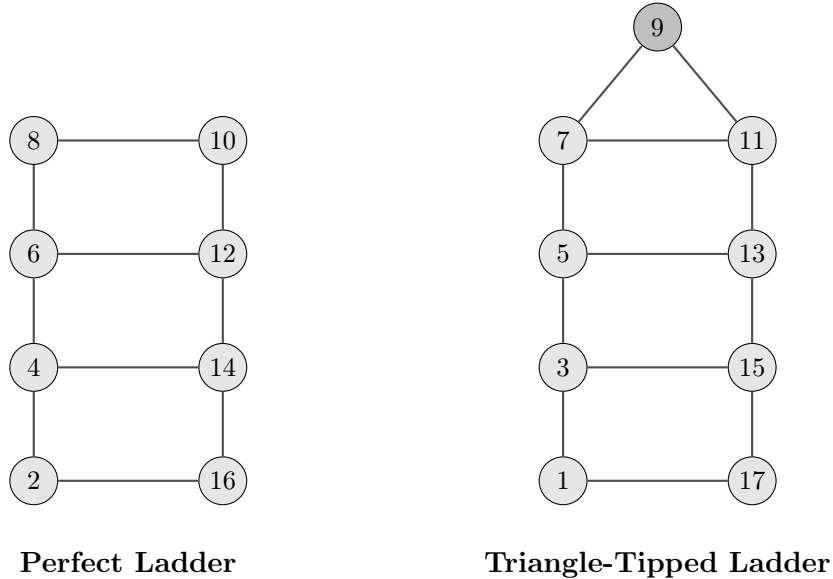

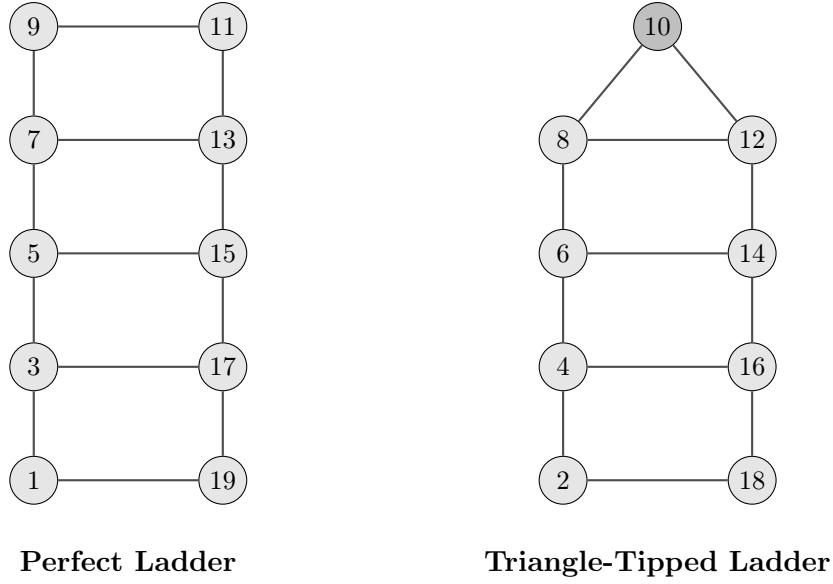
\begin{figure}[htbp]
    \centering
    \begin{tikzpicture}[
        vertex/.style={circle, draw, fill=black!10, inner sep=1pt, minimum size=18pt, font=\small},
        edge/.style={thick, black!70}
    ]

        \begin{scope}[shift={(0,0)}]
            \node[vertex] (L1) at (0,0) {1};
            \node[vertex] (L2) at (0,1.5) {3};
            \node[vertex] (L3) at (0,3) {5};
            \node[vertex] (L4) at (0,4.5) {7};
            \node[vertex] (L5) at (0,6) {9};

            \node[vertex] (R1) at (2.5,0) {19};
            \node[vertex] (R2) at (2.5,1.5) {17};
            \node[vertex] (R3) at (2.5,3) {15};
            \node[vertex] (R4) at (2.5,4.5) {13};
            \node[vertex] (R5) at (2.5,6) {11};

            \draw[edge] (L1) -- (L2) -- (L3) -- (L4) -- (L5);
            \draw[edge] (R1) -- (R2) -- (R3) -- (R4) -- (R5);

            \draw[edge] (L1) -- (R1);
            \draw[edge] (L2) -- (R2);
            \draw[edge] (L3) -- (R3);
            \draw[edge] (L4) -- (R4);
            \draw[edge] (L5) -- (R5);

            \node[below=0.8cm of {(1.25,0)}] {\textbf{Perfect Ladder} };
        \end{scope}

        \begin{scope}[shift={(7,0)}]
            \node[vertex] (L1) at (0,0) {2};
            \node[vertex] (L2) at (0,1.5) {4};
            \node[vertex] (L3) at (0,3) {6};
            \node[vertex] (L4) at (0,4.5) {8};

            \node[vertex] (R1) at (2.5,0) {18};
            \node[vertex] (R2) at (2.5,1.5) {16};
            \node[vertex] (R3) at (2.5,3) {14};
            \node[vertex] (R4) at (2.5,4.5) {12};

            \node[vertex, fill=black!25] (C) at (1.25, 6) {10};

            \draw[edge] (L1) -- (L2) -- (L3) -- (L4) -- (C);
            \draw[edge] (R1) -- (R2) -- (R3) -- (R4) -- (C);

            \draw[edge] (L1) -- (R1);
            \draw[edge] (L2) -- (R2);
            \draw[edge] (L3) -- (R3);
            \draw[edge] (L4) -- (R4); 

            \node[below=0.8cm of {(1.25,0)}] {\textbf{Triangle-Tipped Ladder} };
        \end{scope}

    \end{tikzpicture}
    \caption{Topological distinction between the two invariant parity sub-graphs for a 
    truncation degree $d=9$ ($D=19$). Because $D+1=20$ is not a multiple of 6, Lemma 
    \ref{LEM_Global_Non_Resonance} prevents parasitic resonances (besides the structural 
    chiral  resonances). The even-index sub-graph (right, 9 vertices) absorbs 
    the self-chiral spectral center $m_c=10$, sealing the ladder with a rigid 3-cycle. The 
    odd-index sub-graph (left, 10 vertices) lacks the central node, forming a perfect 
    bipartite ladder graph.}
    \label{fig:ladder_topology_d9}
\end{figure}

 The sector lacking $m_c$ is a perfect bipartite ladder on $d$ vertices, with two symmetric chiral branches 
 flat-ended at the spectral boundary. Depending on the parity of $d$, the sector containing $m_c$ has $d+1$ vertices if $d$ is even or $d$ vertices if $d$ is odd; its two branches meet 
 at $m_c$, and the adjacent nodes $d{-}1$ and $d{+}3$, which are chiral partners, together with $m_c$ form 
 a 3-cycle (a "triangular tip") closing the ladder (Figures \ref{fig:ladder_topology_d8} and 
\ref{fig:ladder_topology_d9}).

 This structural discrepancy breaks the bipartite nature of the routing graph. The routing protocol therefore begins with a preliminary depopulation step on the triangular tip.

\begin{lemma}[Centroid Depopulation on Triangular Tips]
\label{LEM_Tip_Depopulation}
Consider the local 3-level dynamics at the triangular tip of the parity sub-graph, governed by $\dot{x}(t) = \Omega(u(t))x(t)$ for $x \in \mathbb{C}^3$, where the generator is:
\begin{equation}
\Omega(u) = \begin{pmatrix} 
0 & u_1(t) & u_1(t) \\ 
-\bar{u}_1(t) & 0 & u_2(t) \\ 
-\bar{u}_1(t) & -\bar{u}_2(t) & 0 
\end{pmatrix}.
\end{equation}
Let the control cost be $J(u) = \frac{2}{\pi}C_{vert} \|u_1\|_{L^1} + \frac{2}{\pi}C_{cross} \|u_2\|_{L^1}$. For any initial state $x(0) \in \mathbb{S}^5 \subset \mathbb{C}^3$, there exists a single-pulse control sequence that depopulates the centroid state, yielding $x(T) = (0, y_2, y_3)^T$, with a bounded cost:
\begin{equation}
J(u) \le \frac{1}{\sqrt{2}} C_{vert}.
\end{equation}
\end{lemma}

\begin{proof}
The physical symmetry of the centroid couplings implies that the field $u_1$ acts identically on the transitions $1 \leftrightarrow 2$ and $1 \leftrightarrow 3$. To isolate the active dynamics, we apply the orthogonal basis transformation $\mathcal{B} = \{e_1, v_+, v_-\}$, where $v_+ = \frac{e_2+e_3}{\sqrt{2}}$ is the symmetric mode and $v_- = \frac{e_2-e_3}{\sqrt{2}}$ is the antisymmetric mode.

In this basis, the generator $\Omega(u)$ restricts the action of $u_1$ to the $(e_1, v_+)$ subspace:
\begin{equation}
\Omega(u_1, u_2=0) \equiv \begin{pmatrix} 
0 & \sqrt{2} u_1 & 0 \\ 
-\sqrt{2} \bar{u}_1 & 0 & 0 \\ 
0 & 0 & 0 
\end{pmatrix}.
\end{equation}
The control $u_2$ is  unnecessary for depopulating $e_1$ and is identically set to zero ($u_2(t) = 0$), nullifying the $C_{cross}$ cost. 

The initial state decomposes as $x(0) = x_1 e_1 + c_+ v_+ + c_- v_-$. The dynamics on the active 
$\mathrm{span}(e_1, v_+)$ subspace form an  $SU(2)$ rotation with an amplified effective coupling 
strength of $\sqrt{2}|u_1|$. For any complex amplitude pair $(x_1, c_+)$, there exists a constant phase 
$\phi = \arg(u_1)$ and a bounded rotation angle $\gamma \le \pi/2$ such that the pulse nullifies 
the $e_1$ projection. 

The required integrated pulse area evaluates to $\sqrt{2} \int_0^T |u_1| dt = \gamma \le \frac{\pi}{2}$. 
The corresponding $L^1$ norm of the control is therefore $\|u_1\|_{L^1} \le \frac{\pi}{2\sqrt{2}}$. 
Multiplying by the penalty weight $\frac{2}{\pi}C_{vert}$ yields the  upper bound 
$J(u) \le \frac{1}{\sqrt{2}} C_{vert} $. Following this operation, the state vector is fully transferred to 
the $(e_2, e_3)$ support, effectively reducing the active topology to a standard bipartite ladder without 
any edge-crossing penalty.
\end{proof}

Applying Lemma~\ref{LEM_Tip_Depopulation} at the beginning of the state-to-state routing protocol
depopulates the centroid level, reducing the active topology to a perfect bipartite ladder on each
parity subspace. The routing can therefore be designed uniformly for both ladders, the only overhead
being the bounded, dimension-independent cost $\frac{1}{\sqrt{2}} C_{vert}$ (or 
$\sqrt{2} C_{vert}$ for the two-stage gather-scatter protocol).

\begin{remark}
The parity of the truncation degree $d$  dictates the relative dimensions of the reduced bipartite ladders. 
For an even $d$, the centroid depopulation symmetrizes the network, yielding two identical ladders of 
$d$ vertices. Conversely, for an odd $d$, this reduction locks the active state spaces 
into two ladders of distinct sizes, specifically $d+1$ and $d-1$, Figure \ref{fig:ladder_topology_d9}.
\end{remark}

\subsection{State-to-State transfer}\label{SEC_Chiral_StS}

\begin{proposition}[Simultaneous State-to-State Controllability within Parity Sectors]
\label{PRO_Simultaneous_State_Transfer}
Let $D = 2d+1$. The truncated state space $\mathbb{C}^D$ algebraically decomposes into two orthogonal, 
dynamically invariant parity subspaces: $\mathcal{V}_{even} = \mathrm{span}\{\phi_m^B \mid m 
\text{ is even}\}$ and $\mathcal{V}_{odd} = \mathrm{span}\{\phi_m^B \mid m \text{ is odd}\}$. 

For any initial state $\Psi_I \in \mathbb{C}^D$ and target state $\Psi_F \in \mathbb{C}^D$ satisfying the 
parity mass conservation:
\begin{equation}
    \|\Pi_{even} \Psi_I\|^2 = \|\Pi_{even} \Psi_F\|^2 \quad \text{and} \quad \|\Pi_{odd} \Psi_I\|^2 = 
    \|\Pi_{odd} \Psi_F\|^2,
\end{equation}
where $\Pi_{even, odd}$ are the respective spectral projectors, there exists a control protocol achieving 
the simultaneous  transfer $\Psi_I \to \Psi_F$ in a time  bounded by: $     \tau(\Psi_I ,\Psi_F) \le 
\frac{160}{d}$.
\end{proposition}
\begin{proof}

\textbf{Preparation.}
To avoid the divergent edge-crossing penalty at the spectral boundaries, the topological centroids for both $\mathcal{V}_{even}$ and $\mathcal{V}_{odd}$ are fixed at their respective innermost modes, which reside at an index distance of at most $\Delta m_c \le 2$ from the spectral center $(D+1)/2 = d+1$.

Let $\alpha^* = \max_{e} \alpha_e \leq 2$ be the maximum harmonic resonance penalty evaluated on these centered routing trees. 

\emph{1. Vertical Sweeps (Intra-branch bound):} 
From Lemma \ref{LEM_time_constants_C_cross_C_vert}, the temporal cost of any individual vertical step is 
bounded by $C_{vert} \le \frac{2d+2}{2d+1}\frac{3\pi^3 \alpha^*}{8(d+1)^2}$.

\emph{2. Chiral Crossing (Inter-branch bound):}
The crossing is executed at the centroid mode $m_c$, where $m_c = d+1 - k$ with $k \in \{1, 2\}$ depending 
on the parity sector. From Lemma \ref{LEM_time_constants_C_cross_C_vert}, the maximum crossing time at the 
center is $C_{cross} \le \frac{2d+2}{2d+1}\frac{16\pi \alpha^*}{(d+1)^2}$.

\emph{3. Total Temporal Bound:}

The preceding analysis does not allow us to rule out rational dependencies in the spectral gaps of the operator 
$B$. To ensure phase alignment, we must therefore proceed as in Proposition 
\ref{PRO_Result_Constructive_State_to_State_detailed} via a gather-scatter approach. 
As established in Remark \ref{REM_Subspace_decoupling}, the interaction-picture generator is 
block-diagonal, so the even and odd sectors evolve independently under the shared control parameter $v(t)$. 
The routing protocol is designed on the pseudo-spanning tree $\tilde{T}$
(Remark \ref{REM_Subspace_decoupling}) with a zero-flux constraint on the virtual edge, ensuring that no 
population crosses between sectors. Each Givens rotation synthesized on an edge of one sector leaves the 
other sector in a state modified only by a diagonal phase, which is tracked and compensated at the end of 
the protocol. The two sectors can therefore be routed sequentially, and the total time is the sum of the 
individual routing times. 

After elimination of the centroid through Lemma \ref{LEM_Tip_Depopulation}, the routing topology for each parity subspace forms an independent, parity-conserving ladder graph with two legs of size $d/2$ for a total number of vertices of $d$ if $d$ is even, or two legs of size $(d-1)/2$ or $(d+1)/2$ for a total number of vertices of $d-1$ or $d+1$ if $d$ is odd.

\textbf{Routing strategy.}
For each parity subsystem (even or odd), we distinguish two steps.

We first consider the case where $N=2d+2$ is not a multiple of 6, so that no parasitic resonance occurs (Lemma~\ref{LEM_Global_Non_Resonance}); the resonant case is treated in the next paragraph.
 
The gathering step targets the top rung on the left leg of the ladder, denoted $c_L$. It occurs in two stages. First, using  at most $d/2$  vertical moves, we collect the population from the right leg onto the chiral twin $c_R$ of $c_L$. A single transverse move along the central rung then transfers the entire population to the left leg of the ladder. Finally, at most $d/2$ vertical moves suffice to gather the entire population onto $c_L$.

\textbf{Avoidance of the parasitic resonances.}
When $N \geq 18$ is a multiple of $6$, the parasitic resonances listed in Lemma~\ref{LEM_Global_Non_Resonance} occur. They form two classes, which the routing protocol neutralizes by two distinct mechanisms, temporal exclusion and spatial isolation.

\textit{1. Temporal exclusion of the recurrent resonances ($N=6q$, $q \ge 3$):}
The co-resonant quadruplet is $\{(3q\!-\!2,3q),\ (3q,3q\!+\!2),\ (q\!-\!2,q\!+\!2),\ (5q\!-\!2,5q\!+\!2)\}$: the two vertical edges $(3q-2,3q)$ and $(3q,3q+2)$ are incident to the central mode $3q=N/2$, and the two companion transitions $(q-2,q+2)$ and $(5q-2,5q+2)$ have index gap $4$. For even $q$ the quadruplet lies in $\mathcal{V}_{even}$, for odd $q$ in $\mathcal{V}_{odd}$; in either case the central mode $3q$ is the centroid of that sector and is emptied by the preliminary tip depopulation (Lemma~\ref{LEM_Tip_Depopulation}). The two vertical edges incident to $3q$ are therefore never driven during the subsequent gather-scatter phases, and the co-resonant companion transitions are never activated. The argument is uniform in the parity of $q$.

\textit{2. Spatial isolation of the exceptional resonance ($N=30$):}
The exceptional quadruplet $\{(8,10),\ (20,22),\ (2,6),\ (24,28)\}$ lies entirely within the bipartite sector and is not incident to the central mode $15$ (which belongs to the opposite, tip-carrying sector), so temporal exclusion does not apply. The resonance is neutralized by the ordering of the gather: the gathering phase sweeps the population inward from the spectral boundaries, depopulating each rung exactly once and leaving it empty thereafter. As the sweep visits the peripheral rungs before the interior ones, the companion rungs $2,6$ (resp.\ $24,28$) are already empty by the time the vertical edge $(8,10)$ (resp.\ $(20,22)$) is driven, so the parasitic generator acts on the kernel of the local density matrix ($R_{parasite}|0\rangle=|0\rangle$) and no parasitic transfer occurs. The symmetric parasitic rotations induced on the opposite chiral leg are unitary and are undone by the time-symmetric scattering phase.

\textbf{Total time estimates.}
The scattering step is analogous. Each of the two steps incurs a cost of 
${d}C_{vert} + C_{cross}$, resulting in a total cost of $2 d C_{vert} + 2 C_{cross}$. 
When including the triangle-tip overhead (Lemma \ref{LEM_Tip_Depopulation}), the total cost over both subspaces is therefore 
\begin{equation}\label{EQ_majoration_temps_SS_rotateur}
    \tau \le (4d+2) \, C_{vert} + 4C_{cross} \le \frac{2d+2}{2d+1} \left[ \frac{3 \pi^3 \alpha^{\ast}}{2 (d+1)} + \frac{64 \pi \alpha^\ast}{(d+1)^2} \right] \leq \frac{160}{d}.
\end{equation}
This upper bound gives the asymptotic scaling $\tau \le \mathcal{O}(1/d)$.
\end{proof}

\begin{proof}[Proof of Proposition \ref{PRO_rotateur_introduction}] The matrices $A_o$ and 
$B_o$ appearing in the statement of Proposition \ref{PRO_rotateur_introduction} are the 
matrices $A_{odd}$ and $B_{odd}$ corresponding to the restriction of the system $(A,B)$ to 
the subspace of odd functions.

 For the upper bound, the result follows from the proof of Proposition 
\ref{PRO_Simultaneous_State_Transfer} (Step 3 for the odd subspace).

For the lower bound, we use Proposition \ref{PRO_minoration_Ehrenfest} with observable $O=B_o^2$ and the 
specific pair where $\psi_0$ and $\psi_f$ are respectively  the first and the second eigenstates of $A_{odd}$ 
(hence both the source and the target belong to the odd 
subsystem). The commutator $C = [B_o^2, A_o]$ is evaluated element-wise. Since $A_o = -i \, 
\text{diag}(k^2)$, the matrix elements of the commutator are given by:
\begin{equation}
C_{j,k} = (B^2 A - A B^2)_{j,k} = -i(k^2 - j^2) (B_o^2)_{j,k}.
\end{equation}

The diagonal elements of $B_o^2$ (where $k=j$) cancel out in the commutator. The non-zero entries of $C$ lie on the second super-diagonal and sub-diagonal (where $|k-j|=2$):
\begin{itemize}
\item For $k = j+2$: $C_{j, j+2} = -i((j+2)^2 - j^2)\left(-\frac{1}{4}\right) = i(j+1)$.
\item For $k = j-2$: $C_{j, j-2} = -i((j-2)^2 - j^2)\left(-\frac{1}{4}\right) = -i(j-1)$.
\end{itemize}

The matrix representation of the commutator is therefore:
\begin{equation}
[B^2, A] = i \begin{pmatrix}
0 & 0 & 2 & 0 & 0 & \dots \\
0 & 0 & 0 & 3 & 0 & \dots \\
-2 & 0 & 0 & 0 & 4 & \dots \\
0 & -3 & 0 & 0 & 0 & \dots \\
\vdots & \vdots & \vdots & \vdots & \vdots & \ddots
\end{pmatrix}.
\end{equation}

This matrix is Hermitian and block-tridiagonal (decoupling even and odd indices). To bound its spectral operator norm $\| C \|_{op}$, we apply the Gershgorin Circle theorem.

For a matrix $C \in \mathbb{C}^{d \times d}$, the Gershgorin Circle theorem states that every eigenvalue $\lambda$ of $C$ lies within at least one Gershgorin disc $D(C_{j,j}, R_j)$ centered at $C_{j,j}$ with radius $R_j$, such that:
\begin{equation}
|\lambda - C_{j,j}| \le R_j = \sum_{k \neq j} |C_{j,k}|.
\end{equation}
Since $C$ is Hermitian, it is normal, and its operator norm equals its spectral radius ($\| C \|_{op} = \max |\lambda|$). Because all diagonal entries are identically zero ($C_{j,j} = 0$), the absolute value of any eigenvalue is bounded by the maximum absolute row sum:
\begin{equation}
\| C \|_{op} \le \max_{1 \le j \le d} \sum_{k \neq j} |C_{j,k}|.
\end{equation}

Using the particular block-tridiagonal structure of $C=[B_o^2, A_o]$, one gets
 $R_1 = |C_{1,3}| =  2$, $R_2 = |C_{2,4}| =  3$, $R_{d-1} =d-2$, $R_d=d-1$ and, for $3 \le j \le d-2$:
\[ R_j = |C_{j, j-2}| + |C_{j, j+2}| = (j-1) + (j+1) = 2j. \]

The operator norm is bounded by the maximum over the set of $R_j$. For $d \ge 5$, the bulk values $R_j = 2j$ ($3 \le j \le d-2$) are strictly increasing, while the boundary values $R_1 = 2$, $R_2 = 3$, $R_{d-1} = d-2$ and $R_d = d-1$ are all strictly smaller than $2(d-2)$; the global maximum is therefore $R_{d-2} = 2d-4$, attained at the last bulk row $j = d-2$. For $d = 4$ the bulk is empty and the maximum is attained at the boundary rows $j \in \{2, d\}$, with value $3 \le 2d-4$.

This provides an algebraic upper bound for any dimension $d \ge 4$:
$$ \| [B_o^2, A_o] \|_{op} \le \max_{1 \le j \le d} R_j \le 2d - 4 \qquad (d \ge 4). $$

Since $2d - 4 < 2d$ for all $d \ge 2$, the simplified universal upper bound holds:
$$ \| [B_o^2, A_o] \|_{op} < 2d. $$

Substituting this bound back into Proposition~\ref{PRO_minoration_Ehrenfest} (where the numerator evaluates to $1/4$), one gets the claimed lower bound:
$\tau_{12}(d) \ge \frac{1/4}{2d - 4} > \frac{1}{8d}. $
\end{proof}

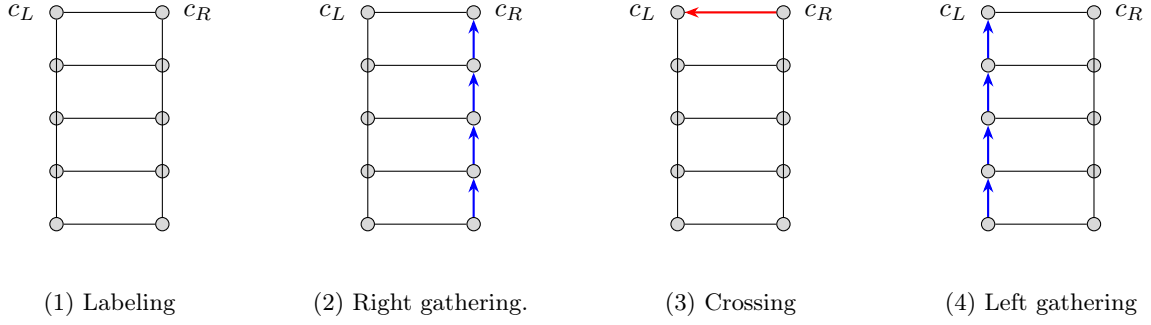
\begin{figure}[htbp]
    \centering
    \tikzstyle{vertex}=[circle, draw, fill=black!15, inner sep=0pt, minimum size=5pt]
    \tikzstyle{activevertex}=[circle, draw, fill=black!60, inner sep=0pt, minimum size=5pt]
    \tikzstyle{edge}=[-{Stealth[scale=0.8]}, thick]
    \tikzstyle{undir}=[-]

    \begin{subfigure}[b]{0.23\textwidth}
        \centering
        \begin{tikzpicture}[scale=0.7]
            \foreach \y in {1,2,3,4,5} {
                \node[vertex] (L\y) at (0,\y) {};
                \node[vertex] (R\y) at (2,\y) {};
            }
            \foreach \y in {1,2,3,4,5} {
                \draw[undir] (L\y) -- (R\y);
            }
            \draw[undir] (L1) -- (L5);
            \draw[undir] (R1) -- (R5);

            \node[left=0.05cm of L5] {\small $c_L$};
            \node[right=0.05cm of R5] {\small $c_R$};

            \node at (1, -0.5) {\footnotesize (1) Labeling};
        \end{tikzpicture}
    \end{subfigure}
    \hfill
    \begin{subfigure}[b]{0.23\textwidth}
        \centering
        \begin{tikzpicture}[scale=0.7]
            \foreach \y in {1,2,3,4,5} {
                \node[vertex] (L\y) at (0,\y) {};
                \node[vertex] (R\y) at (2,\y) {};
            }
            \foreach \y in {1,2,3,4,5} {
                \draw[undir] (L\y) -- (R\y);
            }
            \draw[undir] (L1) -- (L5);
            
            \draw[edge, blue, thick] (R1) -- (R2);
            \draw[edge, blue, thick] (R2) -- (R3);
            \draw[edge, blue, thick] (R3) -- (R4);
            \draw[edge, blue, thick] (R4) -- (R5);
			\node[left=0.05cm of L5] {\small $c_L$};
            \node[right=0.05cm of R5] {\small $c_R$};

            \node at (1, -0.5) {\footnotesize (2) Right gathering.};
        \end{tikzpicture}
    \end{subfigure}
    \hfill
    \begin{subfigure}[b]{0.23\textwidth}
        \centering
        \begin{tikzpicture}[scale=0.7]
            \foreach \y in {1,2,3,4,5} {
                \node[vertex] (L\y) at (0,\y) {};
                \node[vertex] (R\y) at (2,\y) {};
            }
            \foreach \y in {1,2,3,4,5} {
                \draw[undir] (L\y) -- (R\y);
            }
            \draw[undir] (L1) -- (L5);
            \draw[undir] (R1) -- (R5);

            \draw[edge, red, thick] (R5) -- (L5);

            \node[left=0.05cm of L5] {\small $c_L$};
            \node[right=0.05cm of R5] {\small $c_R$};

            \node at (1, -0.5) {\footnotesize (3) Crossing};
        \end{tikzpicture}
    \end{subfigure}
    \hfill
    \begin{subfigure}[b]{0.23\textwidth}
        \centering
        \begin{tikzpicture}[scale=0.7]
            \foreach \y in {1,2,3,4,5} {
                \node[vertex] (L\y) at (0,\y) {};
                \node[vertex] (R\y) at (2,\y) {};
            }
            \foreach \y in {1,2,3,4,5} {
                \draw[undir] (L\y) -- (R\y);
            }
            \draw[undir] (R1) -- (R5);

            \draw[edge, blue, thick] (L1) -- (L2);
            \draw[edge, blue, thick] (L2) -- (L3);
            \draw[edge, blue, thick] (L3) -- (L4);
            \draw[edge, blue, thick] (L4) -- (L5);

        		\node[left=0.05cm of L5] {\small $c_L$};
            \node[right=0.05cm of R5] {\small $c_R$};

            \node at (1, -0.5) {\footnotesize (4) Left gathering};
        \end{tikzpicture}
    \end{subfigure}

\caption{Routing protocol for the proof of Proposition \ref{PRO_Simultaneous_State_Transfer}. The sequence illustrates the concentration of 
probability mass on a central node following the depopulation of the self-chiral center 
in a triangle-tipped ladder. (1) The ladder architecture identifies centroids $c_L$ and 
$c_R$. Although all horizontal rungs are depicted to illustrate the chiral symmetry of 
the spectrum, the active protocol exclusively uses the topologically central rung 
(visually positioned at the top of the diagram), as the lower transitions are 
prohibitively slow. (2) Right-side gathering concentrates mass at $c_R$ while mirroring 
operations on the left leg. (3) The total mass is transferred to the left leg via the 
central rung. (4) Mass is concentrated at $c_L$ before a symmetric scattering phase 
routes it to the target.}
    \label{fig:routing_sequence}
\end{figure}

\subsection{Unitary Propagator Synthesis and the Chiral Penalty}
\label{SEC_Chiral_unitary_synthesis}

Due to spectral chiral symmetry, the gathering-pruning strategy of Section~\ref{SEC_Contructive_bounds_Propagators} cannot be applied directly. We start with a technical lemma on simultaneous control at the center of the ladder; its ideas mirror the dynamical decoupling of Section III.A in \cite{PhysRevA.58.2733}.
 
\begin{lemma}\label{LEM_center_ladder_simultaneous_decoupling}
Consider the bilinear control system on the special unitary group $SU(4)$ described by the differential equation $\dot{X}(t) = \Omega(u(t))X(t)$, where $X(t) \in SU(4)$ with $X(0)=I_4$ and the parameter-dependent matrix $\Omega(u)$ is given by 

$$\Omega(u)=\begin{pmatrix}
0 & u_1(t)  & u_2(t) & 0 \\
-\bar{u}_1(t) & 0 & 0 & u_3(t) \\
 -\bar{u}_2(t) & 0 & 0& u_1(t) \\
0 & -\bar{u}_3(t) &-\bar{u}_1(t) & 0 
\end{pmatrix}$$

with the control functions 
$u_k \in L^1_{\text{loc}}(\mathbb{R}_+, \mathbb{C})$ for $k \in \{1, \dots, 3\}$. 
Let the control cost be defined as 

$$J(u) = \frac{2}{\pi} C_{vert} \|u_1\|_{L^1} + \frac{2}{\pi} C_{cross} \|u_2\|_{L^1}
+ \frac{2}{\pi} C_{cross}' \|u_3\|_{L^1} $$

where $C_{vert}>0$, $C_{cross}>0$ and $C_{cross}'>0$ are fixed weighting coefficients . Let $e_k$ denote the $k$-th canonical basis vector of $\mathbb{C}^4$ .

For every $r \in [-\pi,\pi]$, every initial phase $\tilde{\theta} \in \mathbb{R}$, and every target 
relative phase $\xi \in \mathbb{R}$, there exists an explicit, piecewise constant control sequence driving 
the system from the initial state $x(0)$ to a final state $x(T)$ such that the simultaneous conditions 
$X(T)e_1 = e^{i\phi_1}e_1$  and $X(T) (0,0,\cos(r),e^{i \tilde{\theta}}\sin(r) )^T = e^{i\phi_3}e_3$ are 
satisfied with the phase synchronization $\phi_1 - \phi_3 \equiv \xi \pmod{2\pi}$. Furthermore, the 
total functional cost of this strategy satisfies the upper bound:

$$J(u) \le 2 C_{vert} + 4 C_{cross}$$

\end{lemma}

\begin{proof}
The unperturbed control field $u_1$ activates transitions on the $(1,2)$ and $(3,4)$ subspaces simultaneously . To decouple these dynamics and inject the target relative phase $\xi$, we use phase conjugation via the transversal control $u_2$ on the $(1,3)$ subspace. The remaining coupling $u_3$ (the $(2,4)$ transition) plays no role in this strategy: it is a deliberate choice, not an accident, that an efficient (though non-optimal) decoupling uses only $u_1$ and $u_2$. We therefore set $u_3\equiv0$ throughout, so the $C_{cross}'\|u_3\|_{L^1}$ term in the cost vanishes.

An elementary $\pi$-pulse on the $(1,3)$ transition with constant complex phase $\alpha$, corresponding to a field $u_2(t)$ such that $\int |u_2| dt = \pi/2$ and $\arg(u_2) = \alpha$, generates the unitary operator $V(\alpha) = e^{i\alpha}E_{13} - e^{-i\alpha}E_{31} + E_{22} + E_{44}$. Composing two such pulses yields a diagonal phase gate $D(\Delta) = V(\alpha+\Delta)V(\alpha)$:
\begin{equation}
    D(\Delta) = -e^{i\Delta}E_{11} + E_{22} - e^{-i\Delta}E_{33} + E_{44} = \mathrm{diag}(-e^{i\Delta}, 1, -e^{-i\Delta}, 1).
\end{equation}
The $L^1$-norm of $u_2$ required to synthesize $D(\Delta)$ is $\pi$, independent of the chosen shift $\Delta$.

We set $\Delta_1 = -\pi/2$ to construct the intermediate decoupling gate $D_{13} = D(-\pi/2) = \mathrm{diag}(i, 1, -i, 1)$. We design an active sequence of two distinct $u_1$ pulses interleaved with $u_2$ gates. Let $u_1^{(1)} = \frac{r}{2} e^{-i\tilde{\theta}}$. The second pulse is phase-shifted by $-\pi/2$, setting $u_1^{(2)} = -i u_1^{(1)}$. Following the second pulse, we apply a return gate $D_{ret} = D(\Delta_2)$. 
    
The total evolution operator is $\mathcal{M} = D_{ret} \exp(\Omega(u_1^{(2)})) D_{13} \exp(\Omega(u_1^{(1)}))$. By inserting the trivial identity $D_{13}^{-1} D_{13} = I$, this factors as:
\begin{equation}
    \mathcal{M} = (D_{ret} D_{13}) \left[ D_{13}^{-1} \exp(\Omega(-i u_1^{(1)})) D_{13} \right] \exp(\Omega(u_1^{(1)})).
\end{equation}
Direct evaluation of the conjugation by the diagonal matrix $D_{13}$ yields $D_{13}^{-1} \Omega(-i u_1^{(1)}) D_{13} = -\Omega(u_1^{(1)})|_{(1,2)} + \Omega(u_1^{(1)})|_{(3,4)}$. 
    
Consequently, the active $u_1$ propagators cancel on the $(1,2)$ subspace, while they constructively add on the $(3,4)$ subspace to form $\exp(2\Omega(u_1^{(1)}))$. The state is subsequently multiplied by the boundary gates $D_{ret} D_{13} = \mathrm{diag}(-i e^{i\Delta_2}, 1, i e^{-i\Delta_2}, 1)$.

For the first initial condition $x_1(0) = e_1$, the population is not stationary on $e_1$ during the sequence (it leaks into $e_2$ while the $(1,2)$ block is driven); however, the two $u_1$ pulses cancel on the $(1,2)$ subspace, so the state returns to $e_1$ and accumulates only the boundary-gate phase $\phi_1 = \Delta_2 - \pi/2$. Thus $x_1(T) = e^{i\phi_1} e_1$.
    
For the second initial condition $x_2(0) = \cos(r)e_3 + e^{i\tilde{\theta}}\sin(r)e_4$, the action of $\exp(2\Omega(u_1^{(1)}))$ generates a complete phased Givens rotation of angle $r$, mapping the state to $e_3$. The terminal gates apply the phase factor $i e^{-i\Delta_2}$, setting $\phi_3 = -\Delta_2 + \pi/2$. Thus $x_2(T) = e^{i\phi_3} e_3$.

The resulting relative phase is $\phi_1 - \phi_3 = 2\Delta_2 - \pi$. Setting the phase shift of the return gate to $\Delta_2 = (\xi + \pi)/2$  satisfies the synchronization constraint $\phi_1 - \phi_3 = \xi$.

The functional cost sums the $L^1$-norms of the used fields. The transversal control synthesizes two diagonal gates ($D_{13}$ and $D_{ret}$), with $\|u_2\|_{L^1} = \pi + \pi = 2\pi$. The continuous control $u_1$ accumulates $\|u_1\|_{L^1} = |u_1^{(1)}| + |u_1^{(2)}| = 2 \times |\frac{r}{2}| = |r|$. Using $\max |r| = \pi$, the global functional cost is bounded by:
\begin{equation}
    J(u) = \frac{2}{\pi} C_{vert} |r| + \frac{2}{\pi} C_{cross} (2\pi) \le 2 C_{vert} + 4 C_{cross}.
\end{equation}
\end{proof}

\begin{remark}
 This explicit trajectory solves the simultaneous state transfer without cross-talk, but is 
 not optimal under the $L^1$ (sub-Finslerian) metric. The loss of optimality is the fixed 
 cost $4 C_{cross}$ of the echo cycle (the phase gate $D_{13}$ and its inverse via the 
 transversal control $u_2$), which forces a zero-sum evolution on the $(1,2)$ subspace.
\end{remark}

\begin{proposition}[Global Unitary Synthesis Time on Parity Subspaces]
\label{PRO_Unitary_Synthesis_Constant_Time}

Let $d$ be an integer such that $d+1$ is not divisible by $3$.
Let $\mathcal{V}_p$ be one of the dynamically invariant parity subspaces (either $\mathcal{V}_{even}$ or $\mathcal{V}_{odd}$) of dimension $N_p \le d+1$. For any target unitary operator $U \in SU(\mathcal{V}_p)$, there exists a constructive control sequence that synthesizes $U$. 

 The total synthesis time $\tau_{\mathrm{US}}$ is bounded by a 
dimension-independent  constant, satisfying:
$$ \tau^{\max}_{\mathrm{US}} \le 345 $$
\end{proposition}

\begin{proof} The proof is done for $d\geq 5$ (which allows to use the estimates 
of Lemma \ref{LEM_time_constants_C_cross_C_vert}). The case $d = 3$ can be checked directly. 

The assumption that $d+1$ is not divisible by 3 implies that $N=D+1=2d+2$ is not divisible by 6. Lemma \ref{LEM_Global_Non_Resonance} ensures that the spectral ladders (both on the even and on the odd space) contain no parasitic resonances. 

\emph{Initialization.} If the active routing sub-graph contains an odd number of vertices, we denote its unique central node by $m^\ast$ and first route $U^{-1}(\phi_{m^\ast}^B)$ to $\phi_{m^\ast}^B$. The node is locked and permanently removed from the active routing graph. This initialization step incurs a temporal cost bounded by $\tau_{init} \le (d+1) C_{vert} + 2 C_{cross}$.

\emph{General procedure.}
The synthesis operates pair-by-pair on the chiral ladder graph, processing $N_{pairs} \le \frac{d+1}{2}$ target chiral pairs. For the pair $(m_j, \bar{m}_j)$ at step $j$, (with $2j\leq d$) the routing sequence comprises four sequential phases:

\emph{1. Asymmetric Parking:} 
Let $\Psi_k = U^{-1}(\phi_{m_j}^B)$. $\Psi_k$ is gathered entirely to $c_L$. Since vertical couplings satisfy $C_{vert} \le \frac{2d+2}{2d+1}\frac{3\pi^3 \alpha^*}{8(d+1)^2}$ and the central rung satisfies $C_{cross} \le \frac{2d+2}{2d+1}\frac{16\pi \alpha^*}{(d+1)^2}$, the gathering time is bounded by $\tau_A \le (d-2j+2) C_{vert} + C_{cross}$.

\emph{2. Orthogonal Gather:} 
Let $\Psi_{\bar{k}} = U^{-1}(\phi_{\bar{m}_j}^B)$. By unitarity, $\langle \Psi_k | \Psi_{\bar{k}} \rangle = 
0$. Because $\Psi_k$ occupies $c_L$ with probability 1, $\Psi_{\bar{k}}$ has zero mass on $c_L$. 
The mass of $\Psi_{\bar{k}}$ is gathered to $c'_R$ (the right node immediately below the centroid). This 
requires gathering the left leg to $c'_L$ and the right leg to $c'_R$, followed by a single transversal 
crossing at the adjacent rung to transfer the remaining left-leg mass. The total cost is bounded by $\tau_B 
\le (d-2j) C_{vert} + C'_{cross}$.

\emph{3. Centroid Alignment and Synchronization:} 
Applying Lemma \ref{LEM_center_ladder_simultaneous_decoupling} with $c_L \mapsto e_1$, $c'_L \mapsto e_2$, $c_R \mapsto e_3$, and $c'_R \mapsto e_4$, both states are simultaneously mapped to the central chiral rung $(c_L, c_R)$ while their relative phase is locked via the transversal control mechanism. The cost is bounded by $\tau_C \le 2 C_{vert} + 4 C_{cross}$.

\emph{4. Parallel Scatter:} 
$\Psi_k$ is routed down the left upright to $m_j$, while $\Psi_{\bar{k}}$ is routed down the right upright to $\bar{m}_j$. 
By the chiral symmetry of the spectrum (Corollary \ref{COR_Chiral_Coupling}), the two transitions $(m,m+2)$ and $(\bar m,\bar m+2)$ have identical gap, identical module, and identical harmonic penalty. A single element $M_{part}(\phi) \in CC(t)$ isolates both chiral transitions simultaneously. 
The full population transfer to the adjacent level requires a rotation angle of $\pi/2$, generating the tensor propagator $R(\pi/2, \phi) \otimes R(\pi/2, \phi)$. Applying this operator shifts the argument of both states by $-\phi + \pi$. This preserves the relative phase locked in Step 3, while the free control parameter $\phi$ provides the exact degree of freedom required to adjust the global absolute phase of the pair to its target value. The descent is executed in one single pass on both uprights, at a time cost of $C_{vert}$ per level, yielding a scattering time bounded by $\tau_E \le \left( \frac{d}{2} - j + 1 \right) C_{vert}$.

\emph{Integration of the Temporal Bound:} 
Summing the individual phases, the total time to process a single pair evaluates to:
$$ \begin{aligned}
    \tau_{pair} &= \tau_A + \tau_B + \tau_C + \tau_E \\
    &\le \left[ (d-2j+2) + (d-2j) + 2 + \frac{d}{2} - j + 1 \right] C_{vert} + (1 + 4) C_{cross} + C'_{cross} \\
    &= \left( \frac{5d}{2} - 5j + 5 \right) C_{vert} + 5 C_{cross} + C'_{cross}.
\end{aligned} $$

Summing over the active graph from $j=1$ to $N_{pairs} \le \frac{d+1}{2}$, and adding $\tau_{init}$, provides the  upper bound:
$$ \begin{aligned}
 \tau^{\max}_{\mathrm{US}} &\le (d+1) C_{vert} + 2 C_{cross} + \sum_{j=1}^{\lceil d/2 \rceil} \left[ \left( \frac{5d}{2} - 5j + 5 \right) C_{vert} + 5 C_{cross} + C'_{cross} \right] \\
 &\le \frac{(d+1)(5d+13)}{8} C_{vert} + \frac{5d+9}{2} C_{cross} + \frac{d+1}{2} C'_{cross}.
\end{aligned} $$

Substituting the analytical bounds $C_{vert} \le \frac{2d+2}{2d+1}\frac{3\pi^3 \alpha^*}{8(d+1)^2}$, $C_{cross} \le \frac{2d+2}{2d+1}\frac{16\pi \alpha^*}{(d+1)^2}$, and $C'_{cross} \le \frac{2d+2}{2d+1}\frac{36\pi \alpha^*}{(d+1)^2}$ yields the algebraic rational form:
$$ \tau^{\max}_{\mathrm{US}} \le \alpha^* \left( \frac{2d+2}{2d+1} \right) \left[ \frac{3\pi^3(5d+13)}{64(d+1)} + \frac{8\pi(5d+9)}{(d+1)^2} + \frac{18\pi}{d+1} \right]. $$

Each term in the bracket is monotonically decreasing with respect to the truncation degree $d \ge 1$, as is the prefactor $\frac{2d+2}{2d+1}$; the right-hand side is therefore maximized at $d=1$ over all $d\ge1$, so its value there is a uniform a priori upper bound for every admissible dimension (the proof being established for $d\ge5$, with $d=3$ checked directly). Evaluating the worst-case geometric configuration at $d=1$ yields a prefactor of $4/3$:
$$ \tau^{\max}_{\mathrm{US}} \le \alpha^* \frac{4}{3} \left( \frac{27\pi^3}{64} + 28\pi + 9\pi \right) = \alpha^* \frac{4}{3} \left( \frac{27\pi^3}{64} + 37\pi \right) \le 2 \times \frac{4}{3} \times 129.32 \le 345. $$
\end{proof}

\begin{proposition}[Global Synthesis of Parity-Preserving Unitaries]
\label{PRO_Full_Parity_Preserving_Synthesis}
Let $D = 2d+1$. The full state space decomposes into orthogonal parity subspaces $\mathbb{C}^D = \mathcal{V}_{even} \oplus \mathcal{V}_{odd}$. Let $SU_{p}(D)$ denote the subgroup of $SU(D)$ consisting of all parity-preserving unitary transformations, such that any $U \in SU_{p}(D)$ takes the block-diagonal form $U = U_{even} \oplus U_{odd}$, with $U_{even} \in SU(\mathcal{V}_{even})$ and $U_{odd} \in SU(\mathcal{V}_{odd})$.

Assume that $D+1$ is not divisible by 6. Using a single scalar control field $u(t)$, the synthesis of any target $U \in SU_{p}(D)$ can be achieved in a bounded macroscopic time $\tau_{total} \leq 690$.
\end{proposition}

\begin{remark}[On the magnitude of the numerical constants]
The numerical constants appearing in the constructive upper bounds (e.g., $160$, $345$, $690$) are conservative. They result from the systematic accumulation of worst-case inequalities required to ensure a dimension-independent validity: uniform bounding of the effective couplings via Gershgorin circles, maximal path lengths on the bipartite ladders, and the maximum theoretical Diophantine penalty ($\alpha^* \le 2$). 

These values must not be interpreted as tight physical Quantum Speed Limits. Their purpose is to establish the asymptotic scaling behaviors of the system --specifically $\mathcal{O}(1/d)$ for state-to-state transfer and $\mathcal{O}(1)$ for global unitary synthesis---without relying on the numerical evaluation of accidental spectral resonances for each specific dimension $d$.
\end{remark}

\begin{proof}
Assumption 1 is satisfied independently on the invariant subspaces $\mathcal{V}_{even}$ 
and $\mathcal{V}_{odd}$. Consequently, Lemma \ref{LEM_approx_cont} applies to the 
restricted subgroups $SU(\mathcal{V}_{even})$ and $SU(\mathcal{V}_{odd})$, ensuring the 
equivalence between exact and approximate synthesis times locally within each block. 
Exact accessibility on the parity-preserving group $SU_p(D)$ is then achieved via the 
sequential concatenation of the exact trajectories from each subspace, since the 
dynamical cross-talk induced on the passive block is limited to a deterministically 
compensable diagonal phase drift.

The target operator is factored as $U = (U_{even} \oplus \mathbb{I}_{odd}) (\mathbb{I}_{even} \oplus U_{odd})$. 
The synthesis protocol is split into two consecutive phases, requiring pre-compensation 
for the deterministic unperturbed phase drifts $A_{D, even}$ and $A_{D, odd}$ occurring 
during the respective inactive windows:
\begin{enumerate}
    \item \textbf{Odd Subspace Synthesis:} The algorithm described in Proposition \ref{PRO_Unitary_Synthesis_Constant_Time} is applied to the odd ladder graph. The 
   required synthesis time $\tau_{even}$ for the subsequent Phase 2 is 
   deterministically known. To offset the free evolution that will occur during Phase 
   2, the dynamically adjusted target for this first phase is set to $U_{odd}' = \exp(-
   \tau_{even} A_{D, odd}) U_{odd}$. Because the resonant frequencies used for the 
   vertical intra-branch sweeps are distinct from those of the even branch, the even 
   subspace undergoes purely diagonal unperturbed phase drift during this step. The 
   time required is bounded by $\tau_{odd} \le 345$.
    \item \textbf{Even Subspace Synthesis:} The same algorithmic strategy is applied to 
    the even ladder graph. To offset the free evolution that accumulated during Phase 
    1, the dynamically adjusted target is $U_{even}' = U_{even} \exp(-\tau_{odd} A_{D, 
    even}) $. The time required is identical, yielding $\tau_{even} \le 345$.
\end{enumerate}
In addition to the diagonal drift compensation above, the residual diagonal phases 
accumulated in the 
passive sector during each Givens rotation (as described in Remark \ref{REM_Subspace_decoupling}) are 
included in the target-phase adjustment. Since these phases are diagonal in the eigenbasis of $B$, 
they are fully determined by the routing sequence and can be compensated by modifying the target 
unitary before each synthesis step.

Due to the commutativity of the sub-algebras and the anticipation of the diagonal drifts, the 
sequential product yields the  global parity-preserving target operator $U$. The total synthesis time is the
algebraic sum of the phase durations:
\begin{equation}
    \tau_{total} = \tau_{odd} + \tau_{even} \le 690.
\end{equation}
The full block-diagonal controllability over the $D$-dimensional system is achieved in time $\mathcal{O}(1)$, independent of the truncation degree $d$.
\end{proof}

\begin{figure}[htbp]
    \centering
    \tikzstyle{vertex}=[circle, draw, fill=black!10, inner sep=0pt, minimum size=5pt]
    \tikzstyle{cutvertex}=[circle, draw, fill=red!20, dashed, inner sep=0pt, minimum size=5pt]
    \tikzstyle{edge}=[-, thick, black!40]
    \tikzstyle{arrow_blue}=[-{Stealth[scale=0.8]}, blue, thick]
    \tikzstyle{arrow_red}=[-{Stealth[scale=0.8]}, red, thick]
    \tikzstyle{arrow_green}=[-{Stealth[scale=0.8]}, green!60!black, thick]

    \newcommand{\drawLadder}{
        \foreach \y in {1,2,3,4,5} {
            \draw[edge] (0,\y) -- (2,\y);
        }
        \draw[edge] (0,1) -- (0,5);
        \draw[edge] (2,1) -- (2,5);
        
        \foreach \y in {1,2,3,4,5} {
            \node[vertex] at (0,\y) {};
            \node[vertex] at (2,\y) {};
        }
        \node[left=0.05cm of {(0,5)}] {\small $c_L$}; \node[right=0.05cm of {(2,5)}] {\small $c_R$};
        \node[left=0.05cm of {(0,4)}] {\small $c'_L$}; \node[right=0.05cm of {(2,4)}] {\small $c'_R$};
    }

    \begin{subfigure}[b]{0.24\textwidth}\centering
        \begin{tikzpicture}[scale=0.6] \drawLadder
            \draw[arrow_blue] (2,1) -- (2,2);
            \draw[arrow_blue] (2,2) -- (2,3);
            \draw[arrow_blue] (2,3) -- (2,4);
            \draw[arrow_blue] (2,4) -- (2,5);
        \end{tikzpicture}
        \caption{Right Gather 1}
    \end{subfigure}
    \begin{subfigure}[b]{0.24\textwidth}\centering
        \begin{tikzpicture}[scale=0.6] \drawLadder
            \draw[arrow_red] (2,5) -- (0,5);
        \end{tikzpicture}
        \caption{Crossing}
    \end{subfigure}
    \begin{subfigure}[b]{0.24\textwidth}\centering
        \begin{tikzpicture}[scale=0.6] \drawLadder
            \draw[arrow_blue] (0,1) -- (0,2);
            \draw[arrow_blue] (0,2) -- (0,3);
            \draw[arrow_blue] (0,3) -- (0,4);
            \draw[arrow_blue] (0,4) -- (0,5);
        \end{tikzpicture}
        \caption{Left Gather 1}
    \end{subfigure}
    \begin{subfigure}[b]{0.24\textwidth}\centering
        \begin{tikzpicture}[scale=0.6] \drawLadder
            \draw[arrow_blue] (0,3) -- (0,4); 
            \draw[arrow_blue] (0,2) -- (0,3); 
            \draw[arrow_blue] (0,1) -- (0,2);
            \draw[arrow_red] (0,4) -- (2,4);             
        \end{tikzpicture}
        \caption{Left Gather 2}
    \end{subfigure}

    \vspace{0.5cm}
    \begin{subfigure}[b]{0.24\textwidth}\centering
        \begin{tikzpicture}[scale=0.6] \drawLadder
            \draw[arrow_blue] (2,3) -- (2,4);
            \draw[arrow_blue] (2,1) -- (2,2); 
            \draw[arrow_blue] (2,2) -- (2,3);
        \end{tikzpicture}
        \caption{Right Gather 2}
    \end{subfigure}
    \begin{subfigure}[b]{0.24\textwidth}\centering
        \begin{tikzpicture}[scale=0.6] \drawLadder
            \draw[arrow_green, bend left=15] (0,5) to (2,5); 
            \draw[arrow_green, bend left=15] (2,5) to (0,5);
            \draw[arrow_green, bend left=15] (0,4) to (0,5);
            \draw[arrow_green, bend left=15] (2,4) to (2,5);
            \draw[arrow_green, bend left=15] (0,5) to (0,4);
            \draw[arrow_green, bend left=15] (2,5) to (2,4);
        \end{tikzpicture}
        \caption{Stabilization}
    \end{subfigure}
    \begin{subfigure}[b]{0.24\textwidth}\centering
        \begin{tikzpicture}[scale=0.6] \drawLadder
            \draw[arrow_blue] (0,3) -- (0,2);
            \draw[arrow_blue] (0,4) -- (0,3);
            \draw[arrow_blue] (0,5) -- (0,4);
            \draw[arrow_blue] (0,2) -- (0,1);
            \draw[arrow_blue] (2,3) -- (2,2);
            \draw[arrow_blue] (2,4) -- (2,3);
            \draw[arrow_blue] (2,5) -- (2,4);
            \draw[arrow_blue] (2,2) -- (2,1);
        \end{tikzpicture}
        \caption{Downward sweep}
    \end{subfigure}
    \begin{subfigure}[b]{0.24\textwidth}\centering
        \begin{tikzpicture}[scale=0.6]
            \foreach \y in {2,3,4,5} {
                \draw[edge] (0,\y) -- (2,\y);
                \node[vertex] at (0,\y) {}; \node[vertex] at (2,\y) {};
            }
            \draw[edge] (0,2) -- (0,5); \draw[edge] (2,2) -- (2,5);
            \draw[dashed, red] (0,1) -- (2,1);
            \node[cross out, draw=red, inner sep=2pt] at (0,1) {};
            \node[cross out, draw=red, inner sep=2pt] at (2,1) {};
        \end{tikzpicture}
        \caption{Removal}
    \end{subfigure}
\caption{Mechanics of the proof for Proposition \ref{PRO_Unitary_Synthesis_Constant_Time}. First, the preimage $\Psi_k = U^{-1}(\phi_{m_j}^B)$ of the target left leaf is routed to the left centroid $c_L$ using the mechanism from Proposition \ref{PRO_Simultaneous_State_Transfer}: (a) the probability mass on the right leg is concentrated at the right centroid $c_R$; (b) the entire mass is transferred to the left leg via the central rung; (c) the mass is concentrated at $c_L$. Next, the preimage $\Psi_{\bar{k}} = U^{-1}(\phi_{\bar{m}_j}^B)$ of the chiral twin is routed to the adjacent right node $c'_R$. This involves (d) gathering the mass on the left leg at $c'_L$ and transferring it to $c'_R$, followed by (e) gathering the remaining mass on the right leg to $c'_R$. The centroid alignment phase (f) applies the decoupling sequence of Lemma \ref{LEM_center_ladder_simultaneous_decoupling} to position both preimages at $c_L$ and $c_R$. Finally, a simultaneous downward sweep (g) routes both states to their target leaves, after which these locked nodes are removed from the active graph (h).}   
\label{FIG_Simultaneous_Ladder_parking}
\end{figure}
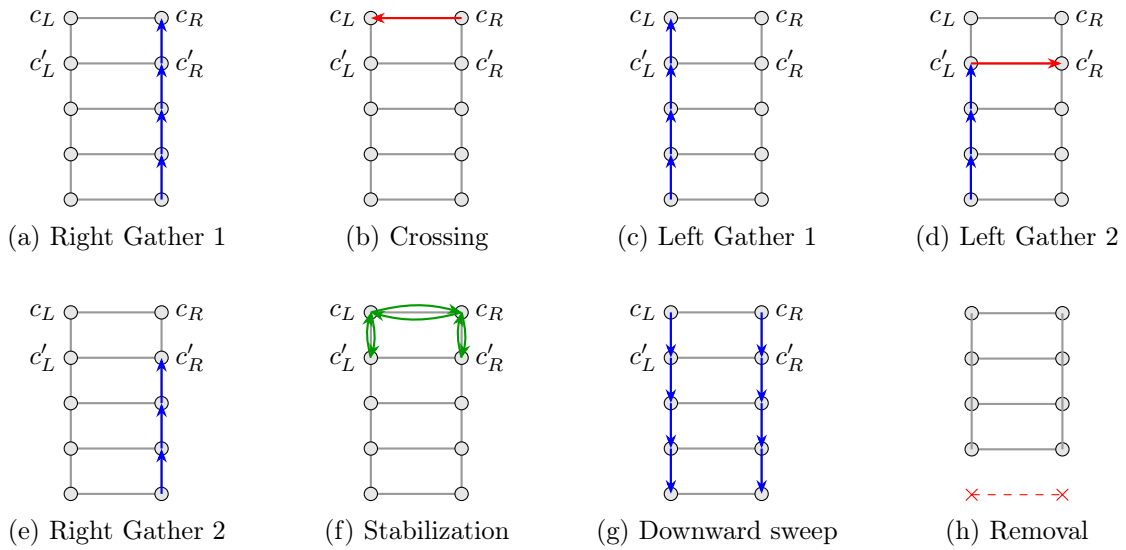

\appendix

\section{Lower bounds}\label{SEC_APP_lower_bounds}

This appendix contains the proofs of Propositions \ref{PRO_Spectral_Bound}, \ref{PRO_Chordal_relaxation} and \ref{PRO_minoration_Ehrenfest}.

\subsection{Spectral bounds}\label{SEC_Spectral_bounds}

\begin{proof}[Proof of Proposition \ref{PRO_Spectral_Bound}]
 Since we operate in the eigenbasis of $B$, the diagonal part of the 
 drift $A_D$ commutes with the control generator, $[A_D, B] = 0$. Let 
 $U_D(t) = \exp\left(t A_D + \int_0^t u(s)B \, \mathrm{d}s\right) \in 
 \mathcal{T}_B$. 

Defining the interaction variable $z(t) = U_D^\dagger(t) \psi(t)$, 
the  dynamics reduce to:
\begin{equation}
    \dot{z}(t) = U_D^\dagger(t) A_{ND} U_D(t) z(t).
\end{equation}

We bound the velocity of $z(t)$ using the unitary invariance of the operator norm. For any $t$, $\|\dot{z}(t)\| \le \|U_D^\dagger(t) A_{ND} U_D(t)\|_{op} \|z(t)\| = \|A_{ND}\|_{op}$, since the state remains on the unit sphere ($\|z(t)\| = 1$).

Integrating the velocity over the interval $[0, \tau]$ gives the length $L$ of the curve $z(t)$ from $z(0) = \psi_0$ to $z(\tau) = U_D^\dagger(\tau) \psi_f$:
\begin{equation}
    L(z) = \int_0^\tau \|\dot{z}(t)\| \, \mathrm{d}t \le \tau \|A_{ND}\|_{op}.
\end{equation}

The length of this continuous path  bounds the geodesic 
Fubini-Study distance between its endpoints, yielding 
$d_{FS}(\psi_0, U_D^\dagger(\tau) \psi_f) \le \tau \|A_{ND}\|_{op}$. 
Recognizing that $U_D^\dagger(\tau) \in \mathcal{T}_B$ and taking the 
infimum over the entire subgroup $\mathcal{T}_B$ provides the lower 
bound.
\end{proof}

\subsection{Chordal relaxation}\label{SEC_Chordal}

\begin{proof}[Proof of Proposition \ref{PRO_Chordal_relaxation}]
By right-invariance of the system, we assume without loss of generality that $X_0=I_d$ (right-multiplying the solution by $X_0^{-1}$).
Let as before $\mathcal{T}_B$ be the set of matrices of $G$ commuting with $B$. By the same Lipschitz argument as Proposition \ref{PRO_Spectral_Bound}, applied to the propagator $z(t)=U_D^\dagger(t)X(t)$ (a matrix rather than a state vector), the transfer time is bounded by:
\[
\tau_{\mathrm{US}}(X_f)
\;\geq\;
\frac{
\inf_{D\in\mathcal{T}_B}\,\|X_0^{-1}X_f D - I\|_F
}{
\|A_{ND}\|_F
},\]
where the Frobenius norm of a matrix $X$ in $\mathfrak{g}$ is defined as 
$\|X\|_F^2 = \mathrm{Tr}(X^\dagger X)$.
For any $D \in \mathcal{T}_B$, let $U = X_0^\dagger X_f D$. Since $U \in G$, its eigenvalues are 
$e^{i\theta_j}$ with $\theta_j \in [-\pi, \pi]$.
The geodesic distance in $G$ with the bi-invariant Frobenius metric is bounded below by the chordal 
distance, since any $C^1$ curve $\gamma$ from $I$ to $U$ satisfies 
$\int_0^1\|\dot\gamma(t)\|_F\,dt\geq\|U-I\|_F$.
Expanding the squared Frobenius norm yields:
\begin{eqnarray}
    \| U - I \|_F^2 & =& \mathrm{Tr}\left( (I-U)^\dagger (I-U) \right) \\
	&= &     \mathrm{Tr}\left( I - U - U^\dagger + U^\dagger U \right) \\
	&=& \mathrm{Tr}\left(2 I - \mathrm{Re} (2U) \right) \\
    &=& 2d - 2 \, \mathrm{Re}\left( \mathrm{Tr}(U) \right).
\end{eqnarray}

Let $M = X_0^\dagger X_f$ and let $D = \mathrm{diag}(e^{i\phi_j})$. The trace is $\sum_{j=1}^d M_{jj} 
e^{i\phi_j}$. To obtain an analytic bound, we compare the optimization over $\mathcal{T}_B$ with that over the maximal torus of $U(d)$, denoted ${\mathcal{T}}_B^{U(d)}$. If $G=SU(d)$, the domain $\mathcal{T}_B$ imposes the global constraint $\sum \phi_j \equiv 0 \pmod{2\pi}$ and $\mathcal{T}_B \subset \mathcal{T}_B^{U(d)}$, so the supremum can only increase when passing to $\mathcal{T}_B^{U(d)}$; if $G=U(d)$, then $\mathcal{T}_B=\mathcal{T}_B^{U(d)}$. In either case the comparison maintains a valid lower bound:
\begin{equation}
    \sup_{D \in \mathcal{T}_B} \mathrm{Re}\left( \sum_{j=1}^d M_{jj} e^{i\phi_j} \right) \le \sup_{D \in 
    {\mathcal{T}}_B^{U(d)}} \mathrm{Re}\left( \sum_{j=1}^d M_{jj} e^{i\phi_j} \right).
\end{equation}
On $\mathcal{T}_B^{U(d)}$, the phases $\phi_j$ are independent. The 
supremum is achieved when each $e^{i\phi_j}$ cancels the 
complex argument of $M_{jj}$, yielding $\sum_{j=1}^d |M_{jj}|$. 
Substituting this analytic supremum back into the chordal inequality 
provides the result.
\end{proof}

\subsection{Algebraic Speed Limit via Ehrenfest theorem}\label{SEC_Ehrenfest}

\begin{proof}[Proof of Proposition \ref{PRO_minoration_Ehrenfest}]
Consider the expectation value of the observable $M(t) = \langle \psi(t), O \psi(t) \rangle$. Its time derivative is given by:
$$ \frac{d}{dt} M(t) = \langle \dot{\psi}(t), O \psi(t) \rangle + \langle \psi(t), O \dot{\psi}(t) \rangle. $$
Substituting the equation of motion $\dot{\psi} = (A + u(t)B)\psi$ yields:
$$ \frac{d}{dt} M(t) = \langle (A + u(t)B)\psi(t), O \psi(t) \rangle + \langle \psi(t), O(A + u(t)B)\psi(t) \rangle. $$
Using the definition of the adjoint and the anti-Hermiticity of the generators ($A^\dagger = -A$ and $B^\dagger = -B$), the first term transforms as $\langle X\psi, O\psi \rangle = \langle \psi, X^\dagger O \psi \rangle = \langle \psi, (-A - u(t)B)O \psi \rangle$. Factoring the inner product gives the sum of commutators:
$$ \frac{d}{dt} M(t) = \langle \psi(t), (OA - AO)\psi(t) \rangle + u(t) \langle \psi(t), (OB - BO)\psi(t) \rangle. $$
$$ \frac{d}{dt} M(t) = \langle \psi(t), [O, A]\psi(t) \rangle + u(t) \langle \psi(t), [O, B]\psi(t) \rangle. $$
By hypothesis, the observable is chosen such that $[O, B] = 0$. This  eliminates the dependence on the unbounded control field $u(t)$:
$$ \frac{d}{dt} M(t) = \langle \psi(t), [O, A]\psi(t) \rangle. $$
Since $A = A_D + A_{ND}$ with $A_D$ diagonal in the eigenbasis of $B$, and since $O$ commutes with $B$ (hence is diagonal in that basis, the spectrum of $B$ being simple), the commutator $[O, A_D]$ vanishes and $[O, A] = [O, A_{ND}]$.
Applying the Cauchy-Schwarz inequality and bounding the Rayleigh quotient by the operator norm, the rate of change is globally bounded by:
\begin{equation}\label{EQ_majoration_vitesse_Ehrenfest}
 \left | \frac{d}{dt} M(t) \right| \le \| [O, A_{ND}] \|_{op}. 
 \end{equation}
These computations were independent of the control, the source and the target. We now 
consider $\psi_0$, $\psi_f$ in the unit sphere $\mathbb{S}^{2d-1}\subset \C^d$, and an 
associated sequence $u_n:[0,t_n]\to \R$ of integrable functions such that $t_n$ tends to $
\tau_{\mathrm{ASS}}(\psi_0,\psi_f)$  and $X^{u_n}_{t_n,0} \psi_0 $ tends to $\psi_f$ 
as $n$ tends to infinity

Integrating differential inequality \eqref{EQ_majoration_vitesse_Ehrenfest} over the control duration $t \in [0, t_n]$ yields:
$$ \left| \langle X^{u_n}_{t_n,0}\psi_0,\, O\, X^{u_n}_{t_n,0}\psi_0 \rangle - \langle \psi_0, O \psi_0 \rangle \right| \le \int_0^{t_n} \| [O, A_{ND}] \|_{op}\, \mathrm{d}t = t_n \| [O, A_{ND}] \|_{op}.$$
Isolating $t_n$ yields 
$$
t_n \geq \frac{ \left| \langle X^{u_n}_{t_n,0}\psi_0,\, O\, X^{u_n}_{t_n,0}\psi_0 \rangle - \langle \psi_0, O \psi_0 \rangle \right|} {\| [O, A_{ND}] \|_{op}}$$
 and the result follows for $\tau_{\mathrm{ASS}}(\psi_0,\psi_f)$ letting $n$ tend to infinity. The result for $\tau_{\mathrm{SS}}(\psi_0,\psi_f)$ ensues trivially. 
\end{proof}

\section{Spectral Properties of the Shift Matrix}\label{SEC_Appendix}
\label{SEC_Tridiagonal}

Let $T$ be the $p \times p$ symmetric tridiagonal matrix defined by zero diagonal entries and off-diagonal entries equal to $1/2$:
\begin{equation}
    T = \frac{1}{2} \begin{pmatrix}
    0 & 1 & 0 & \cdots & 0 \\
    1 & 0 & 1 & \ddots & \vdots \\
    0 & 1 & \ddots & \ddots & 0 \\
    \vdots & \ddots & \ddots & 0 & 1 \\
    0 & \cdots & 0 & 1 & 0
    \end{pmatrix}
\end{equation}

\begin{lemma}\label{LEM_eigenvectors_shift_matrix}
The eigenvalues of $T$ are given by
$ \displaystyle 
    \lambda_m = \cos\left( \frac{m\pi}{p+1} \right)$, for $ m = 1, \ldots, p,
$
with corresponding unnormalized eigenvectors $x^{(m)} \in \mathbb{R}^p$ whose $k$-th component is:
\begin{equation}
    x_k^{(m)} = \sin\left( \frac{km\pi}{p+1} \right), \quad k = 1, \ldots, p.
\end{equation}
Consequently, the spectral radius of $T$ is $\|T\|_2 = \sup_{1\leq m \leq p} |\lambda_m|= \cos\left(\frac{\pi}{p+1}\right)$.
\end{lemma}

\begin{proof}
Let $\lambda$ be an eigenvalue and $x = (x_1, \ldots, x_p)^T$ be an associated eigenvector. The eigenvalue equation $Tx = \lambda x$ translates to the system of linear difference equations:
\begin{equation}\label{EQ_App_Difference}
    x_{k-1} + x_{k+1} = 2\lambda x_k, \quad k = 1, \ldots, p
\end{equation}
subject to the Dirichlet boundary conditions $x_0 = 0$ and $x_{p+1} = 0$.

This is a second-order linear homogeneous recurrence relation. Its characteristic equation is $r^2 - 2\lambda r + 1 = 0$. To satisfy the boundary conditions non-trivially, the roots must be complex conjugates on the unit circle. We parametrize the eigenvalue as $\lambda = \cos(\theta)$ for some $\theta \in (0, \pi)$, yielding the roots $r_{\pm} = e^{\pm i\theta}$.

The general solution to \eqref{EQ_App_Difference} takes the form:
\begin{equation}
    x_k = A \cos(k\theta) + B \sin(k\theta).
\end{equation}
The boundary condition $x_0 = 0$ forces $A = 0$, reducing the solution to $x_k = B \sin(k\theta)$. 
The second boundary condition $x_{p+1} = 0$ requires $\sin((p+1)\theta) = 0$. Since $x$ must be non-zero, $B \neq 0$, which implies:
\begin{equation}
    (p+1)\theta = m\pi, \quad m \in \mathbb{Z}.
\end{equation}
Restricting to $\theta \in (0, \pi)$ to obtain distinct eigenvalues yields $p$ roots:
\begin{equation}
    \theta_m = \frac{m\pi}{p+1}, \quad m = 1, \dots, p.
\end{equation}
Substituting $\theta_m$ back into the parametrization gives the spectrum $\lambda_m = \cos\left(\frac{m\pi}{p+1}\right)$ and the corresponding eigenvector components $x_k^{(m)} = \sin\left(\frac{km\pi}{p+1}\right)$. 
The spectral norm corresponds to the maximum absolute eigenvalue, achieved for $m=1$ (and $m=p$), proving $\|T\|_2 = \cos\left(\frac{\pi}{p+1}\right)$.
\end{proof}

\section*{Acknowledgments}
The author gratefully acknowledges the financial support for this work.  
It was funded, in whole or in part, by the French National Research Agency (ANR) under the QuBiCCS project (grant ANR‑24‑CE40‑3008‑01).  The research is also part of the EIPHI Graduate School initiative (grant ANR‑17‑EURE‑0002) and received additional backing from the Bourgogne Franche‑Comté Region.

It is a pleasure for the author  to thank Christiane Koch and Fernando Gago Encinas for proposing the
research question and for the stimulating discussions that helped shape the study.

The author also acknowledges the use of large language model services in the preparation of this manuscript. The commercial Google Gemini Pro 3.1 model was used to translate the original French notes, to draw the TikZ figures, and to prepare first drafts. The open-weight Qwen 3.8 model, self-hosted at the University of Strasbourg (Unistra), was used extensively for proofreading, for correcting and improving the mathematical proofs, and for refining the language. While these tools assisted with the writing and visualization, the author takes full responsibility for the scientific integrity, mathematical correctness, and final verification of the article.

\bibliographystyle{siam}
\bibliography{Biblio}

\begin{thebibliography}{10}

\bibitem{zbMATH05117619}
{\sc A.~Agrachev and T.~Chambrion}, {\em An estimation of the controllability
  time for single-input systems on compact {Lie} groups}, ESAIM, Control Optim.
  Calc. Var., 12 (2006), pp.~409--441.

\bibitem{Agrachev2004}
{\sc A.~A. Agrachev and Y.~L. Sachkov}, {\em Attainable Sets}, Springer Berlin
  Heidelberg, Berlin, Heidelberg, 2004, pp.~109--119.

\bibitem{ALBERTINI2002213}
{\sc F.~Albertini and D.~D'Alessandro}, {\em The lie algebra structure and
  controllability of spin systems}, Linear Algebra and its Applications, 350
  (2002), pp.~213--235.

\bibitem{Bang-Jensen2002}
{\sc J.~Bang-Jensen and G.~Gutin}, {\em Digraphs: Theory, Algorithms and
  Applications}, Springer London, London, 2002.

\bibitem{BasilewitschPRR2024}
{\sc D.~Basilewitsch, C.~Dlaska, and W.~Lechner}, {\em Comparing planar quantum
  computing platforms at the quantum speed limit}, Phys. Rev. Res., 6 (2024),
  p.~023026.

\bibitem{zbMATH05001889}
{\sc U.~Boscain, T.~Chambrion, and G.~Charlot}, {\em Nonisotropic 3-level
  quantum systems: complete solutions for minimum time and minimum energy},
  Discrete Contin. Dyn. Syst., Ser. B, 5 (2005), pp.~957--990.

\bibitem{zbMATH02122172}
{\sc U.~Boscain, G.~Charlot, J.-P. Gauthier, S.~Gu{\'e}rin, and H.-R. Jauslin},
  {\em Optimal control in laser-induced population transfer for two- and
  three-level quantum systems.}, J. Math. Phys., 43 (2002), pp.~2107--2132.

\bibitem{zbMATH05154205}
{\sc U.~Boscain and P.~Mason}, {\em Time minimal trajectories for a spin 1/2
  particle in a magnetic field}, J. Math. Phys., 47 (2006), pp.~062101, 29.

\bibitem{BOUSSAID2020108412}
{\sc N.~Boussaïd, M.~Caponigro, and T.~Chambrion}, {\em Regular propagators of
  bilinear quantum systems}, Journal of Functional Analysis, 278 (2020),
  p.~108412.

\bibitem{doi:10.1137/0125025}
{\sc R.~W. Brockett}, {\em Lie theory and control systems defined on spheres},
  SIAM Journal on Applied Mathematics, 25 (1973), pp.~213--225.

\bibitem{CanevaPRL2009}
{\sc T.~Caneva, M.~Murphy, T.~Calarco, R.~Fazio, S.~Montangero, V.~Giovannetti,
  and G.~E. Santoro}, {\em Optimal control at the quantum speed limit}, Phys.
  Rev. Lett., 103 (2009), p.~240501.

\bibitem{PhysRevLett.96.060503}
{\sc A.~Carlini, A.~Hosoya, T.~Koike, and Y.~Okudaira}, {\em Time-optimal
  quantum evolution}, Phys. Rev. Lett., 96 (2006), p.~060503.

\bibitem{zbMATH05509716}
{\sc T.~Chambrion, P.~Mason, M.~Sigalotti, and U.~Boscain}, {\em
  Controllability of the discrete-spectrum {Schr{\"o}dinger} equation driven by
  an external field}, Ann. Inst. Henri Poincar{\'e}, Anal. Non Lin{\'e}aire, 26
  (2009), pp.~329--349.

\bibitem{Clements:16}
{\sc W.~R. Clements, P.~C. Humphreys, B.~J. Metcalf, W.~S. Kolthammer, and
  I.~A. Walmsley}, {\em Optimal design for universal multiport
  interferometers}, Optica, 3 (2016), pp.~1460--1465.

\bibitem{Conway1976}
{\sc J.~Conway and A.~Jones}, {\em Trigonometric diophantine equations (on
  vanishing sums of roots of unity)}, Acta Arithmetica, 30 (1976),
  pp.~229--240.

\bibitem{Dumitrescu2017}
{\sc B.~Dumitrescu}, {\em Positive Trigonometric Polynomials and Signal
  Processing Applications}, Springer International Publishing, 2017.

\bibitem{zbMATH03515496}
{\sc R.~C. Entringer, D.~E. Jackson, and D.~A. Snyder}, {\em Distance in
  graphs}, Czech. Math. J., 26 (1976), pp.~283--296.

\bibitem{FarmanianPRA2024}
{\sc A.~Farmanian and V.~Karimipour}, {\em Quantum speed limits for
  implementation of unitary transformations}, Phys. Rev. A, 110 (2024),
  p.~042611.

\bibitem{fejer1915uber}
{\sc L.~Fej{\'e}r}, {\em {\"U}ber trigonometrische {P}olynome}, Journal f{\"u}r
  die reine und angewandte Mathematik, 146 (1915), pp.~53--82.

\bibitem{PhysRevA.92.062106}
{\sc M.~Gajdacz, K.~K. Das, J.~Arlt, J.~F. Sherson, and T.~c.~v. Opatrn\'y},
  {\em Time-limited optimal dynamics beyond the quantum speed limit}, Phys.
  Rev. A, 92 (2015), p.~062106.

\bibitem{zbMATH07308540}
{\sc J.-P. Gauthier and F.~Rossi}, {\em A universal gap for non-spin quantum
  control systems}, Proc. Am. Math. Soc., 149 (2021), pp.~1203--1214.

\bibitem{GoerzJPB11}
{\sc M.~H. Goerz, T.~Calarco, and C.~P. Koch}, {\em The quantum speed limit of
  optimal controlled phasegates for trapped neutral atoms}, JPB, 44 (2011),
  p.~154011.
\newblock arXiv:1103.6050.

\bibitem{GoerzNPJQI17}
{\sc M.~H. Goerz, F.~Motzoi, K.~B. Whaley, and C.~P. Koch}, {\em Charting the
  circuit-qed design landscape using optimal control theory}, npj Quantum Inf.,
  3 (2017), p.~37.
\newblock arXiv:1606.08825.

\bibitem{zbMATH07755548}
{\sc C.~Gorodski, C.~Lange, A.~Lytchak, and R.~A.~E. Mendes}, {\em A diameter
  gap for quotients of the unit sphere}, J. Eur. Math. Soc. (JEMS), 25 (2023),
  pp.~3767--3793.

\bibitem{GuentherAVSQS2023}
{\sc S.~Günther and N.~A. Petersson}, {\em A practical approach to determine
  minimal quantum gate durations using amplitude-bounded quantum controls}, AVS
  Quantum Science, 5 (2023), p.~043802.

\bibitem{HegerfeldtPRL13}
{\sc G.~C. Hegerfeldt}, {\em Driving at the quantum speed limit: Optimal
  control of a two-level system}, Phys. Rev. Lett., 111 (2013), p.~260501.

\bibitem{JURDJEVIC1972313}
{\sc V.~Jurdjevic and H.~J. Sussmann}, {\em Control systems on lie groups},
  Journal of Differential Equations, 12 (1972), pp.~313--329.

\bibitem{Kato1953}
{\sc T.~Kato}, {\em Integration of the equation of evolution in a banach
  space}, Journal of the Mathematical Society of Japan, 5 (1953), pp.~208--234.

\bibitem{PhysRevA.65.032301}
{\sc N.~Khaneja, S.~J. Glaser, and R.~Brockett}, {\em Sub-riemannian geometry
  and time optimal control of three spin systems: Quantum gates and coherence
  transfer}, Phys. Rev. A, 65 (2002), p.~032301.

\bibitem{khaneja2005optimal}
{\sc N.~Khaneja, T.~Reiss, C.~Kehlet, T.~Schulte-Herbr{\"u}ggen, and S.~J.
  Glaser}, {\em Optimal control of coupled spin dynamics: design of nmr pulse
  sequences by gradient ascent algorithms}, Journal of magnetic resonance, 172
  (2005), pp.~296--305.

\bibitem{Naseri2024}
{\sc M.~Naseri, C.~Macchiavello, D.~Bruß, P.~Horodecki, and A.~Streltsov},
  {\em Quantum speed limits for change of basis}, New Journal of Physics, 26
  (2024), p.~023052.

\bibitem{Nielsen_Chuang_2010}
{\sc M.~A. Nielsen and I.~L. Chuang}, {\em Quantum Computation and Quantum
  Information: 10th Anniversary Edition}, Cambridge University Press, 2010.

\bibitem{PhysRevLett.73.58}
{\sc M.~Reck, A.~Zeilinger, H.~J. Bernstein, and P.~Bertani}, {\em Experimental
  realization of any discrete unitary operator}, Phys. Rev. Lett., 73 (1994),
  pp.~58--61.

\bibitem{riesz1915uber}
{\sc F.~Riesz}, {\em {\"U}ber ein problem des {H}errn {C}arath{\'e}odory},
  Journal f{\"u}r die reine und angewandte Mathematik, 146 (1915), pp.~83--87.

\bibitem{Sachkov}
{\sc Y.~Sachkov}, {\em Controllability of invariant systems on lie groups and
  homogeneous spaces}, Journal of Mathematical Sciences, 100 (2000),
  p.~2355–2427.

\bibitem{PhysRevA.58.2733}
{\sc L.~Viola and S.~Lloyd}, {\em Dynamical suppression of decoherence in
  two-state quantum systems}, Phys. Rev. A, 58 (1998), pp.~2733--2744.

\bibitem{Williams2011}
{\sc C.~P. Williams}, {\em Quantum Gates}, Springer London, London, 2011,
  pp.~51--122.

\bibitem{zbMATH03328106}
{\sc B.~Zelinka}, {\em Medians and peripherians of trees}, Arch. Math. (Brno),
  4 (1968), pp.~87--95.

\end{thebibliography}
\end{document}